\RequirePackage{fix-cm}
\documentclass[smallextended]{svjour3}       
\smartqed  
\usepackage{graphicx}
\usepackage{amsmath}
\usepackage{amsfonts}
\usepackage{amssymb}
\usepackage{epsfig}

\usepackage{amsmath}
\usepackage{amsfonts}
\usepackage{amssymb}
\usepackage{dsfont}
\usepackage{mathrsfs}
\usepackage{bbm}
\usepackage{amsmath}
\usepackage{amsfonts}
\usepackage{amssymb}
\usepackage{multirow}
\usepackage{mathrsfs}
\usepackage[dvipsnames,usenames]{xcolor}

\newenvironment{pack_item}{
\begin{itemize}
  \setlength{\itemsep}{1pt}
  \setlength{\parskip}{0pt}
  \setlength{\parsep}{0pt}}{\end{itemize}
}

\newcommand{\bE}{\mathbb{E}}

\newcommand{\bR}{\mathbb{R}}

\newcommand{\cF}{\mathcal{F}}

\newcommand{\cP}{\mathcal{P}}

\newcommand{\fB}{\mathfrak{B}}
\newcommand{\fC}{\mathfrak{C}}

\newcommand{\fF}{\mathfrak{F}}

\newcommand{\fI}{\mathfrak{I}}

\newcommand{\fM}{\mathfrak{M}}

\newcommand{\fP}{\mathcal{P}}

\newcommand{\fW}{\mathfrak{W}}
\newcommand{\fX}{\mathfrak{X}}

\newcommand{\fZ}{\mathfrak{Z}}

\newcommand{\step}[1]{\mathbf{1}\left\{#1\right\}}
\newcommand{\prob}[1]{\text{Pr}\left\{#1\right\}}

\newcommand{\dd}{\mathsf {d\kern -0.07em l}} 

\newcommand{\bgeqn}{\begin{eqnarray}}
	\newcommand{\edeqn}{\end{eqnarray}}
\newcommand{\bgeq}{\begin{eqnarray*}}
	\newcommand{\edeq}{\end{eqnarray*}}
\newcommand{\inmat}[1]{\mbox{\rm {#1}}}
\def\bbe{{\bE}}

\newcommand{\vt}{{\vartheta}}

\usepackage{enumitem}

\usepackage{amsmath}
\allowdisplaybreaks[4]
\usepackage{verbatim}

\usepackage{wrapfig}
\usepackage{rotating}

\def\bbe{{\mathbb{E}}}
\renewcommand\thesubsubsection{\Alph{section}}
\newcommand{\R}{{\rm I\!R}}
\newtheorem{assumption}{Assumption} 

\def\prob{{\rm{Prob}}}

\newcommand{\tabincell}[2]{\begin{tabular}{@{}#1@{}}#2\end{tabular}}
\renewcommand\thesubsubsection{\arabic{section}.\arabic{subsection}.\arabic{subsubsection}}

\usepackage{float}
\usepackage[title]{appendix}
\usepackage{makecell} 
\usepackage{boldline}
\usepackage{diagbox}
\usepackage{threeparttable}
\usepackage[dvipsnames]{xcolor}
\usepackage{adjustbox}
\usepackage{booktabs}

\usepackage[colorlinks=true]{hyperref}
 \hypersetup{urlcolor=blue,
 citecolor=blue,linkcolor=blue}
   
\makeatletter
\g@addto@macro{\UrlBreaks}{\UrlOrds}
\makeatother

\begin{document}

\title{
Distributional 
Utility Preference Robust Optimization Models in 
Multi-Attribute Decision Making
\thanks{This project is supported by RGC grant 14500620 and CUHK Start-up grant.}
}

\titlerunning{Multi-Attribute Utility Preference Robust Optimization}    

\author{Jian Hu       \and
        Dali Zhang    \and
        Huifu Xu      \and 
        Sainan Zhang
}


\institute{
Jian Hu  \at
Department of Industrial and Manufacturing Systems Engineering, University of Michigan - Dearborn, Dearborn, MI 48128.
\\
\email{jianhu@umich.edu}           
\and
Dali Zhang \at
Sino-US Global Logistics Institute, Antai College of Economics and Management, Shanghai Jiao Tong University, Shanghai 200030, China. \\
\email{zhangdl@sjtu.edu.cn}  
\and
Huifu Xu, 
Sainan Zhang \at
Department of Systems Engineering \& Engineering Management,
			The Chinese University of Hong Kong, Shatin, N.T., Hong Kong. \\
\email{hfxu@se.cuhk.edu.hk}, 
\email{snzhang@se.cuhk.edu.hk}
}

\date{Received: date / Accepted: date}

\maketitle

\begin{abstract}
	Utility preference robust optimization (PRO) has recently been proposed to deal with optimal decision-making problems where the decision maker's 
 (DM's) preference over gains and losses is ambiguous. In this paper, we take a step further to investigate the case that the DM's preference is random.
 We propose 
to use a random utility function to describe 
the DM's preference and develop 
distributional
utility preference robust optimization (DUPRO) 
models
when 
the distribution of the random utility function
 is ambiguous.
We concentrate on 
 data-driven problems where 
samples of the random parameters are obtainable but 
 the sample size  
 may be relatively small.    
 In the case when the random utility functions are of piecewise linear structure, 
 we propose 
  a bootstrap method to construct the ambiguity set and 
 demonstrate how the resulting DUPRO can be solved 
 by 
  a mixed-integer linear program. 
The piecewise linear structure 
 is versatile in its
ability to incorporate classical non-parametric utility assessment methods into the sample generation of a random utility function.
Next, we expand 
the proposed DUPRO models and computational schemes to address general cases 
where the random utility functions are not necessarily piecewise linear. 
We 
show how the DUPRO models with piecewise linear random utility functions can 
serve as approximations for 
the DUPRO models with general random utility functions and allow us to quantify the approximation errors.
 Finally, we carry out some performance studies of the proposed bootstrap-based DUPRO model and report the preliminary numerical test results.
The paper is the first attempt to use distributionally robust optimization methods for PRO problems.
\end{abstract}




\maketitle

	\section{Introduction} 
	\label{sec:introduction}

In decision-making under uncertainty,
a utility representation characterizes the decision maker's (DM's) risk attitude towards the risk arising from systemic randomness (\cite{von_neumann_theory_1947,sugden_alternatives_1998}). 
It is traditionally assumed that the DM makes consistent choices and hence the utility function representing 
the DM's preference is deterministic, i.e.,  
the DM always makes the same decision in identical choice situations unless 
the DM is exactly indifferent between alternatives \cite{blavatskyy_stochastic_2007}. However, these assumptions are persistently violated in practice \cite{allais_comportement_1953,simon_rational_1956,tversky_intransitivity_1969,tversky_framing_1981}. Indeed, it is usually observed that the DM exhibits mutable and ambivalent preference, particularly without 
complete information of 
complex problems and uncertain environments
in the real world \cite{camerer_experimental_1989,wu_empirical_1994,hey_investigating_1994,starmer_developments_2000}.
Ambivalent preference
may also
be caused by other reasons, e.g, 
the DM's preference does not satisfy transitivity in multi-attribute decision-making; 
the independence axiomatic property in the von Neumann-Morgenstern's expected utility theory fails to hold 
\cite{allais_comportement_1953}; 
decision-making is state dependent \cite{anscombe_definition_1963,karni1983state}
and 
 there are inconsistencies 
in
responses
in the preference elicitation process
 \cite{bertsimas2013learning}.
This prompts us to adopt
a random utility function 
to describe 
the DM's preference
where each scenario 
corresponds to the DM's particular preference.
This kind of random utility 
differs from the well-established 
random utility theory in discrete choice models where the latter 
is used to describe a group of customers' preferences
to a commodity (product) and the randomness   
characterizes variability of customers' preferences
and/or 
idiosyncratic product specific random shock, see e.g.~ \cite{fishburn_stochastic_1998,koppen_characterization_2001} and 
references therein.
With the random utility function, 
we may consider the expected (mean) value of 
the DM's utility from a modeller's perspective if we know the probability distribution of the random parameter. 
The mean utility value captures
the DM's average utility preference 
for the future decision-making.
We refer readers to \cite{fishburn_stochastic_1998} and references therein  
for a thorough discussion of random utility representations.
In the absence
of complete information of 
the probability distribution of the random parameter,
 we propose a distributionally robust model where the worst-case 
distribution from an ambiguity set is to be considered for the calculation of the mean of the random utility.
The next example explains the basic idea.

\begin{example} \label{ex:hannah} 
Hannah wants to buy a jacket and 
hesitates
to choose between $A$ and $B$ while considering price, style and color. 
Her wavering 
preference can be reasonably attributed to the assumption of
random utility. 
Denote by $u_1$ and $u_2$ respectively the utility functions which characterize Hannah's two possible inconsistent preferences over price, style
and color. 
Suppose that $u_1(A) = 0.9$, $u_2(A) = 0.4$, and $u_1(B) = u_2(B) = 0.5$.
Hannah's hesitance is explained mathematically 
by $u_1(A) > u_1(B)$ but $u_2(A) < u_2(B)$. 
This vague preference can be represented as a random utility function, which is $u_1$ with probability $p$ and $u_2$ with probability $1 - p$. Assume that, based on Hannah's past clothing shopping choices, $u_1$ could be her most likely taste and thus $p \ge 0.5$. A natural way is to use the mean utility function, $\bE_p[u] = p u_1 + (1-p) u_2$, to reconcile the inconsistency.
When the exact $p$ is unknown, 
a distributionally robust model would suggest or predict Hannah to choose $A$ since
$$
\min_{p \in [0.5, 1]} \bE_p [u(A)] = 0.65 >  \min_{p \in [0.5, 1]} \bE_p [u(B)] = 0.5.
$$
Without any prior knowledge of the probability distribution, it leads to a deterministic worse-case utility approach over the entire sample space $\{u_1, u_2\}$ as follows:  
$$
\max_{x \in \{ A, B \}} \min_{p \in [0, 1]} \bE_p [u(x)] = \max_{x \in \{ A, B \}} \min \{u_1(x), u_2(x)\}
 = \max\{0.4,0.5\}=0.5,
$$
which 
means that 
Hannah 
should or would go with $B$. 
In our view, 
the mean utility approach is more reasonable, given the information on Hannah's past shopping choices. 
\end{example}


Research on decision-making based on the worst-case utility may  be traced back to 
earlier work 
by Maccheroni \cite{maccheroni2002maxmin} 
who considers a worst-case expected utility model for 
 	a conservative DM with an unclear evaluation of the different outcomes when facing lotteries. Hu and Mehrotra \cite{hu_robust_2012} (also see \cite{hu_robust_2015}) look into the issue 
	from robust optimization perspective which is 
	later on known 
as preference robust optimization (PRO). They propose a moment-type framework
		for constructing an ambiguity set of a DM's utility preference which covers a number of important preference elicitation approaches including certainty equivalent and pairwise comparison.
To solve the resulting PRO model, they develop  
a step-like approximation scheme for the functions 
in the moment conditions and carry out some convergence analysis
for the justification of the approximation.

Armbruster and Delage \cite{armbruster_decision_2015}  consider an ambiguity
	set of utility functions which meet some criteria such as preferring certain lotteries over other lotteries and being risk averse, $S$-shaped, or prudent. 
	Instead of trying to identify a single utility function satisfying the criteria, they develop  a maximin  PRO model 
	where the optimal decision  is based on the worst
	utility function from the ambiguity set and demonstrate how the maximin optimization can be reformulated as a finite-dimensional linear program.
Guo et al.~\cite{guo2023utility}
take it further by proposing  a
piecewise linear approximation of the 
utility function
and 
 quantifying the 
 approximation error and its propagation 
 to the optimal value and optimal solutions of the PRO model.
 	In the case that a DM
  has a nominal utility 
 	function but lacks complete information about 
 	whether the nominal utility function is the true one,
	Hu and Stepanyan \cite{hu_optimization_2017} propose 
	to construct an ambiguity set of 
	utility functions 
	neighbouring the nominal.

Over the past few years, research on the PRO models has expanded in several directions.
Haskell et al.~\cite{haskell2016ambiguity} 
propose a robust optimization model where
the optimal decision is based on the worst-case 
utility function and the worst-case probability 
distribution of the underlying exogenous uncertainty amid both the true utility and the true probability distribution are ambiguous. 
Haskell et al.~\cite{HaskellHuangXu2018,HaskellXuHuang2022} 
extend the PRO model from single-attribute
decision-making 
to multi-attribute decision-making and from 
utility preference robust models 
to
a broader class of
preference robust quasi-concave choice function models.  
The latter gets around  the Allais paradox  
concerning the utility-based PRO models.  
Liu et al.~\cite{liu2021multistage}
present a multistage utility preference robust optimization model and 
demonstrate how the
state-dependent ambiguity set of utility functions 
may be constructed and how the state dependence may affect 
time-consistency of the dynamic model. 

The PRO models are also effectively extended to risk management problems
from convex risk measures by Delage and Li \cite{delage2017minimizing,delage2022shortfall}
to spectral coherent risk measures by Wang and Xu \cite{wang2020robust}.
These works deal with the issue that a DM
may have 
several risk measures at hand but lack complete information 
as to which one should be used for his/her decision-making. 
The authors propose various approaches
for constructing the ambiguity set and demonstrate how  
to solve the resulting PRO models efficiently. 
More recently,
Li~\cite{li2021inverse} proposes an inverse optimization approach
for solving 
PRO
models when the 
DM's
risk preference 
is obtained in a learning process.
Most of these PRO models 
consider the worst-case preference.
For decision-making with ambiguous utility functions,
minimax regret is an alternative criterion.
In the multi-attribute 
utility case,
 Boutilier et al.~\cite{boutilier2006constraint} seek decisions that achieve minimum worst-case regret, that is, regret experienced by a posteriori
 once the true utility function is revealed.
Vayanos et al.~\cite{VMYDR20} consider the worst-case absolute regret to mitigate the conservatism of classical robust optimization (which maximizes worst-case utility).

Our approach fundamentally differs from the existing literature research. We introduce the distributional utility PRO (DUPRO) model, which is motivated by the random utility theory. In contrast, those existing PRO models are grounded in von Neumann–Morgenstern's utility theorem \cite{von_neumann_theory_1947}. This distinction in motivation has far-reaching implications for our understanding of preference robustness. Von Neumann-Morgenstern’s utility theory, along with its underlying PRO models, is not suitable for the problems like Example \ref{ex:hannah} where a  deterministic utility function does not exist 
for characterizing the DM's inconsistent preferences. 
By contrast, the random utility 
effectively addresses the inconsistency issue caused by the random preferences.
In Example \ref{ex:hannah}, Hannah's vague choice between jackets $A$ and $B$ can be described as $\text{Prob}(u(A) > u(B)) = p$
and $\text{Prob}(u(A) < u(B)) = 1 - p$. However, this inconsistent pairwise comparison, used in the PRO models, results in the emptiness of the ambiguity set of utility functions. 
As such, our DUPRO model 
complements the existing research on PRO models.
The main contributions of this paper are summarized as follows.
	\begin{pack_item}
 		\item We investigate the case that the DM's preference is not only ambiguous but also potentially inconsistent or even displaying some kind of randomness. 
	Differing from the classical worst-case utility maximization approaches, 
	we propose a  
 \ref{DUPRO}
 framework where the DM's preference is represented by a random utility function and the ambiguity is described by a set of probability distributions of the random utility.
An obvious advantage is that the \ref{DUPRO} model
can 
accommodate 
inconsistencies and randomness 
encountered 
in the preference elicitation process.

		\item In practice, a utility function is often elicited with a piecewise linear structure. In this case
		we show how the randomness can be represented by the increments of the
		each linear piece and subsequently 
		propose two statistical approaches which suit  
		data-driven problems
		for constructing an ambiguity set of 
		the distributions of random increments.
		One is to
		 construct an ellipsoidal  confidence region 
		with sample mean and sample covariance matrix which is widely used in the literature of distributionally robust optimization (DRO,  see e.g. \cite{delage_distributionally_2010});
		the other is to specify a nonparametric percentile-$t$ bootstrap confidence region.
While the bootstrap approach is widely used in statistics, 
we have not seen the approach 
used to construct
an ambiguity set of probability distributions in  DRO models. 
Here we adopt it because 
in practice it is often difficult
to obtain samples of a DM's preference and hence
the sample size is relatively small in our model.
Moreover, compared to the ellipsoidal method, bootstrap has two potential advantages 
for constructing an ambiguity set. 
Unlike
the ellipsoidal method which requires an estimation of the size of the ellipsoid for a given confidence level which could be unusually conservative (very large when sample size is small), the bootstrap-based ambiguity set is uniquely determined by 
the specified confidence level and the set of resamples. 
Although the two concepts of
confidences differ slightly, they are close when the 
size of original sample is sufficiently large. 
Moreover, we demonstrate the DUPRO model based on the bootstrap ambiguity set can be reformulated as an MILP which can be solved fairly efficiently.

\vspace{0.3cm}

		\item 
 For the \ref{DUPRO} model with general random utility functions, we demonstrate how to use a piecewise linear random utility function to approximate it and derive 
error bound  which depends
 only on the largest distance between two consecutive grid points
 of the piecewise linear random utility function 
 for 
 the optimal value 
 and
 convergence of the optimal solutions
 (Proposition \ref{P-utility-appr-PL}).
With the piecewise linearly approximated 
 utility maximization problem in place, 
 we move on to discuss how to solve it with standard sample average approximation method in data-driven setting and derive exponential rate of convergence of the optimal values and almost sure convergence of the optimal solutions as the sample size increases (Proposition \ref{C-u-appro-PL-DRO}).
Moreover, we consider DUPRO approaches 
for the piecewise linearly approximated model and derive exponential rate of convergence of the optimal values 
when the ambiguity set is constructed with ellipsoid method 
(Theorem \ref{thm:DUPRO-Ellip-optim})
and almost sure convergence when the ambiguity set is constructed with the bootstrap method (Theorem \ref{thm:DUPRO-Bootstrap-convg}).
  We have also conducted some numerical tests
  on the proposed \ref{DUPRO} model 
 and computational schemes 
 via both an 
 academic example and a practically oriented 
 case study. 
 The test results show that \ref{DUPRO} performs well as expected in the theoretical analysis 
  in terms of asymptotic convergence, computational time, and out-of-sample performance. 

\vspace{0.3cm}

	\end{pack_item}
	
	The rest of the paper 
	are
	organized as follows. Section~\ref{sec:Random_additive_PLUF} structures random multi-attribute utility functions with a piecewise linear additive 
 structure. Section~\ref{sec:Data-Driven DUPRO} proposes two DUPRO Models.
The tractable reformulation and solution method of the bootstrap model are developed in Section~\ref{sec:Reformulations and Solution Methods}. 
Section \ref{sec:PLA_random utility} extends 
 the DUPRO models
 to 
 general random utility function cases
 without a piecewise linear structure. 
 Section~\ref{sec:Case Studies} 
 reports preliminary numerical test results.
		Section \ref{sec:Conclusions} concludes.
  
	\section{Random additive piecewise linear multi-attribute utility functions}
\label{sec:Random_additive_PLUF}	
\begin{figure}[!htbp]
 \vspace{-0.3cm}
\minipage{0.5\textwidth}
 \centering
  \includegraphics[width=6.5cm]
  {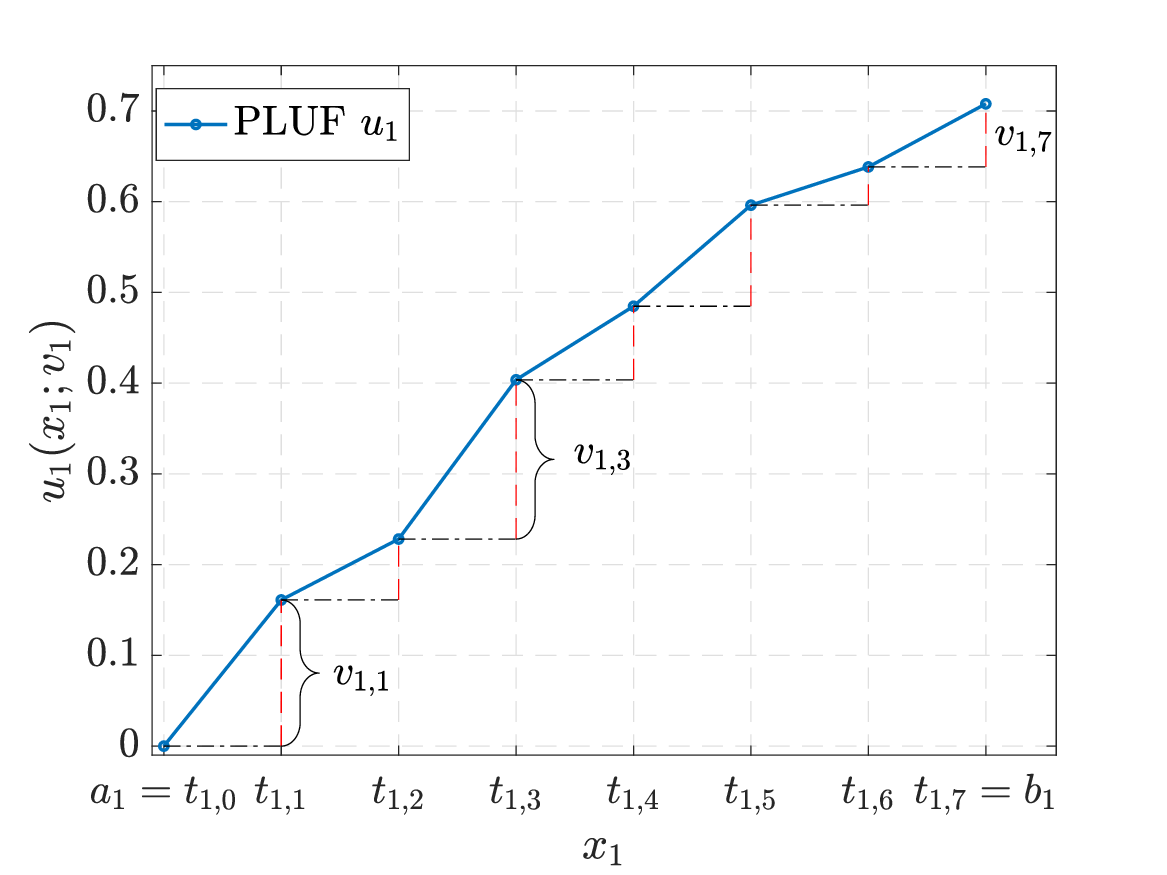}
  \text{ (a) $u_1(x_1;v_1)$ with $I_1=7$, $a_1=0,b_1=0.7$}
\endminipage\hfill
\minipage{0.5\textwidth}
  \centering
  \includegraphics[width=6.5cm]{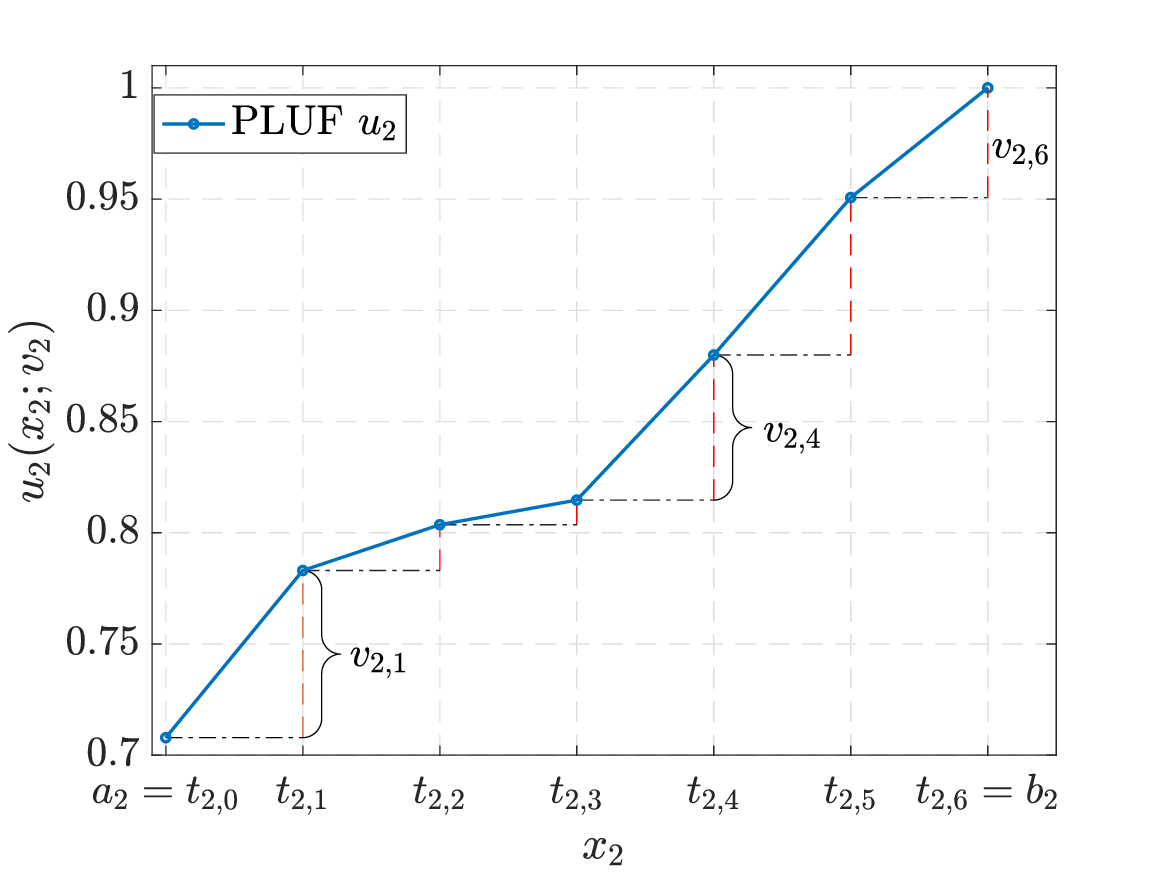}
  \text{ (b) $u_2(x_2;v_2)$ with $I_2=6$, $a_2=0,b_2=0.6$}
\endminipage
\hfill
\caption{\footnotesize PLUF $u(x;v)=u_1(x_1;v_1)+u_2(x_2;v_2)$ with $x=(x_1,x_2)^T$, $v=(v_1^T,v_2^T)^T$,
$v_1=(v_{1,1},v_{1,2},\cdots,v_{1,7})^T$,
$v_2=(v_{2,1},v_{2,2},\cdots,v_{2,6})^T$.}
\label{fig:example_M2}
\end{figure}

	We consider a multi-attribute decision-making problem where the DM's preference
is 
described by an additive
utility function.
Additive multi-attribute utility functions are widely used in the literature of 
decision analysis and behavioural economics, see e.g.~\cite{pollak1967additive,huber1974multi,torrance1995multi,ghaderi2017linear,kadzinski2013robust} and references therein.
Here we use it primarily because 
the tractable reformulations
and  the theoretical results 
of the DUPRO models to be developed in the forthcoming discussions  
rely heavily on the additive structure.
	Let $x_m$ be the performance of attribute $m$, for $m \in \fM :=\{1, \dots, M\}$, and
	$x := (x_1, \cdots, x_M)^T\in\mathbb{R}^M$ 
	be the vector of all attribute performances.
	Denote by $u_m:\mathbb{R}\to\mathbb{R}$ the 
single-attribute
 utility function of
 attribute $m$.
	We begin by considering a piecewise linear utility  function (PLUF) 
 of
 each attribute and then 
	move on to discuss general utility functions in Section \ref{sec:PLA_random utility}. 
 In this setup, we assume the DM's preference can be represented by a PLUF.
 	This is not only because the PLUF will facilitate us to derive tractable 
	DUPRO models but also because the DUPRO models with the PLUF can be conveniently applied in 
 real-world data-driven problems. Indeed, the PLUF structure is versatile in its ability to incorporate classical non-parametric utility assessment methods \cite{farquhar_utility_1984,wakker_eliciting_1996} into the random sampling processes required by the DUPRO models to generate observations of a random utility function.
	Let $u_m$ be a PLUF
	defined over interval $[a_m, b_m]$ with breakpoints
	\begin{align} \label{def:break_points}
		a_m = t_{m, 0} < \cdots < t_{m, I_m} = b_m,
	\end{align}
 where $t_{m,i}$ is a certain level of attribute $m$ and we assume that the DM's utility value at $t_{m,i}$ can be elicited. 
Assume without loss of generality that $u_m(a_m) = 0$ 
and let $v_{m,i} := u_m(t_{m,i}) - u_m(t_{m, i-1})$ be the 
	overall increment of $u_m$ over $[t_{m,i-1}, t_{m,i}]$  
	for $i \in \fI_m := \{1, \dots, I_m\}$ and
	$v_m := (v_{m,1}, \cdots v_{m, I_m})^T\in \mathbb{R}^{I_m}$.
	If we regard $v_m$ as a vector of parameters, then 
	we can 
	obtain a class of 
 piecewise linear functions parameterized by $v_m$ 
	for the $m$-th
	attribute as follows:
	\begin{align} \label{equ:u_m}
		u_m (x_m; v_m) := \sum_{i \in \fI_m} \left[
		\frac{v_{m,i}}{t_{m, i} - t_{m, i-1}} (x_m - t_{m, i-1}) + \sum_{j = 1}^{i-1} v_{m,j} \right] \step{t_{m, i-1} < x_m \le t_{m, i}},
	\end{align}
	where
	$\step{\cdot}$ is an indicator function. 
 In other words, $u_m (x_m; v_m)$
 defines a class of 
PLUFs
parameterized by a vector of increments ${v}_m$.
	Let $v :=(v^T_1, \cdots, v^T_M)^T$ and 
	$I := \sum_{m \in \fM} I_m$.
	We 
	can then 
	define 
the aggregate
 multi-attribute 
 PLUF
$u:\mathbb{R}^M\times \mathbb{R}^I\to \mathbb{R}$
	with an additive form as
	\begin{align} \label{MAPLU}
		  u(x; v) := \sum_{m \in \fM} u_m (x_m; v_m).
	\end{align}
	Without loss of generality, we assume that the utility function is nondecreasing and normalized to $[0,1]$ (i.e., $u(a;v)=0$ and $u(b;v)=1$ with $a=(a_1,\cdots,a_M)^T$, $b=(b_1,\cdots,b_M)^T$)
 because 
	the normalization does not affect its representation of the DM's preference.
	Under this assumption,  the utility function $u_m$
	is uniquely determined by vector $v_m$ for  $m \in \fM$
	for a given set of breakpoints\footnote{In the numerical tests, these points are generated randomly.}. 
Figure~\ref{fig:example_M2} depicts the
 PLUF
 when $M=2$.
 Moreover, the vector $v$ 
	lies in 
 the simplex of $\mathbb{R}^I$ with
	$v \ge 0$ and $e_I^T v = 1$, where $e_I$ is the $I$-dimensional vector of all ones. Note that a traditionally defined additive multi-attribute utility function is the weighted sum of the normalized 
single-attribute
 utility functions of
 each attribute. The normalization of $u$ and the definition of the utility function in \eqref{equ:u_m}-\eqref{MAPLU} ensure that the criterion weight of the $m$-th attribute is included but hidden in $u_m$.
	
	In practice, the vector $v$ of increments is
	determined 
 by elicited information on the DM's utility scoring at certain levels of the attributes.
 As we discussed in the Introduction, 
 we concentrate on the case that the DM's 
 preferences vary 
 randomly
 and they are
 potentially inconsistent,  
 which means that 
 the values of $v$ obtained as such 
 cannot be used to construct a single deterministic utility function, rather they may be treated as the realizations of a random vector $V$.
	This motivates us to replace
	the deterministic vector of parameters $v$ 
	with a random vector $V:(\Omega, \mathscr{F}, P)\to\mathbb{R}^I$
	in \eqref{MAPLU} and subsequently 
 consider a 
	random 
 PLUF
 mapping from $\mathbb{R}^M$ to $\mathbb{R}$ as
	\begin{align} \label{Random-MAPLU}
		u(x; V) = \sum_{m \in \fM} u_m (x_m; V_m),
	\end{align}
 where each realization of $V$
 may be interpreted as a ``state variable''
signifying the state of the DM (the DM's mood, situation, etc.). 
	This effectively randomizes 
 the class of 
utility functions
 parameterized by random vector $V$  
	with the support set
	\begin{align} \label{set:fV}
\mathcal{V}	:= \left \{
		v \in \bR^I \ |\
		e_I^T v =1, \
		v \ge 0
		\right \}.
	\end{align}
	In this setup, the DM's preference is described by 
a random 
PLUF
rather than a deterministic function. It means that the DM's utility may change randomly and is uncertain. The uncertainty may arise from the DM's inconsistent preferences or from incomplete and inaccurate observations of the DM's preference. 

	In some cases, we may require that a utility function 
	has
	properties on its shape. For example, the DM's risk-averse preference should be characterized as an increasing concave utility function. The support set in the case of increasing concave utility functions is described as
	\begin{align} \label{set:fV_C}
		\mathcal{V}_c:= \left \{
		v \in \bR^I
		\left |\;
		e_I^T v = 1, \
		A v \ge 0
		\right.
		\right \},
	\end{align}
	where matrix
	\begin{align*}
		A := \left[
		\begin{array}{ccc}
			A_1 & & \\
			& \ddots & \\
			& & A_M
		\end{array}
		\right]_{(I - M) \times I},
	\end{align*}
	in which blocks $A_m$ for $m \in \fM$ are
	\begin{align*}
		A_m := \left[
		\begin{array}{ccccc}
			\frac{1}{t_{m,1} - t_{m,0}} & \frac{-1}{t_{m,2} - t_{m,1}} & & & \\
			& \frac{1}{t_{m,2} - t_{m,1}} & \frac{-1}{t_{m,3} - t_{m,2}} & & \\
			&  & \ddots & \ddots & \\
			& & & \frac{1}{t_{m,I_m-1} - t_{m,I_m-2}} & \frac{-1}{t_{m,I_m} - t_{m,I_{m-1}}}
		\end{array}
		\right]_{(I_m-1) \times I_m}.
	\end{align*}

Before concluding this section, we note that there are a number of ways to represent a piecewise linear utility function. 
The specific approach that we adopt here offers significant advantages in deriving the DUPRO models discussed in the next two sections. As demonstrated in Proposition \ref{pp:equ_models}, the utility form $u(x; v)$ that we introduce in \eqref{MAPLU} can be expressed as a linear function of the increment vector $v$ for a given $x$. This characteristic facilitates the tractable formulations of the 
DUPRO models. 
Moreover, 
the specific form of the 
utility function differs 
from the conventional 
way of representing a multi-attribute utility function as a weighted sum of normalized single-attribute utility functions. When eliciting each normalized single-attribute utility function and assessing the trade-offs between them from the same DM, we posit that the criterion weights and these single-attribute utility functions are correlated random quantities. The random utility $u(x;V)$ described in \eqref{Random-MAPLU} captures this correlated relationship by integrating the criterion weights into the single-attribute utility functions.

\section{Data-driven distributional preference robust model}
	\label{sec:Data-Driven DUPRO}
	
	We consider a multi-attribute decision-making problem which aims to maximize the overall utility of attribute performances.
	Since the DM's utility function is stochastic, we consider the expected utility $\bE_P [u(x; V)]$ where
	the expectation is taken
  with respect to (w.r.t) 
 the probability distribution $P$ of $V$.
 Moreover, since the true probability distribution $P$ is unknown in some
	data-driven problems,  we may have to rely on
 partially available information such as empirical data,
 computational simulations,
 or
	subjective judgment to construct an ambiguity set of distributions.

	\subsection{The \ref{DUPRO} Models}
	Denote by ${\cal P}$
	the ambiguity set of the distributions of $V$ and
	consider the optimal decision based on the worst-case  
 distribution from the set as follows:
	\begin{align} \tag{DUPRO} \label{DUPRO}
		\max_{x \in \fX} \min_{P \in 
		{\cal P}
		} \bE_P [u(x; V)],
	\end{align}
	where $\fX \subseteq \otimes_{m \in \fM} [a_m, b_m] := [a_1, b_1] \times \dots \times [a_M, b_M]$ is the joint region of possible attribute values. 
	 There are some fundamental differences between \ref{DUPRO} and the existing PRO models in the literature. 
First, the true utility function
in \ref{DUPRO} is random whereas the true utility function in the existing PRO models 
is deterministic.
It means the former is about ambiguity of the probability distribution of the true random utility function whereas the latter is about 
ambiguity of the true unknown (deterministic) utility function.
Second, the random utility functions at different scenarios do not have to be consistent when they are applied to 
a given pair of prospects. In contrast, the existing PRO models are built on von Neumann-Morgenstern's expected utility theory \cite{von_neumann_theory_1947} which means 
the true utility function exists and is unique up to positive linear transformation although in practice it is unknown.
Thus the DM's utility preference must be consistent. Of course, in practical applications, 
the DM's answers to certain questionnaires may be inconsistent due to mistaken answers or 
errors in measurements and there is a way to handle this kind of inconsistency, see \cite{armbruster_decision_2015,bertsimas2013learning}.
However,
the two models are different in nature.
Third, in \ref{DUPRO}, 
we assume that 
a random sample of $V$ is obtainable explicitly
by elicitation on the DM's utility scoring at a specified set of levels of attributes. 
If we interpret  
$V=v$ as a scenario of $V$, it means  that   
the 
DM's utility
$u(x;V)$
at the scenario is obtainable, 
there is no ambiguity about 
 $u(\cdot;v)$. 
The only ambiguity lies in the true probability distribution of $V$.
Fourth, \ref{DUPRO} should be distinguished from distributionally robust expected utility optimization models where the ambiguity lies 
in the probability distribution of 
exogenous random parameters, see e.g.~\cite{natarajan2010tractable}.
Fifth, the random utility model has been widely used 
to describe customer's choices of certain products, see e.g.~\cite{berry1995automobile,cascetta2009random,aksoy2013price,sun2017saa}.
\cite{mishra2012choice} seems to be the first 
to consider a distributionally robust 
choice model where the expected utility of a selected product
is maximized under the ambiguity of correlations between 
customer's utility valuations of different products.
There are at least three differences between our models and theirs.
(a) 
In the discrete choice models,
the random utility function
represents a group
of 
DMs' (customers') utility evaluations
of certain products
where each realization of the random utility
 represents some customers' evaluation
of a product. In our models, 
the random utility function
represents a single DM's utility evaluations 
in a multi-attribute decision-making problem where
each realization of the random utility
 represents the DM's evaluation of the attributes in a particular ``state''.
(b) We consider the expected total additive utilities 
rather than the expected value of 
the maximum utility of the $M$ attributes.
(c) The ambiguity in our model is on randomness of the 
single-attribute
utility function of
each attribute
rather than correlation between these 
utility functions
(corresponding to correlation of utility of different products).

	In practice, $x$ is often a vector-valued function of some action denoted by $z$, that is, $x = h(z)$. Let $\fZ$ be a feasible action space. It follows that 
	\begin{align*}
		\fX = \{x \in \mathbb{R}^M
		\;|\ x= h(z),  \ z \in \fZ\}.
	\end{align*} 
	Consequently, we rewrite this variation 
 of \ref{DUPRO} as
	\begin{align} \tag{DUPRO-1}  \label{DUPRO-1}
		\max_{z \in \fZ} \min_{P \in 
		{\cal P}
		} \bE_P [u(h(z); V)]. 
	\end{align}
	In the case that the relation between $x$ and $z$ is affected by 
	some exogenous uncertainties denoted by a vector of random variables $\xi$
	with distribution $Q$, we may obtain a further variation 
	\begin{align} \tag{DUPRO-2} \label{DUPRO-2}
		\max_{z \in \fZ} \min_{P \in 
		{\cal P}
		} \bE_P [\bE_Q[u(h(z,\xi); V)]]. 
	\end{align}
	We next focus on discussing \ref{DUPRO}. All the model configurations and corresponding solution methods presented in this paper can be straightforwardly extended to \ref{DUPRO-1} and \ref{DUPRO-2}.

	\subsection{Construction of the ambiguity set}
	\label{subsec:Construction of ambiguity set}
	
	A key component of \ref{DUPRO} is the ambiguity set. 
	In the current main stream research 
	on
 DRO models,
the ambiguity is concerned with the probability distribution of exogenous uncertainty,
see \cite{goh2010distributionally,calafiore2006distributionally,wiesemann2014distributionally,hanasusanto2015distributionally,chen2018data,long2022robust,chen2019distributionally} and the references therein. Here, the ambiguity lies in the probability distribution of 
endogenous uncertainty (DM's utility preference).
	We 
	outline two main approaches for 
	constructing the ambiguity set ${\cal P}$
	in a data-driven environment. Denote by $\mu$ and $\Sigma$ the mean and covariance matrix of $V$. 
	Here we consider a situation where the true $\mu$ and $\Sigma$ are unknown 
	but it is possible to obtain an approximation with sample data.
	Let $V^1, \dots, V^N$ be an independent and identically distributed (i.i.d. for short) 
	random sample of $V$ and denote by $\bar{V}$ the sample mean and by $S$ a sample covariance matrix 
	approximating their true counterparts $\mu$ and $\Sigma$. 
 In practice, the random samples can be 
 generated using classical non-parametric utility assessment methods. In Section \ref{subsec:Project Investment}, we delve into the utilization of a generative logistic regression method to assemble a random dataset of consumer preference in automobile market.

Note that traditional PRO models are not suitable 
in this case. 
This is primarily because 
there may not exist a deterministic utility function which fits into 
the random samples. 
For instance, let $u(\cdot \; ; V^1), \dots, u(\cdot \; ; V^N)$ be the $N$ observations of  
Hannah's vague preference in Example \ref{ex:hannah}.
It is highly likely that one observation 
supports the choice of jacket $A$, i.e., $u(A; V^i) > u(B; V^i)$,
whereas another 
favors $B$ 
with $u(A; V^j) < u(B; V^j)$. 
Those inconsistent observations result in 
emptiness of the ambiguity set of utility functions in 
the existing PRO models \cite{armbruster_decision_2015,hu_robust_2015}.

	In what follows, we address how to construct a set of moment uncertainty using 
	a random sample. The support set $\mathcal{V}$
	defined as
	in \eqref{set:fV} indicates that the components of $V$ are linearly dependent, i.e, $e_I^T V = 1$ almost surely. Therefore, we can reduce by one dimension to consider the first $I-1$ components of $V$ in the construction of the ambiguity set. Accordingly, $V_{M, I_M}$, which is the $I$-th component of $V$, is 
	equal to $1$ less
	the sum of the other components. 
	Let 
	\begin{align}
		\label{eq:C=matrix}
		C := \left[\begin{array}{cccc}
			1 & & & 0 \\
			& \ddots &  & \vdots \\
			& & 1 & 0	
		\end{array}
		\right]_{(I-1) \times I}.
	\end{align}
	Then $C V\in \R^{I-1}$ is the vector consisting of the first $I-1$ components of $V$. Denote by $\mu_{I-1}$ and $\Sigma_{I-1}$ the mean and covariance matrix of $CV$ and by $\bar V_{I-1}$ and $S_{I-1}$ the sample counterparts of $\mu_{I-1}$ and $\Sigma_{I-1}$. Note that 
	$$\mu_{I-1} = C \mu, \ \Sigma_{I-1} = C \Sigma C^T, \ \bar V_{I-1} = C \bar V, \ \text{and} \  S_{I-1} = C S C^T.$$
	To facilitate our analysis in the forthcoming discussions, we make a blanket assumption on $\Sigma_{I-1}$.
	
	\begin{assumption}
		\label{Assu-Sigma}
		The covariance matrix $\Sigma_{I-1}$ is nonsingular.
	\end{assumption}

	Under this assumption, $S_{I-1}$ is nonsingular for a sufficiently large sample size $N$, which is needed in the following constructions of the ambiguity set ${\cal P}$.
	The assumption requires that the components of $CV$ are 
	linearly independent. If the assumption is not satisfied, then we will repeat the procedure of dimension reduction until remaining components are linearly independent. 

\subsubsection{Ellipsoid approach} 
A 
popular
approach 
in the literature of
distributionally robust optimization 
is to specify 
	the ambiguity set by moment conditions
	\cite{scarf_min_max_1958,bertsimas_optimal_2005,yue_expected_2006,popescu_robust_2007,delage_distributionally_2010}. 
	Here we 
	consider an ambiguity set 
	with ellipsoidal structure of the first moment:
	\begin{align} \label{set:fP(bar V, S, gamma)}
{\cal P}_A
		(\gamma) :=
		\left\{
		P \in 
		\mathscr{P}(\mathbb{R}^I)
		\left |
		\begin{array}{l}
			P (V \in 
			\mathcal{V}
			) = 1  \\
			\left\| S_{I-1}^{-1/2}  \left(C \bE_P [V] - \bar V_{I-1} \right)  \right\|^2 \le \gamma
		\end{array}
		\right.
		\right\},
	\end{align}
	where $	\mathscr{P}(\mathbb{R}^I)$ denotes
	the set of all probability measures on the space 
	$\mathbb{R}^I$, $\|\cdot\|$ denotes the Euclidean norm and 
	$S_{I-1}^{1/2}$ is a full rank matrix 
	with $S_{I-1} = \left(S_{I-1}^{1/2} \right)^T S_{I-1}^{1/2}$ and $S_{I-1}^{-1/2} = \left(S_{I-1}^{1/2} \right)^{-1}$.
	The second condition in the set
	${\cal P}_A(\gamma)$
	specifies 
	the range of 
	$C \bE_P [V]$, that is, we consider only the candidate probability 
	distributions with associated means value of $V$ falling within the specified
	ellipsoid centered at the sample mean $\bar V_{I-1}$.
 {The parameter}
 $\gamma$ determines 
	the size of the
	ellipsoid and its choice depends on the DM's confidence in the sample data. The whole idea 
	of the ambiguity set is built on
	the confidence region for multivariate random variables. 
 The following proposition
 states conditions under which the true probability distribution of $V$ falls in the ambiguity set ${\cal P}_A(\gamma)$
 by 
	a proper setting of $\gamma$ 
	for the fixed sample size $N$ and confidence
	level $1-\alpha$.

	\begin{proposition} 
	\label{Prop:bootstrap-confidence}
 For a given $\alpha\in (0, e^{-2}(2-e^{-2}))$, 
 let
	\bgeqn
	\label{eq:gamma_N^2}
	\gamma_N^{(1)} = \frac{
	32R^2 e^2\ln^2(1/(1-\sqrt{1-\alpha}))}{N},
	\edeqn
	where $R := \left\| \Sigma_{I-1}^{-1/2} \right\|$.
Then		there exists a positive integer $N_0$ depending on $\alpha$ 
		such that, for the true probability distribution $P$ of the random vector $V$,
		\begin{equation*}
			\prob\left(P\in  
			{\cal P}_A(\gamma_N^{(1)})
			\right) \geq 1-\alpha
		\end{equation*}
		for all $N\geq N_0$.
	\end{proposition}

\proof
Let $\beta:=1-\sqrt{1-\alpha}$.
As we discussed earlier, under Assumption \ref{Assu-Sigma}, we have 
$$
\left\| \Sigma_{I-1}^{-1/2} (C \bE[V] - \bar V_{I-1})  \right\| 
\le (I-1)R < \infty \ a.s..
$$
Hence, the bounded support guarantees the ``Condition (G)'' in \cite{so_moment_2011}, i.e.,
for any $p \ge 1$, 
$$
\bE \left[ \left\|\Sigma_{I-1}^{-1/2}(C\bbe[V]-\bar{V}_{I-1}) \right\|^p\right] \leq 
[
4R^2 p]^{p/2}.
$$
It follows by Proposition 4 in \cite{so_moment_2011} that
$$
\prob\left(
\left\| \Sigma_{I-1}^{-1/2} (C \bE[V] - \bar V_{I-1})  \right\|
\leq \gamma_N^{(1)}/2   \right)
\geq 1-\beta.
$$
On the other hand, 
since 
$S_{I-1}$ converges to $\Sigma_{I-1}$ weakly by the weak law of large number, 
and $\left\|C \bE[V] - \bar V_{I-1}  \right\|$
is bounded by 
$2$ a.s., then 
there exists $N_0$ depending on $\alpha$ such that 
		$$
		 \prob\left(\left\| S_{I-1}^{-1/2} (C \bE[V] - \bar V_{I-1})  \right\| 
		\leq 	2\left\| \Sigma_{I-1}^{-1/2} (C \bE[V] - \bar V_{I-1})  \right\|
		\right) \geq 1-\beta
		$$
		for all $N\geq N_0$.
Consequently 		
			\bgeq
&&	 \prob\left(
	\left\| S_{I-1}^{-1/2} (C \bE[V] - \bar V_{I-1})  \right\| \le \gamma_N^{(1)}\right)\\
	&\geq&  \prob\left(
	\left\| \Sigma_{I-1}^{-1/2} (C \bE[V] - \bar V_{I-1})  \right\|
	\leq \gamma_N^{(1)}/2  \right) (1-\beta)
	\geq (1-\beta)^2 =1-\alpha,
	\edeq
		which completes the proof.
	\hfill $\Box$

	The 
	proposition states that for a fixed sample size $N$ how large the parameter $\gamma$
should be 
so that the true distribution 
falls in ${\cal P}_A(\gamma)$
with a 
specified confidence level.
This means that, with probability 
at least $1-\alpha$, 
the optimal value of the 
\ref{DUPRO} model computed with such an ellipsoidal ambiguity set provides 
a lower bound for the true optimal value of the expected utility maximization problem, see Theorem 4
in \cite{delage_distributionally_2010} for a similar argument.
From 
a 
decision-making perspective, we want the ambiguity set to be small (a smaller $\gamma$) so that the resulting decision based on the worst distribution from the ambiguity set 
is less conservative. On the other hand, we also hope the true probability distribution 
lies within
the ambiguity set with a specified confidence level
so that 
the model is more relevant. 
Nevertheless, 
a smaller $\gamma$ 
would make this less likely to happen.
	If we are able to recover the true probability distribution with
	a sample of size $N$, then we may set $\gamma=0$. 
	In practice, however, this is very unlikely particularly when the sample size is small. 
	A potential drawback of this approach is that 
	when we use sample covariance 
	to approximate the true covariance, $R$ could be very large when the covariance matrix is nearly singular (the data are nearly correlated), this will result in a large 
	$\gamma_N^{(1)}$.
		We refer readers to 
	discussions 
	in \cite{delage_distributionally_2010}, \cite{so_moment_2011} and our discussions in 
	Section \ref{subsec:The distributionally robust model} for 
	appropriate setting of $\gamma$.  This prompts us to consider the bootstrap approach
 for constructing an ambiguity set.

\subsubsection{The bootstrap approach}
 We
propose to use a percentile-$t$ bootstrap method to specify a confidence region 
of $\mu_{I-1}$.
Specifically, we generate $K$ nonparametric bootstrap resamples of the first $I-1$ components, denoted by $V_{I-1}^{(k,1)}, \dots, V_{I-1}^{(k,N)}$, $k = 1, \dots, K$, based on the empirical distribution (i.e., the uniform distribution on original sample data $V_{I-1}^1,\dots,V_{I-1}^N$).
Let ${\bar V_{I-1}}^{*,k}$ and
	${S_{I-1}^{*,k}}$ denote 
	the sample mean and the sample covariance matrix of $V_{I-1}^{(k,1)}, \dots, V_{I-1}^{(k,N)}$. 
	We consider the  studentized statistics
	$$
	T := \sqrt{N} S_{I-1}^{-1/2} (\bar V_{I-1} - \mu_{I-1})
	$$ 
	and its bootstrap counterparts 
	\bgeqn
	\label{eq:Tk*}
	T^*_k 
	:= \sqrt{N} \left(S_{I-1}^{*,k}\right)^{-1/2} ({\bar V}_{I-1}^{*,k} - \bar V_{I-1}), \quad k=1,\cdots,K.
	\edeqn
	The literature \cite{hall_bootstrap_1987,yeh_balanced_1997,battista_multivariate_2004} addresses multivariate bootstrap confidence regions using 
 data depth and likelihood.
	Here we state the procedure given by \cite{yeh_balanced_1997}. Recall that Tukey's depth of a point $x$ under some distribution $\cF$ is defined as
        \begin{align*}
		{\rm TD}(\cF, x) := \inf_{\|s\| = 1} \int \step{s^T (y - x) \ge 0} d \cF(y).
	\end{align*}
 Let $\cF^*_K$ be the empirical cumulative distribution function built using 
	$(T^*_1,\cdots,T^*_K)$
	and calculate $d_k := {\rm TD}(\cF^*_K, T^*_k)$,
	$k=1,\dots,K$. 
	Denote by $d_{(1)}, \dots, d_{(K)}$ in 
	the decreasing order of the quantities 
	of $d_{1}, \dots, d_{K}$, 
	and let $T^*_{(k)}$ be the statistics such that $d_{(k)} = {\rm TD}(\cF^*_K, T^*_{(k)})$. For a given $\alpha \in (0, 1)$, let $\hat T_{1-\alpha} := \left [T^*_{(1)}, \dots, T^*_{(\lceil (1-\alpha) K \rceil)}\right]$ be an $M \times \lceil (1-\alpha) K \rceil$ dimensional matrix.
	Here, $\lceil \cdot \rceil$ is the ceiling function.
		We construct a convex hull, denoted by $\fW^*_{1-\alpha}$, of $T^*_{(k)}$ for $1 \le k \le \lceil (1-\alpha) K \rceil$, i.e.,
	\bgeqn
	\label{eq:set-fW-Boots}
	\fW^*_{1-\alpha} := \left\{
	\hat T_{1-\alpha} w \in \mathbb{R}^{I-1}\ \Big|\   e_{\lceil (1-\alpha) K \rceil}^T w = 1, \ w \ge 0 \right \}.
	\edeqn
	Then a $100(1-\alpha) \%$ bootstrap confidence region for $\mu_{I-1}$ is obtained as
	\begin{align}
	\label{eq:set:fC}
		\fC^*_{1-\alpha} :=  \left\{ \bar V_{I-1} -  S_{I-1}^{1/2} \tilde w / \sqrt{N}  \ \Big|\  \tilde w \in   \fW^*_{1-\alpha} \right\}.
	\end{align}
	 	On this basis, we define an ambiguity set of the first moment as
	\begin{align} \label{set:fP_B}
{\cal P}_B
		(\alpha) :=
		\left\{
		P \in \mathscr{P}(\mathbb{R}^I) \ |\ P(V \in \mathcal{V}
		) = 1, \ C \bE_P [V] \in \fC^*_{1-\alpha}
		\right\}.
	\end{align}
In this bootstrap approach, $\alpha$ is the critical value of the confidence region. A larger $\alpha$ means less sample points from the bootstrap resampling process are included. Subsequently, $\fC^*_{1-\alpha}$ and 
${\cal P}_B(\alpha)$
are smaller.

\begin{figure}[!htbp]
\vspace{-0.3cm}
\minipage{0.33\textwidth}
 \centering
  \includegraphics[width=4.5cm]{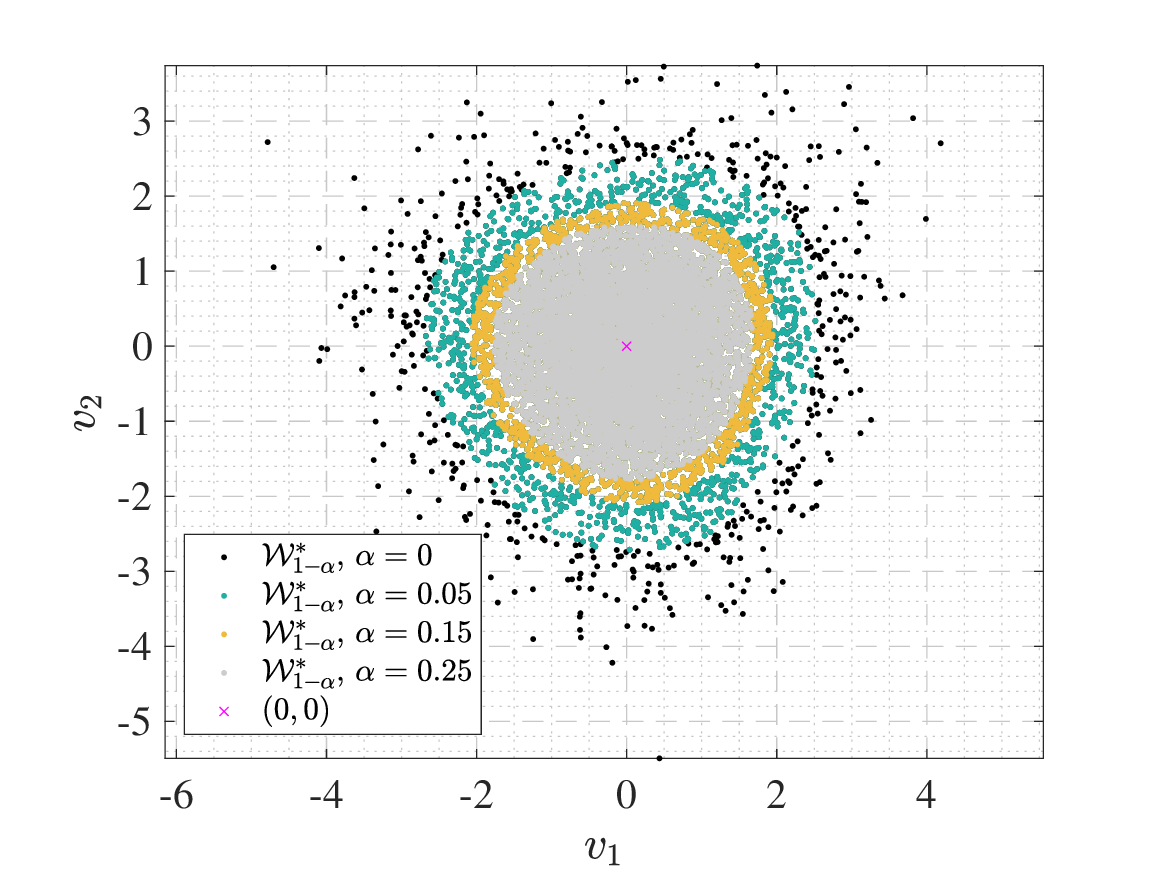}
  \textbf{\tiny (a)}
\endminipage\hfill
\minipage{0.33\textwidth}
  \centering
  \includegraphics[width=4.5cm]{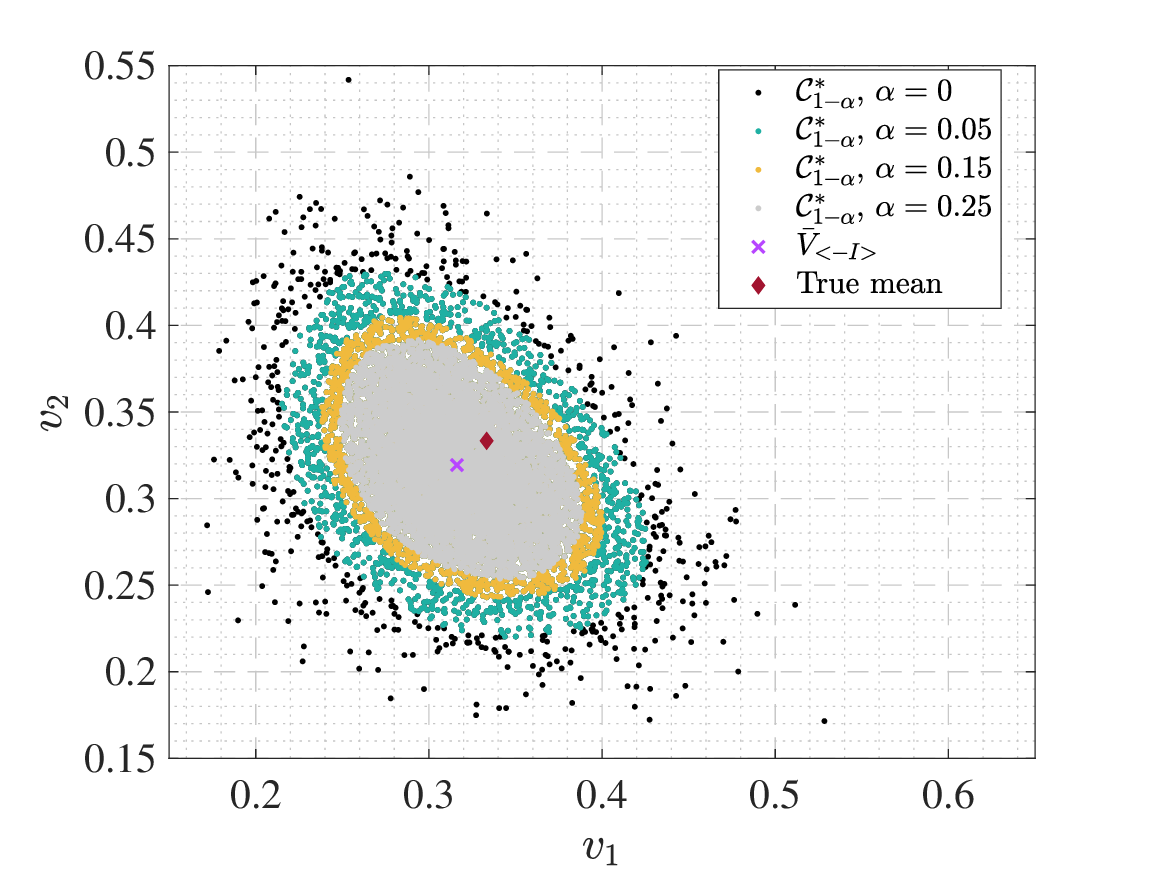}
  \textbf{\tiny (b)}
\endminipage
\hfill
\minipage{0.33\textwidth}
  \centering
  \includegraphics[width=4.5cm]{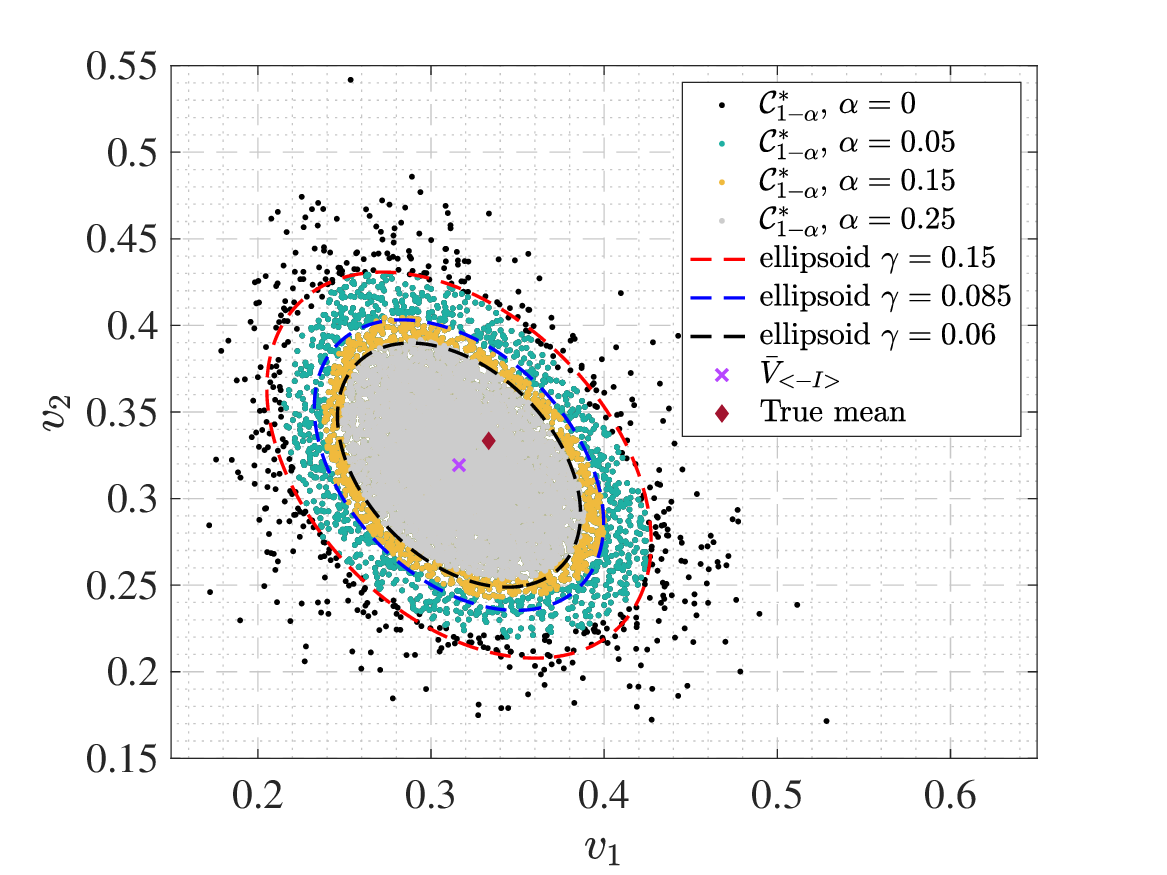}
  \textbf{\tiny (c)}
\endminipage
\caption{\small {$K=1000$, $N=50$}.
\textbf{(a)} convex hull $\fW^*_{1-\alpha}$.
\quad \textbf{(b)} $100(1-\alpha)\%$ bootstrap confidence region $\fC^*_{1-\alpha}$.\quad
\textbf{(c)} $\fC^*_{1-\alpha}$ v.s. ellipsoid 
$\left\{v_{<3-1>}\in \mathbb{ R}^2_+: \|S^{-1/2}_{<3-1>}(v_{<3-1>}-\bar{V}_{<3-1>})\|^2\leq \gamma \right\}$.}
\label{Cstar_vs_ellipsoid}
\vspace{-0.5cm}
\end{figure}

 Observe that for fixed $K$,
$\fW_{1-\alpha}^*$ is a bounded set 
and 
by Corollary 6 in \cite{shawe2003estimating},
$S^{1/2}_{I-1}$ converges to $\Sigma^{1/2}_{I-1}$ 
at an exponential rate
with the increase of $N$.
Thus
$S^{1/2}_{I-1}\tilde{w}/\sqrt{N}\to 0$ as $N\to \infty$ w.p.1. which implies that the ambiguity set
$\fC^*_{1-\alpha}$ 
  shrinks to a singleton 
 at an exponential rate
  as $N\to\infty$.
That is, for a fixed $\alpha$, the 
size of the ambiguity set
is only determined by $N$ (original sample size).
In contrast, the ellipsoidal
method requires one to choose an appropriate size by adjusting $\gamma$ with 
regard to $N$,
such as $\gamma_N^{(1)}$ in (\ref{eq:gamma_N^2}). 
When $N$ is small, $\gamma_N^{(1)}$ could be very large to secure
  the true probability distribution to fall into ${\cal P}_A(\gamma_N^{(1)})$ by Proposition~\ref{Prop:bootstrap-confidence}. However, such a large ambiguity set ${\cal P}_A(\gamma_N^{(1)})$ is too conservative in practice.
  Indeed, for a given $N$, the existing literature lacks a practical method for determining the size of an ellipsoid ambiguity set. Differently,     the bootstrap approach offers the convenience of selecting an appropriate critical value $\alpha$ with a probabilistic interpretation.
Here we illustrate how the confidence region constructed by the bootstrap looks like
in the case that the dimension of $V_{I-1}$ is $2$, the original sample size is 
$N=50$ and the number of the resamples is $K=10,000$. 
Figures~\ref{Cstar_vs_ellipsoid}~(a) and (b) depict $\fW^*_{1-\alpha}$
and $\fC_{1-\alpha}^*$ for different $\alpha$ values 
ranging from $0$ to $0.25$.
Figure~\ref{Cstar_vs_ellipsoid}~(c) depicts the minimal ellipsoid containing
$\fC_{1-\alpha}^*$ for different $\alpha$ values.
Since each ellipsoid may be regarded as a confidence region, we can see
from the figures the sizes of the confidence regions with specified 
confidences of $95\%$, $85\%$ and $75\%$.
 Moreover, we can see that the sizes of confidence regions
 are considerably smaller than those defined via
 ${\cal P}_A(\gamma)$.
For example, 
when $N=50$, and $\alpha=0.25$,
	$\gamma_N^{(1)}=70$.

In this paper, we adopt the bootstrap approach
to construct the ambiguity set not only 
because the size of the ambiguity set is relatively smaller but also 
it is often difficult to obtain a large number of observations of $V$. Moreover, we will also see in the next section that 
 \ref{DUPRO} can be easily 
reformulated 
at an MILP which is relatively easier to solve.

	\section{Reformulation of the \ref{DUPRO}
 model
 }
	\label{sec:Reformulations and Solution Methods}
	In this section, we develop computational approaches for solving \ref{DUPRO} with the ambiguity sets ${\cal P}_A(\gamma)$ and ${\cal P}_B(\alpha)$
	defined  as in \eqref{set:fP(bar V, S, gamma)} and \eqref{set:fP_B}, respectively. Section~\ref{subsec:General Increasing Utility Function} addresses the case with increasing utility functions, while Section~\ref{subsec:Increasing Concave Utility Function} discusses the case with risk-averse utility functions.
		We begin by presenting
		a reformulation of the objective function
$\bbe_P[u(x;V)]$
  of \ref{DUPRO} as a linear function of the 
expectation of random vector $V$ in the next proposition.
To ease the notation,
we 
  recall some of the key notations introduced in Section~\ref{sec:Random_additive_PLUF}.
 $
 \step{t_{m, i-1} < x_m \le t_{m, i}}
 $ 
 represents an indicator function of interval 
 $(t_{m, i-1}, t_{m, i}]$,
 $x=(x_1,\cdots,x_M)^T\in \R^M$,
 $u_m$ is a PLUF defined over $[a_m,b_m]$ with breakpoints $\{t_{m,i}:i\in \{0\}\cup \fI_m\}$, for $m\in \fM$, $v=(v_1^T,\cdots,v_M^T)^T\in \R^{I}$ with $\fI_m=\{1,\cdots,I_m\}$, $\fM=\{1,\cdots,M\}$, $I=\sum_{m\in \fM}I_m$, and 
$u(x;v)=\sum_{m\in \fM} u_m(x_m;v_m)$. 
	\begin{proposition} \label{pp:equ_models}
		Let 
		$$	g_{m,i}(x_m)  := \step{t_{m, i-1} < x_m \le t_{m, i}} \quad \inmat{and}\quad
			h_{m,i}(x_m)  := \left( \frac{x_m - t_{m, i-1}}{t_{m, i} - t_{m, i-1}} \right)  g_{m,i}(x_m),
		$$
		for $i\in \fI_m, m\in \fM.$
		Let
		\begin{align*}
			f_m(x_m) := \left[
			\begin{array}{c}
				h_{m,1}(x_m) + \sum_{j \in \fI_m \setminus \{1\}} g_{m,j}(x_m) \\
				h_{m,2}(x_m) + \sum_{j \in \fI_m \setminus \{1, 2\}} g_{m,j}(x_m) \\
				\vdots \\
				h_{m, I_m}(x_m)
			\end{array}
			\right]_{I_m \times 1}
		\end{align*}
	and 
			$f(x) := \left (f_1(x_1)^T, \dots, f_M(x_M)^T \right)^T$.
		Then
		\begin{itemize}
	\item[(i)]	$	u_m (x_m; v_m) = v_m^T f_m(x_m)$,
	for $m\in \fM$
	and 
	\begin{align}
			u (x; v) = v^T f(x).
			\label{eq:u-vec-rep}
		\end{align}
	\item[(ii)] \ref{DUPRO} is equivalent to
		\begin{align} \label{mod:PPO-linear-V}
			\max_{x \in \fX} \min_{P \in 
			{\cal P}
			} \bE_P[V^T] f(x)
		\end{align}
		and
		\begin{align} \label{mod:PPO}
			\max_{x \in \fX} \min_{v \in \fF(
			{\cal P})
			} \left\{ u(x; v) = v^T f(x) \right\},
		\end{align}
		where $\fF({\cal P})
		:= \{\bE_P(V) \;|\;  P \in 
		{\cal P}
		\}$.
	\end{itemize}
	\end{proposition} 
Formulation (\ref{mod:PPO}) implies that \ref{DUPRO} 
can be converted into a maximin robust optimization problem where the 
ambiguity set $\fF({\cal P})$
is composed of the expected values of $V$
each of which corresponds to an expected utility function.
The worst-case expected utility function 
$\min_{P \in 
{\cal P}
} \bE_P[V^T] f(x)$
is 
determined by the minimum of $\{v^Tf(x): v \in \fF(
{\cal P})
\}$.
Thus, we may solve \ref{DUPRO} by solving problem (\ref{mod:PPO}).
 This is indeed one of the main advantages
for us to adopt the particular form
of PLUF in (\ref{Random-MAPLU}).
However, the structure of $f(x)$ is still too complex 
 for us to derive a tractable formulation of
 \ref{DUPRO}.
To address the issue,
we derive an alternative representation of $f(x)$.
Let 
$y = (y^T_1, \dots, y^T_M)^T \in \bR^I$ and $z = (z^T_1, \dots, z^T_M)^T \in \{0, 1\}^I$, where $y_m = (y_{m,1}, \dots y_{m, I_m})^T \in \bR^{I_m}$ and $z_m = (z_{m,1}, \dots z_{m, I_m})^T \in \{0, 1\}^{I_m}$ for $m \in \fM$. 	Let 
	\begin{align}\label{eq:YZ}
		Y := \left[
		\begin{array}{ccc}
			Y_1 & & \\
			& \ddots & \\
			& & Y_M
		\end{array}
		\right]_{I \times I},
		\quad 
		Z := \left[
		\begin{array}{ccc}
			Z_1 & & \\
			& \ddots & \\
			& & Z_M
		\end{array}
		\right]_{I \times I}, 
	\end{align}
	where the blocks $Y_m$ and $Z_m$ for $m \in \fM$ are
	\begin{align*}
		Y_m := \left[
		\begin{array}{ccc}
			\frac{1}{t_{m,1} - t_{m,0}} & & \\
			& \ddots & \\
			& & \frac{1}{t_{m, I_m} - t_{m, I_{m-1}}}
		\end{array}
		\right]_{I_m \times I_m},
		\quad
		Z_m := \left[
		\begin{array}{cccc}
			\frac{-t_{m,0}}{t_{m,1}-t_{m,0}} & 1 & \dots & 1 \\
			& \frac{-t_{m,1}}{t_{m,2}-t_{m,1}} & \dots & 1 \\
			&  & \ddots & \vdots \\
			&  &  & \frac{-t_{m,I_m-1}}{t_{m,I_m}-t_{m,I_m-1}}
		\end{array}
		\right]_{I_m \times I_m}.
	\end{align*}
	Denote five other block-diagonal matrices
 \begin{subequations}
  \label{eq:BD}
	\begin{align}
		& B := \left[
		\begin{array}{ccc}
			B_1 & & \\
			& \ddots & \\
			& & B_M
		\end{array}
		\right]_{I \times I},
		\quad
		D := \left[
		\begin{array}{ccc}
			D_1 & & \\
			& \ddots & \\
			& & D_M
		\end{array}
		\right]_{I \times M},
		\quad
		E := \left[
		\begin{array}{ccc}
			e_{I_1} & & \\
			& \ddots & \\
			& & e_{I_M}
		\end{array}
		\right]_{I \times M}, \\
		&	H^- := \left[
		\begin{array}{ccc}
			H^-_{1} & & \\
			& \ddots & \\
			& & H^-_{M}
		\end{array}
		\right]_{I \times I},
		\quad
		H^+ := \left[
		\begin{array}{ccc}
			H^+_{1} & & \\
			& \ddots & \\
			& & H^+_{M}
		\end{array}
		\right]_{I \times I}, 
	\end{align}
 \end{subequations}
	where the blocks $B_m$, $D_m$, $H^-_m$, and $H^+_m$ for $m \in \fM$ are
 \begin{subequations} \label{eq:B_D_H}
	\begin{align}
 \label{eq:B_mD_m}
		& B_m := \left[
		\begin{array}{ccccc}
			1 &  & & & \\
			1 & 1 &  & & \\
			\vdots  & \vdots & \ddots & & \\
			1 & 1 &\cdots & 1  &
		\end{array}
		\right]_{I_m \times I_m},
		\qquad
		D_m := \left[
		\begin{array}{c}
			t_{m,1} \\
			\vdots\\
			t_{m,I_m}
		\end{array}
		\right]_{I_m \times 1}, \\
		&	H^-_m := \left[
		\begin{array}{ccc}
			t_{m,0} & & \\
			& \ddots & \\
			&  & t_{m,I_m-1}
		\end{array}
		\right]_{I_m \times I_m},
		\quad
		H^+_m := \left[
		\begin{array}{ccc}
			t_{m,1} & & \\
			& \ddots & \\
			& & t_{m,I_m}
		\end{array}
		\right]_{I_m \times I_m}.
	\end{align}
 \end{subequations}
With all these new notations, we are ready to 
represent $f(x)$ as follows:	
\bgeqn
f(x)= Yy+Zz
\label{eq:f(x)}
\edeqn
under the conditions below
 \bgeqn 
  E^Tz=e_M,\quad
            E^Ty =x,\quad
            H^{-}z-y\leq 0,\quad
            H^+z-y \geq 0,\quad
            z\in \{0,1\}^I.
 \label{eq:f(x)-Y-Z}
 \edeqn
It is important to note that 
for a fixed $x$,
$(y,z)$ satisfying (\ref{eq:f(x)-Y-Z}) is unique. 
 Representation  of $f(x)$
via (\ref{eq:f(x)})-(\ref{eq:f(x)-Y-Z}) is a crucial step towards tractable formulation of (\ref{mod:PPO})
in the rest of this section. This is 
another main advantage that we adopt the particular form of PLUF in (\ref{Random-MAPLU}).

	\subsection{
 PLUF without concavity
 }
	\label{subsec:General Increasing Utility Function}

	We 
 consider the case that the random utility functions are piecewise linear but are not necessarily concave, and 
 discuss 
 reformulation of \ref{DUPRO} with the ambiguity set 
	$\mathcal{P}_B(\alpha)$
	(see (\ref{set:fP_B}))
		\bgeqn
	\label{eq:DUPRO_FB_bootstrap}
	\max_{x\in\fX} \min_{v \in \fF(
\mathcal{P}_B(\alpha)
	)} v^T f(x),
	\edeqn
where 
	\begin{align}
	\label{eq:fF_PB_alpha}
		\fF(
		\mathcal{P}_B(\alpha)
		)  = \left \{
		v \in \bR^I
		\left |
		\begin{array}{l}
			v \ge 0 \\
			e_I^Tv = 1 \\
			Cv=\bar V_{I-1}-S^{1/2}_{I-1}\tilde{w}/\sqrt{N},\; \tilde w \in \fW_{1-\alpha}^* 
		\end{array}
		\right.
		\right \}.
	\end{align}
 The
 next proposition
 states that the maximin problem
(\ref{eq:DUPRO_FB_bootstrap})	can be reformulated as a single MILP.

	\begin{proposition}\label{prop:PRO_B_equ}
		\ref{DUPRO} with the ambiguity set 
		$\mathcal{P}_B(\alpha)$
		defined as in \eqref{set:fP_B} is equivalent to
		\begin{subequations} \label{pp:PRO_B_equ}
			\begin{align}
				\max_{x, y, z, \eta, \pi, \tau}  \ &  \bar V_{I-1}^T \eta + \pi + \tau\\
				{\rm s.t.}  \quad~ &  C^T\eta + e_I \tau \le Yy+Zz,\\  
				& N^{-1/2} \hat T_{1-\alpha}^T \left(S_{I-1}^{1/2}\right)^T \eta + e_{\lceil (1-\alpha) K \rceil} \pi \le 0, \label{pp:PRO_B_equ-2}\\
				& E^T z = e_M, \label{pp:PRO_B_equ-3}\\
				& E^T y = x, \label{pp:PRO_B_equ-4}\\
				& H^- z - y \le 0, \label{pp:PRO_B_equ-5}\\
				& H^+ z - y \ge 0, \label{pp:PRO_B_equ-6}\\
				& z \in \{0, 1\}^I, \label{pp:PRO_B_equ-7}\\
				& x \in \fX, \label{pp:PRO_B_equ-8}
			\end{align}
		\end{subequations}
		where $Y,Z,H^-,H^+$ are defined as in (\ref{eq:YZ})-(\ref{eq:B_D_H}), and $\hat T_{1-\alpha}$ is defined 
		as in (\ref{eq:set-fW-Boots}).
	\end{proposition}
\proof
	By Proposition \ref{pp:equ_models}, \ref{DUPRO} with the ambiguity set 
	$\mathcal{P}_B(\alpha)$
	is equivalent to \eqref{mod:PPO} with
	$\fF(\mathcal{P}_B(\alpha))$.
	Also, for any $Cv \in \fC^*_{1-\alpha}$, there exists a $w \in \bR^{\lceil (1-\alpha)K \rceil}$ such that $Cv = \bar V_{I-1} - N^{-1/2} S_{I-1}^{1/2} \hat T_{1-\alpha} w$, $e_{\lceil (1-\alpha) K \rceil}^T w = 1$, and $w \ge 0$. Hence, the inner minimization problem of \eqref{mod:PPO} is written as
	\begin{align*}
		\min_{v, w} \ & v^T f(x) \\
		\text{s.t.} \
		& Cv + N^{-1/2} S_{I-1}^{1/2} \hat T_{1-\alpha} w = \bar V_{I-1},\\
		& e_{\lceil (1-\alpha) K \rceil}^T w = 1, \\
		& e_I^T v = 1, \\
		& v \ge 0, \ w \ge 0.
	\end{align*}
	Let $\eta\in \mathbb{R}^{I-1}$, $\pi\in \mathbb{R}$, and $\tau\in \mathbb{R}$ be the dual variables regarding the first three constraints in the problem above. We obtain the dual problem as
	\begin{align*}
		\max_{\eta, \pi, \tau}  \ &  \bar V_{I-1}^T \eta + \pi + \tau\\
		\text{s.t.}  \ &  C^T\eta + e_I \tau \le f(x),\\  
		& N^{-1/2} \hat T_{1-\alpha}^T \left(S_{I-1}^{1/2}\right)^T \eta + e_{\lceil (1-\alpha) K \rceil} \pi \le 0.
	\end{align*}
The conclusion follows
by using the unique representation of $f(x)$ via (\ref{eq:f(x)})-(\ref{eq:f(x)-Y-Z}).

\hfill $\Box$
	
From (\ref{pp:PRO_B_equ}), we can see that as $\alpha$ increases, there will be fewer inequality constraints in (\ref{pp:PRO_B_equ-2}) and consequently 
 the optimal value of (\ref{pp:PRO_B_equ}) will increase. 
Problem \eqref{pp:PRO_B_equ} has $I$ binary variables. Its scalability is contingent on the number of attributes and the number of pieces in each piecewise linear single-attribute utility function. The number of pieces depends on the complexity of the problem, empirical data, expert judgment, or various mathematical and statistical techniques. In some cases with relatively straightforward and homogeneous preferences, a 
PLUF
may have just a few pieces or segments to capture the essential variations in utility. In other cases, analysts may choose to use more complex 
PLUFs
with multiple segments to better capture the nonlinear nature of preferences. In Section \ref{case_study1}, we conduct an experiment to test the computational time of problem \eqref{pp:PRO_B_equ}. The test results exhibit that the computational time has an exponential growth rate as $I$ increases.

	\subsection{
Random concave PLUF  case}
	\label{subsec:Increasing Concave Utility Function}
	We now turn to discuss the case that 
 the random utility function of each
 single-attribute 
is concave. In this case, 
	$\mathcal{V}_c$
	defined as in \eqref{set:fV_C} is the support set of the random vector $V$ of increments. In this section we substitute 
	$\mathcal{V}_c$ for $\mathcal{V}$
	in the ambiguity set
	$\mathcal{P}_B(\alpha)$ given in
	\eqref{set:fP_B}.
	We discuss a reformulation of \ref{DUPRO} with the ambiguity set
	$\mathcal{P}_B(\alpha)$
	in the proposition below.

	\begin{proposition} \label{pp:Cav_Reform_B}
		\ref{DUPRO} with the ambiguity set
		$\mathcal{P}_B(\alpha)$
		is equivalent to
		\begin{subequations}
		\label{eq:Cav_Reform_B}
			\begin{align}
				\max_{x, \eta, \pi, \tau, \zeta, \lambda} \ & \bar V_{I-1}^T \eta + \pi + \tau \\
				\text{\rm s.t.}~~~ \
				& C^T \eta + e_I \tau + A^T \zeta - B^T \lambda \le 0,\\
				& N^{-1/2} \hat T_{1-\alpha}^T \left(S_{I-1}^{1/2}\right)^T \eta + e_{\lceil (1-\alpha) K \rceil} \pi \le 0,\\
				& D^T \lambda \le x,\\
				& E^T \lambda \le e_M, \\
				& \zeta \ge 0, \ \lambda \ge 0,\\
				& x \in \fX,
			\end{align}
		\end{subequations}
  where $A$ is defined as in (\ref{set:fV_C}),
  $B$, $D$, and $E$ are defined as in (\ref{eq:BD})-(\ref{eq:B_D_H}).
	\end{proposition}
\proof
		We first prove that the inner minimization problem of \ref{DUPRO} can
	be represented as
	\begin{subequations}
	\begin{align}
		\min_{y, z, v, w} \ & x^T y + e_M^T z \label{eq:concave_inner-a}\\
		\text{s.t.}~ \
		& Cv + N^{-1/2} S_{I-1}^{1/2} \hat T_{1-\alpha} w = \bar V_{I-1},\label{eq:concave_inner-b}\\
		& e_{\lceil (1-\alpha) K \rceil}^T w = 1, \label{eq:concave_inner-c}\\
		& e_I^T v = 1,\label{eq:concave_inner-d}\\
		& A v \ge 0,\label{eq:concave_inner-e} \\
		& D y + E z - B v \ge 0, \label{eq:concave_inner-f}\\
		& v \ge 0, \ w \ge 0, \ y \ge 0, \ z \ge 0.\label{eq:concave_inner-g}
	\end{align}
		\end{subequations}
The $\fF$-mapping of $
\mathcal{P}
	_B(\alpha)$ with the support 
	$\mathcal{V}_c$
	is written as
	\begin{align*}
		\fF(
	\mathcal{P}_B(\alpha)
		) =
		\left\{
		v \in \bR^I_+  \;\Big|\;
		\text{$v$ satisfies conditions \eqref{eq:concave_inner-b} - \eqref{eq:concave_inner-e}}
		\right\}.
	\end{align*}
	For any given $v \in \fF(
	\mathcal{P}_B(\alpha)
	)$, $u_m(x_m; v_m)$ is a piecewise linear increasing concave function with breakpoints $t_i$ and values 
	$\sum_{j=1}^i v_{m,j}$, $i = 1,
	\dots, I_m$. It can be written as the minimization of all subgradients at points $\left(t_i, \sum_{j=1}^i v_{m,j} 
	\right)$ as
	\begin{align*}
		u_m(x_m; v_m) = \min_{y_m, z_m,v_m} \ & x_m y_m + z_m \\
		\text{s.t.} ~~~\ 
		& D_m y_m + z_m e_{I_m} - B_m v_m \ge 0,\\
		& y_m \ge 0, \ z_m \ge 0,
	\end{align*}
where $D_m$ and $B_m$ are given in (\ref{eq:B_D_H}).
	Therefore, we combine all $u_m$ to the multi-attribute utility function
	\begin{align*}
		u(x; v) = \sum_{m \in \fM} u_m(x_m; v_m)  = \min_{y, z, v} \ & x^Ty + e^T_Mz \\ 
		\text{s.t.} \ & D y + E z - B v \ge 0, \\
		& y \ge 0, \ z \ge 0.
	\end{align*}
	This gives the objective \eqref{eq:concave_inner-a} and constraints \eqref{eq:concave_inner-f} and \eqref{eq:concave_inner-g}.
	
	Letting $\eta\in \mathbb{R}^{I-1}$, $\pi\in \mathbb{R}$, $\tau\in \mathbb{R}$, $\zeta\in \mathbb{R}^{I-M}$, and $\lambda\in \mathbb{R}^I$ 
	be the dual variables regarding the constraints, respectively, in the problem above, we can derive the Lagrange dual
	\begin{align*}
		\max_{\eta, \pi, \tau, \zeta, \lambda} \ & \bar V_{I-1}^T \eta + \pi + \tau \\
		\text{s.t.}~~ \
		& C^T \eta + e_I \tau + A^T \zeta - B^T \lambda \le 0,\\
		& N^{-1/2} \hat T_{1-\alpha}^T \left(S_{I-1}^{1/2}\right)^T \eta + e_{\lceil (1-\alpha) K \rceil} \pi \le 0,\\
		& D^T \lambda \le x,\\
		& E^T \lambda \le e_M, \\
		& \zeta \ge 0, \ \lambda \ge 0.
	\end{align*}
	A combination of the above dual program with the outer maximization problem yields (\ref{eq:Cav_Reform_B}). \hfill $\Box$

As in Proposition 
\ref{prop:PRO_B_equ}, the optimal value of (\ref{eq:Cav_Reform_B}) increases as $\alpha$ increases.
Note that Proposition~\ref{prop:PRO_B_equ} states the MILP reformulation of \ref{DUPRO} when considering general increasing piecewise linear single-attribute utility functions. Proposition \ref{pp:Cav_Reform_B}
shows the LP reformulation of  \ref{DUPRO} if all the single-attribute utility functions are increasing and concave.	
This phenomenon is consistent with the reformulations 
in the literature of  deterministic 
PRO models albeit our DUPRO framework 
is significantly different from deterministic PRO frameworks. 
This is primarily because
we use the increment-based piecewise linear random utility function, which allows us to recast \ref{DUPRO}
as a deterministic maximin robust optimization problem (\ref{mod:PPO}).
In the case when the piecewise linear random utility function is concave,
we can use the support function method, 
as in the deterministic PRO reformulations,
to reformulate the problem
$\min_{v \in \fF(
			{\cal P})
			}  u(x; v),
$
	where $\fF({\cal P})
		:= \{v=\bE_P(V) \;|\;  P \in 
		{\cal P}
		\}$, as an LP.

\section{
	Quantification of modelling errors of the PLUF-based \ref{DUPRO}
	}
	\label{sec:PLA_random utility}
	
	In the preceding sections, we focus on 
	\ref{DUPRO} with 
	the DM's true utility function of each attribute having a piecewise linear structure.
	In practice, the true random utility function is unknown and does not necessarily have a piecewise linear structure. 
What one can usually 
 do is to use nonparametric utility assessment methods 
 to obtain the DM's 
 evaluation (scoring) 
at different levels of attributes
and
then construct a 
PLUF as an approximation.
This means that we may treat the proposed 
PLUF-based \ref{DUPRO} 
as an approximation of general expected random utility maximization problems.
In this section, we quantify the model (approximation)
errors.
	Specifically, we consider the following
	expected utility 
	maximization problem with deterministic attributes
	\begin{align} \label{mod:PRO-general}
	\vt:=\max_{x \in \fX} 
		\bE [U(x)],
	\end{align}
	where $U(x) = \sum_{m\in  \fM} U_m(x)$,
	$U_m$ is the DM's true 
single-attribute utility function of attribute $m$.
	Without loss of generality, we assume that
 $U_m$
 is a general continuous and nondecreasing
	random 
 utility function of attribute $m$ 
 defined on the domain $[a_m, b_m]$ for attribute $m \in \fM$ but it does not necessarily have a piecewise linear structure. 
	Let 
	$\fX \subseteq \otimes_{m \in \fM} [a_m, b_m] = [a_1, b_1] \times \dots \times [a_M, b_M]$ be a compact set. 
	To ease the exposition, we
	assume that $U$ is normalized with   $U(a_1, \dots, a_M) = 0$ and $U(b_1, \dots, b_M) = 1$ almost surely.

	\subsection{Static piecewise linear approximation }
	\label{subsec:Static piecewise linear approximation}
	We begin by discussing the piecewise linear approximation of $U(x)$ and its impact on the optimal value and optimal solutions.
	For each fixed $m\in  \fM$, let $t_{m,i}$, for $i \in \fI_m$, be defined 
	as in (\ref{def:break_points}). 
	Let 
	\begin{align*}
		V_{m, i} := U_m(t_{m, i}) - U_m(t_{m, i-1}), \quad i \in \fI_m,
	\end{align*}
	be the increment of
	$U_m$ over interval $[t_{m,i-1},t_{m,i}]$
	and $V_m=(V_{m, i},\cdots,V_{m,I_m})^T$.
 Throughout this section,
 we let the set of the breakpoints
 $\{t_{m,i}: i \in \fI_m, m\in  \fM\}$ 
 be fixed.
	We construct a piecewise utility function
	over $[a_m,b_m]$, denoted by $u_m(\cdot; V_m)$,
	with breakpoints $t_{m,i}$,  $i \in \fI_m$, 
	and
	\begin{align} \label{eq:u_approximation}
		u_m(t_{m,i}; V_m) = U_m(t_{m,i}), \quad \text{for} \; i \in \fI_m
	\end{align}
	and  use $u_m(\cdot; V_m)$ to approximate $U_m$
	for $m\in  \fM$.
	Let $V:=(V_1^T,\cdots,V_M^T)^T$ and 
    $$
	u(x;V):= 
	\sum_{m\in \fM} u_m(x_{m}; V_m).
	$$
	Note that the dimension of vector $V$ is
    $I =\sum_{m=1}^M I_m$.
	We 
	propose to obtain an 
	approximated optimal value and optimal solution of problem (\ref{mod:PRO-general}) by solving the following piecewise linear approximated expected utility maximization problem:
	\begin{align} \label{mod:PRO-PLA}
	\vt_I:=	\max_{x \in \fX} 
		\bE [u(x; V)].
	\end{align}
 	Let  $X^*$ and $X_I$ be the
 	respective 
	sets of optimal solutions of problems (\ref{mod:PRO-general}) and (\ref{mod:PRO-PLA}).
	The subscript $I$ indicates the optimal value depends 
	on the total number of breakpoints $I$ in the piecewise linear approximation. Of course, it also depends 
	on the location of these points.
	Let
	\begin{equation}
		\Delta := \max_{m \in \fM} \max_{i \in \fI_m} (t_{m,i}-t_{m,i-1}).
	\end{equation}
	Obviously in order to secure a good approximation of $\vt$ by $\vt_I$,
	we need $\Delta$ to be sufficiently small. In the case when
	the breakpoints are evenly spread, this is equivalent to 
	setting $I$ a large value.
To see this more clearly,
let us consider the simplest case that
$I_1=I_2=\cdots=I_M$ where $I_m$ is the number of breakpoints in piecewise 
approximation of the random utility function of attribute $m$. In that case
$
I_m=\frac{I}{M},
$
where $I$ denotes the total number of breakpoints of piecewise linear utility functions of all attributes. 
The requirement on $\Delta\leq \epsilon$ means that
\bgeq 
\max_{m\in \fM} \frac{b_m-a_m}{I_m}=
\max_{m\in \fM} \frac{b_m-a_m}{I/M}\leq \epsilon.
\edeq 
The latter is equivalent to 
$
 I \geq \max_{m\in \fM} 
 \frac{M(b_m-a_m)}{\epsilon}$, 
 which means $I=O(1/\epsilon)$.
	Unless specified otherwise, we assume in the rest of discussions that
	$I\to \infty $ ensures $\Delta\to 0$.
The following proposition
 addresses the approximation of problem (\ref{mod:PRO-general}) by problem (\ref{mod:PRO-PLA}) in terms of the optimal value and optimal solutions.

	\begin{proposition}
		\label{P-utility-appr-PL}
		Suppose that for 
		$ m\in  \fM$, $U_m$ is Lipschitz continuous over
		$[a_m,b_m]$ with random Lipschitz modulus $L_m$ almost surely, where
		$\bE [L_m] <\infty$ 
		and the expectation is taken w.r.t the probability distribution of the random factors
		underlying $U(x)$.
		Then the following assertions hold. 
		\begin{itemize}
			
			\item[(i)] For any $\epsilon > 0$,
			\begin{equation}
	|\vartheta-\vartheta_I| \leq \bbe[L]\Delta, 
						\label{eq:utility-approx}
			\end{equation}
			where 
	$L := \sum_{m \in \fM} L_m$.

			\item[(ii)] Let $\{x_I\}$ be a sequence of optimal solutions obtained from solving problem 
			(\ref{mod:PRO-PLA}). Then every cluster point of the sequence is an optimal solution of problem (\ref{mod:PRO-general}), that is,
			\begin{equation}
				\displaystyle{\lim_{I\to \infty}}
			\mathbb{D}(X_I, X^*)
			=0, 
			\end{equation}
			where $\mathbb{D}(A,B)$ denotes the access distance of set $A$ over set $B$.
			
		\end{itemize}
		
	\end{proposition}

\proof
	Part (ii)  follows directly from Part (i) and the well-known stability results in parametric programming, see e.g.
	 \cite[Lemma 3.8]{liu_stability_2013}.
	We only prove Part~(i). 
	Observe first that 
	$$
	|\vartheta-\vartheta_I| \leq \sup_{x\in 
		\fX}
|\bbe[u(x; V)]-\bbe[U(x)]|.
	$$
	By the monotonicity, Lipschitz continuity of the utility function, and equality \eqref{eq:u_approximation},
	$$
	|u_m(x_m; V_m)- U_m(x_m)| \leq |U_m(t_{m,i})-U_m(t_{m,i-1})|
	\leq L_m |t_{m,i}-t_{m,i-1}|
	\leq L_m \Delta, \ a.s.
	$$
	for any $x_m\in [t_{m,i-1},t_{m,i}]$, which ensures that
	$$
	\sup_{x_m\in [a_m,b_m]}
|\bbe[u_m(x_m; V_m)]- \bbe[U_m(x_m)]| \leq \bbe[ L_m] \Delta,
	$$
	for $m \in \fM$ and hence,
	$$
	\sup_{x\in \fX}
	|\bbe[u(x; V)]- \bbe[U(x)]| 
	\le
	\sup_{x\in \otimes_{m \in \fM} [a_m, b_m]} 
	|\bbe[u(x; V)]- \bbe[U(x)]| 
	\leq 
	\sum_{m \in \fM} 
	\bbe[L_m] \Delta = \bbe[L] \Delta.
	$$
The proof is complete. \hfill $\Box$

	The result provides a theoretical guarantee for using a piecewise linear utility function model 
 parameterized by a vector of random increments to approximate the true
 random utility function in absence of complete information on the latter.
  In other words, one may use a random PLUF as we discussed in Section~\ref{sec:Random_additive_PLUF} to represent 
  a DM's true random utility function provided that 
  it is Lipschitz continuous 
  over the specified compact set and the breakpoints of the PLUF are sufficiently dense in the set.
  Lipschitz continuity over a compact is fulfilled by most utility functions in the literature.

	\subsection{Sample average approximation} 
	\label{subsec:Sample average approximation}
	
	In practice, 
	the true probability distribution of $V$ is often unknown,
	but it is possible to obtain some 
	observations from empirical data
 such as scoring.
This means that instead of solving problem (\ref{mod:PRO-PLA}), 
we often solve the sample average approximation of the problem.
 Research on SAA is well-documented,
	see e.g.~\cite{ShH00,Hom08,shapiro2021lectures}.
	Let $V^1,\cdots,V^N$ denote an iid 
	random sample of $V$ and 
	$\bar V$ the sample mean.
 	By Proposition \ref{pp:equ_models}, we propose to approximate $\bE [U(\cdot)]$ using the sample average
	$$
	\frac{1}{N} \sum_{n=1}^N u(\cdot; V^n) = \frac{1}{N} \sum_{n=1}^N f(\cdot)^T V^n = f(\cdot)^T \bar{V} = 
	u(\cdot; \bar{V}).
	$$
The following proposition
 gives a qualitative description 
	of such approximation.

	\begin{proposition}
		\label{P-u-appro-N} 
		Under the settings and conditions of Proposition \ref{P-utility-appr-PL},
		for any $\epsilon > 0$ and $\delta > 0$, there exists $N_0 > 0$ such that
		\begin{equation}
			\prob\left(
			\sup_{x\in \fX}
			|u(x; \bar{V}) - \bE[U(x)] |
			\geq \epsilon\right) \leq \delta,
			\label{eq:utility-approx-N-0}
		\end{equation}
		for all $N\geq N_0(\epsilon,\delta) =
O(\ln \delta/\epsilon^2)$ 
and $\Delta\leq \frac{\epsilon}{2\bE[L]}$,
		where 
		$\Delta$ is defined as in Proposition \ref{P-utility-appr-PL}.
	\end{proposition}

\proof
	We write $\bar{V} = (\bar V_1^T, \dots, \bar V_M^T)^T$,
	where $\bar{V}_m$ is the
	associated sample average for the vector of
	increments regarding the 
single-attribute utility function of attribute $m$.
	By the triangle inequality, we have
	\bgeqn
	|u_m(x; \bar V_m) - \bbe[U_m(x)]| &\leq&
	|u_m(x; \bar V_m) - \bbe[u_m(x; V_m)]| \nonumber\\
	&&+
	|\bbe[u_m(x; V_m)] - \bbe[U_m(x)]|, \; \text{for}\; m=1,\cdots,M.\quad
\label{eq:u_mEU_m}
	\edeqn
	By the monotonicity, Lipschitz continuity of the utility function,
	and equality \eqref{eq:u_approximation},
	$$
	|u_m(x_m; V_m)- U_m(x_m)| \leq |U_m(t_{m,i})-U_m(t_{m,i-1}))|
	\leq L_m |t_{m,i}-t_{m,i-1}|
	\leq L_m \Delta
	\ a.s.,
	$$
	for all $x_m\in [t_{m,i-1},t_{m,i}]$, which implies that
	$$
	\sup_{x_m\in [a_m,b_m]}  |\bbe[u_m(x; V_m)] - \bbe[U_m(x)]| \leq \sup_{x_m\in [a_m,b_m]}  \bbe |u_m(x; V_m) -U_m(x)| \le \bbe[L_m]\Delta.
	$$
	Thus
	\bgeqn
	\sup_{x\in 
		\fX
	}  |\bbe[u(x; V)] - \bbe[U(x)]|\leq
	\bbe\left[\sum_{m \in \fM} \sup_{x_m\in [a_m,b_m]} |u_m(x; V_m)-U_m(x)|\right] \leq \bbe[L]\Delta. \qquad
	\label{eq:Unfm-app-u-EU}
	\edeqn
	Let $\Delta$ be sufficiently small such that  $\bbe[L]\Delta\leq \epsilon/2$. 
	By (\ref{eq:u_mEU_m}) and (\ref{eq:Unfm-app-u-EU}), we have
	\bgeq
	\sup_{x\in\fX} 	|u(x; \bar{V}) - \bbe[U(x)]|
	&\leq& 
		\sup_{x\in\fX} |u(x; \bar{V}) - \bbe[u(x; V)]| +\bbe[L]\Delta\\
	&\leq& 	\sup_{x\in\fX} |u(x; \bar{V}) - \bbe[u(x; V)]| +\frac{\epsilon}{2}.
	\edeq
	Together with
	Proposition \ref{pp:equ_models}, we obtain
	\bgeq
	\prob \left(
	\sup_{x\in\fX} |u(x; \bar{V}) - \bbe[U(x)]| \geq \epsilon \right)
	&\leq&
	\prob\left(\sup_{x\in 
		\fX
	} |u(x; \bar{V}) - \bbe[u(x; V)]| \geq \frac{\epsilon}{2}\right)\\
	&\leq& \prob\left(\| \bar{V}-\bbe[V]\|\sup_{x\in
		\fX
	}\|f(x)\| \geq \frac{\epsilon}{2}\right).
	\edeq
	By the definition of $f$ in Proposition \ref{pp:equ_models}, we  know that $\sup_{x\in 
		\fX
	}\|f(x)\| <\infty$.
	Thus by 
	Cram\'er's large deviation theorem,
 there exists a positive integer $N_0$ and positive constant 
	$\Upsilon(\epsilon)$ (depending on $\epsilon$ with $\Upsilon(0)=0$)
	 such that for all $N\geq N_0$
	\bgeqn
	\prob\left(\|\bar{V}-\bbe[V]\|\sup_{x\in 
		\fX
	} \|f(x)\| \geq \frac{\epsilon}{2}\right)\leq e^{-\Upsilon(\epsilon)N}.
	\label{eq:VNm-appr-EVm-CLDP}
	\edeqn
Specifically, it follows by
 Hoeffding's inequality (see e.g.~\cite[Theorem 2]{shawe2003estimating}), 	
for fixed $\delta\in (0,1)$, we can set $N_0(\epsilon,\delta) := 
-\frac{\ln \delta}{\epsilon^2}$
and subsequently obtain
	\bgeqn
	\prob\left(\|\bar{V}-\bbe[V]\|\sup_{x\in 
		\fX
	} \|f(x)\| \geq \frac{\epsilon}{2}\right)\leq \delta
	\label{eq:VNm-appr-EVm}
	\edeqn
	for all $N\geq N_0(\epsilon,\delta)$.
 	The proof is complete. \hfill $\Box$

	With the theoretical 
justification of approximation of $\bbe[U(x)]$ with $u(x; \bar{V})$,
	we consider the sample data-based utility maximization problem
	\begin{align} \label{mod:PRO-SAA-1}
		\max_{x \in \fX}
		u(x; \bar{V}).
	\end{align}
Differing from standard sample average approximation
scheme in stochastic programming, 
problem (\ref{mod:PRO-SAA-1}) consists of two layers
of approximation: piecewise linear approximation
of utility function $U(x)$ by $u(x;V)$ 
and sample average approximation of 
$\bbe[u(x;V)]$.
 Let $\vartheta_{N,I}$ denote the optimal value of problem \eqref{mod:PRO-SAA-1} and $X_{N,I}$ the set of the optimal solutions.
	Here the subscripts $N$ and $I$ are used to indicate that these 
	quantities depend on both the sample size $N$ and the number of breakpoints $I$.
	By
	combining Propositions \ref{P-utility-appr-PL} and  \ref{P-u-appro-N}, we can
	establish 
 convergence of $u(x; \bar V )$ to  $\bE[U(x)]$ 
	and associated optimal values and
 optimal solutions
	as both $N$ and $I$ go to infinity. 
The following proposition
 addresses this.

	\begin{proposition}
		\label{C-u-appro-PL-DRO}
		Assume the settings and conditions of Propositions \ref{P-utility-appr-PL} and  \ref{P-u-appro-N}.
		Then the following assertions hold.

  \begin{itemize}
\item[(i)] 
Let $I$ (the set of breakpoints) be fixed.  For any $\epsilon > 0$ and $\delta > 0$, there exists $N_0(\epsilon,\delta) =
O(\ln \delta/\epsilon^2)$ such that
		\begin{equation}
   \prob\left( |\vartheta_{N,I}-\vartheta|\geq \epsilon\right)\geq 
   \delta,
			\label{eq:utility-approx-N}
		\end{equation}
		for all $N\geq N_0(\epsilon,\delta)$ 
and $\Delta\leq \frac{\epsilon}{2\bE[L]}$.

\item[(ii)]
  For any positive number $\epsilon$,
		\begin{equation}
			\lim_{ N,I \to \infty} \prob\left( |\vartheta_{N,I}-\vartheta|\geq \epsilon\right)
			\label{eq:sampleu-convg-trueu}
		\end{equation}
		and
		\begin{equation}
			\displaystyle{\lim_{N, I\to\infty}
				\prob(\mathbb{D}(X_{N,I}, X^*)\geq \epsilon)=0.}
			\label{eq:optimal-slu-SAA}
		\end{equation}
			\end{itemize}

 \end{proposition}
	
	\textbf{Proof.} 
Observe that
 \bgeqn
 |\vartheta_{N,I}-\vartheta|
			\leq  \sup_{x \in 
				\fX
			} |u(x; \bar V ) - \bE[U(x)] |.
\label{eq:vt-vt_IN-sup-u}
 \edeqn
 Part (i) follows directly from 
(\ref{eq:utility-approx-N-0}) and  (\ref{eq:vt-vt_IN-sup-u}).
 Part (ii). 
 Equality (\ref{eq:sampleu-convg-trueu})
	follows from (\ref{eq:utility-approx}), (\ref{eq:utility-approx-N}) 
 and (\ref{eq:vt-vt_IN-sup-u})
 whereas 
	(\ref{eq:optimal-slu-SAA}) follows from (\ref{eq:sampleu-convg-trueu}) and classical stability results (see e.g. \cite{BoS00}).
	\hfill $\Box$

The proposition provides a theoretical guarantee that 
one can solve problem~(\ref{mod:PRO-SAA-1})
to 
obtain an approximation
of the optimal value and optimal solution with specified 
confidence provided that the sample size is sufficiently large.

	\subsection{
	\ref{DUPRO}
	approximation and convergence}
	\label{subsec:The distributionally robust model}
	The convergence results established in Proposition \ref{C-u-appro-PL-DRO} are 
 based on the assumption
	that the sample size $N$ can be arbitrarily large. In some practical data-driven problems, 
	this assumption may not be fulfilled. This motivates us 
  to use \ref{DUPRO} with a piecewise linear random utility function to approximate problem (\ref{mod:PRO-general}).
  In Section~\ref{subsec:Construction of ambiguity set}, we propose two approaches to 
  construct an ambiguity set of  \ref{DUPRO}:
  an ellipsoid moment region 
  $\mathcal{P}_A$
  in \eqref{set:fP(bar V, S, gamma)} 
  and a bootstrap confidence region 
  $\mathcal{P}_B$
  in \eqref{set:fP_B}. 
  
  We begin with 
  $\mathcal{P}_A$
  based  \ref{DUPRO} 
  and move on to the one based on 
  $\mathcal{P}_B$.
   To  facilitate reading, we repeat the definition of $\mathcal{P}_A(\gamma_N)$
		\begin{align} \label{set:fP_again}
	\mathcal{P}_A(\gamma_N)
		:=
		\left\{
		P \in 
		\mathscr{P}(\mathbb{R}^I)
		\left |
		\begin{array}{l}
			P (V \in \mathcal{V}
			) = 1  \\
			\left\| S_{I-1}^{-1/2}  \left(C \bE_P [V] - \bar V_{I-1} \right)  \right\|^2 \le \gamma_N
		\end{array}
		\right.
		\right\}.
	\end{align}
By Proposition \ref{pp:equ_models},
we can 
reformulate
\ref{DUPRO} with 
$\mathcal{P}_A(\gamma_N)$
as 	\begin{align} \label{mod:PRO-SAA-DRO-1}
	\vt_{\text{elp}}:=	\max_{x \in \fX}
		\min_{	v\in \fF(
		\mathcal{P}_A(\gamma_N)
		)}
		u(x; v),
	\end{align}
where
			\bgeqn
		\fF(\mathcal{P}_A(\gamma_N)
		) :=\left\{v\in \mathcal{V}
		\;|\;
		(Cv - \bar{V}_{I-1})^T S_{I-1}^{-1}
		(Cv - 
		\bar{V}_{I-1}) \le \gamma_N \right\}.
		\label{eq:set-fV}
		\edeqn
  The reformulation is sensible 
in that the objective function only depends on
$\bbe[V]$ and so does the moment constraint
in the ambiguity set.
Analogous to 
$\gamma_N^{(1)}$ defined as in (\ref{eq:gamma_N^2}), we
require that $\gamma_N$ decreases to 0 as $N$
 goes to $\infty$.
 Then
	the ambiguity set $\fF(
	\mathcal{P}_A(\gamma_N)
	)$ converges to the singleton 
 $\{ \mu = \bE[V] \}$.
We investigate 
convergence of 
$\vt_{\text{elp}}$
to $\vt$ (the optimal value of problem (\ref{mod:PRO-general}))
as $N$ and $I$ go to $\infty$ and  $\gamma_N$ is driven to 0. The next theorem states the convergence. 

	\begin{theorem}[\ref{DUPRO} with ellipsoidal ambiguity (\ref{eq:set-fV})]
	\label{thm:DUPRO-Ellip-optim}
	Let $\gamma_N$ in \eqref{set:fP_again} monotonically decrease to $0$ as $N \to \infty$. 
	Then for any 
   $\alpha\in (0,1)$ and
	small $\epsilon > 0$, there exist positive numbers  $I_0(\epsilon), N_0(\epsilon,\alpha)$  such that
		\bgeqn
		\label{eq:convg-opti-value-thm2}
		\prob\left(
		|	\vt_{{\rm elp}}
		-\vt|\geq \epsilon
		\right)
		\leq \alpha
		\edeqn
		for all 
 $N\geq N_0(\epsilon,\alpha)
= 
O\left((\ln \alpha)/\epsilon^2\right)$,\, 
$I\geq I_0(\epsilon)= O(1/\epsilon)$.
	\end{theorem}

	\proof
For fixed $v\in \fF(
\mathcal{P}_A(\gamma_N)
	)$, 
	it follows by (\ref{eq:u-vec-rep}) that 
	\bgeqn
	|u(x;v)  -\bbe[U(x)]|
	&\leq &
	|u(x;v) - u(x;\bar{V})| +
	|u(x; \bar{V}) - \bbe[U(x)]|\nonumber\\
	&=&
		|f(x)^T(v -\bar{V})| + |u(x; \bar{V}) - \bbe[U(x)]|.
\label{eq:ineqlty-u-appr}
	\edeqn
Let 
 $E$ denote the event that 
$
S_{I-1}^{-1}\succ
(2\Sigma_{I-1})^{-1}$
and $F$ denote the event that 
$\sup_{x\in \fX}|u(x; \bar{V}) - \bbe[U(x)]|\leq \epsilon/2$.
Since $
2\Sigma_{I-1} \succ \Sigma_{I-1}$,
by Corollary 6 in \cite{shawe2003estimating}, 
there exists 
$\hat{N}_0(\alpha)=
O(-\ln\alpha)$ 
such that
\bgeqn 
\prob(E) =\prob\left(S_{I-1}^{-1}\succ
(2\Sigma_{I-1})^{-1}\right)\geq 1-\alpha
\label{eq:E}
\edeqn
for all 
$N\geq \hat{N}_0(\alpha)$.
Moreover, it follows by Proposition~\ref{P-u-appro-N} and its proof,
there exists 
$\tilde{N}_0(\epsilon/2,\alpha)
=O((\ln \alpha)/\epsilon^2)
$
such that 
\bgeqn 
\prob(F)=\left(\sup_{x\in \fX}|u(x; \bar{V}) - \bbe[U(x)]|\leq \epsilon/2\right) \geq \alpha
\label{eq:F}
\edeqn 
for all 
$N\geq 
\tilde{N}_0(\epsilon/2,\alpha)$.
Thus 
\bgeqn 
\prob(E\cap F) \geq 1-\alpha
\label{eq:EF}
\edeqn 
for all 
 $N\geq N_0(\epsilon,\alpha):=
 \max\{\hat{N}_0(\alpha),\tilde{N}_0(\epsilon/2,\alpha)\}
=O((\ln \alpha)/\epsilon^2)
 $.
			By (\ref{eq:ineqlty-u-appr}),
	\bgeqn
		& & \prob\left(
		|\vt_{{\rm elp}}-\vt|\geq \epsilon
		\right) \nonumber \\
	&\leq&
	\prob\left(
	\sup_{x\in 
		\fX, v \in 	\fF \left(
		\mathcal{P}_A
		\left(\gamma_N \right) \right)}
	|u(x;v)  -\bbe[U(x)]|\geq \epsilon\right)\nonumber\\
&=&	\prob\left(
	\sup_{x\in 
		\fX, v \in 	\fF\left(
	\mathcal{P}_A
		\left(\gamma_N \right) \right)}
	|u(x;v)  -\bbe[U(x)]|\geq \epsilon\bigcap (E\cap F)\right)\nonumber\\
	&& \qquad +
	\prob\left(
	\sup_{x\in 
		\fX, v \in 	\fF\left(
	\mathcal{P}_A
		\left(\gamma_N \right) \right)}
	|u(x;v)  -\bbe[U(x)]|\geq \epsilon\bigcap \overline{E\cap F}\right)\nonumber\\
&\leq& \prob\left(
	\sup_{x\in 
		\fX, v \in 	\fF\left(
	\mathcal{P}_A
		\left(\gamma_N \right) \right)}
	|u(x;v) - u(x,\bar{V})| +
	|u(x; \bar{V}) - \bbe[U(x)]| \geq \epsilon \;\Bigg| E\cap F\right)\prob(E\cap F)\nonumber\\
	&& \qquad + \prob(\overline{E\cap F})\nonumber\\
&\leq& \prob\left(
	\sup_{x\in 
		\fX, v \in 	\fF\left(
	\mathcal{P}_A
		\left(\gamma_N \right) \right)}
	|u(x;v) - u(x;\bar{V})| \geq \frac{\epsilon}{2}\right)\prob(E\cap F)+
\prob(\overline{E\cap F})
 \quad (\inmat{By} \; (\ref{eq:F})) 
\nonumber\\
	&\leq &
	\prob\left(\sup_{x\in 
		\fX, v \in 	\fF\left(
	\mathcal{P}_A
		\left(\gamma_N \right) \right)
	}|f(x)^T(v -\bar{V})|\geq \frac{\epsilon}{2} \right)+\alpha. \quad (\inmat{By} \; (\ref{eq:EF})) 
	\label{eq:Thm2-proof-1}
	\edeqn
Let $f_{I-1}(x)$, $v_{I-1}$ and $\bar{V}_{I-1}$ denote respectively the 
vectors which consist of the first $I-1$ components of $f(x), v$ and $\bar{V}$, 
let $f_I(x), v_I, \bar{V}_I$ denote the $I$-th component of $f(x)$, $v$ and $\bar{V}$, and $e_{I-1}$
denotes an $(I-1)$-dimensional vector with all ones.
Then  
$$
v_I-\bar{V}_I = 1-e^T_{I-1}v_{I-1} - (1-e^T_{I-1}\bar{V}_{I-1})=e_{I-1}^T(\bar{V}_{I-1} - v_{I-1})
$$
and subsequently
\bgeq 
f(x)^T(v -\bar{V})= f_{I-1}(x)^T( v_{I-1}-\bar{V}_{I-1}) +
f_I(x) e_{I-1}^T(\bar{V}_{I-1} - v_{I-1}).
\edeq
By the H\"older inequality
\bgeqn 
|f(x)^T(v -\bar{V})|
&\leq&
\|
(2\Sigma)_{I-1}^{1/2}f_{I-1}(x)\|
\|
(2\Sigma)_{I-1}^{-1/2}( v_{I-1}-\bar{V}_{I-1})\|\nonumber\\ 
&&+
|f_I(x)|\|
(2\Sigma)_{I-1}^{1/2}e_{I-1}\|
\|
(2\Sigma)_{I-1}^{-1/2}(v_{I-1}-\bar{V}_{I-1})\|\nonumber\\
&\leq & (I-1)\|
(2\Sigma)_{I-1}^{1/2}\|\|f(x)\|\|
(2\Sigma)_{I-1}^{-1/2}( v_{I-1}-\bar{V}_{I-1})\|.
\label{eq:v_I-and-V}
\edeqn
Let 
$N\geq N_0(\epsilon,\alpha)$
be sufficiently large such that
\bgeqn
\frac{\epsilon}{2
(I-1)\|
(2\Sigma)_{I-1}^{1/2}\|\sup_{x\in\fX}\|f(x)\|} 
> 
\sqrt{\gamma_N}
\label{eq:Thm2-proof-Opti-value-4}
\edeqn
for all 
$N\geq N_0(\epsilon,\alpha)$.
Combining (\ref{eq:Thm2-proof-1})-(\ref{eq:Thm2-proof-Opti-value-4}),
we have
\bgeq
&&
	\prob\left(
		|
			\vt_{\text{elp}}
		-\vt|\geq \epsilon
		\right)\nonumber\\
	&\leq&
\prob\left(\sup_{x\in 
		\fX, v \in 	\fF(
	\mathcal{P}_A
		(\gamma_N))
	}|f(x)^T(v -\bar{V})|\geq \frac{\epsilon}{2} \right)+\alpha\nonumber\\
&\leq&
\prob\left(\sup_{x\in 
		\fX, v \in 	\fF(
	\mathcal{P}_A
		(\gamma_N)
		)
	}
	(I-1)\|
	(2\Sigma)_{I-1}^{1/2}\|\|f(x)\|\|
	(2\Sigma)_{I-1}^{-1/2}( v_{I-1}-\bar{V}_{I-1})\|
	\geq \frac{\epsilon}{2}\right)+\alpha\nonumber\\
&\leq& 
\prob\left(\sup_{v \in 	\fF(
\mathcal{P}_A
(\gamma_N)}\|
(2\Sigma)_{I-1}^{-1/2}( v_{I-1}-\bar{V}_{I-1})\|
	\geq \frac{\epsilon}{2\sup_{x\in 
		\fX 
	}
	(I-1)\|
	(2\Sigma)_{I-1}^{1/2}\|\|f(x)\|}\right)+\alpha
\nonumber\\
& \le & 
\prob\left(\sup_{v \in 	\fF(
\mathcal{P}_A
(\gamma_N))}\|
(2\Sigma)_{I-1}^{-1/2}( v_{I-1}-\bar{V}_{I-1})\|
	> \sqrt{\gamma_N}\right)+\alpha
\nonumber\\
&\leq& 
\prob\left(\sup_{v \in 	\fF(
\mathcal{P}_A
(\gamma_N))}
		(Cv - \bar{V}_{I-1})^T 
		S_{I-1}^{-1}
		(Cv - 
		\bar{V}_{I-1}) > \gamma_N\right)+\alpha\\
						&=&0+\alpha=\alpha,
		 \edeq
%
where the  last inequality is 
 due to (\ref{eq:E})
and the last equality is because of (\ref{eq:set-fV}).
\hfill $\Box$

The theorem provides a theoretical guarantee that
the optimal value obtained from solving 
problem (\ref{mod:PRO-SAA-DRO-1}) 
approximates the true one
with a specified confidence when the sample size
goes to infinity and $\gamma_N\to 0$. This kind of convergence result should be distinguished 
from those in the literature
of 
DRO
(see e.g.~\cite{sun2016convergence}) in that here problem 
(\ref{mod:PRO-SAA-DRO-1}) is a robust optimization model instead of a 
DRO
model, and the ambiguity set $\fF(\mathcal{P}_A(\gamma_N)$ in terms of $v$
in (\ref{eq:set-fV}) shrinks to a singleton whereas 
the ambiguity set in terms of $P$ in (\ref{set:fP_again})
does not. 
One of the main challenges that we have to tackle in the proof is that the ellipsoid 
in (\ref{eq:set-fV}) is defined for 
$v_{I-1}$ rather than $v$ whereas the objective function $u(x;v)$ depends on $v$.
The mismatch requires us
to derive an upper bound for 
$|f(x)^T(v -\bar{V})|$ in terms of 
$S_{I-1}^{-1/2}(v_{I-1} - \bar{V}_{I-1}) 
		$.

We now move on to investigate 
convergence of 
the optimal value of \ref{DUPRO}
to $\vt$ when the ambiguity set is constructed by the bootstrap samples.
Recall that $\fW^*_{1-\alpha}$ defined as in \eqref{eq:set-fW-Boots}
is the convex hull of points 
$T^*_{(1)}, \dots, T^*_{(\lceil (1-\alpha) K \rceil)}$.
Here we consider the case that $K=\infty$ and use 
$\fW_{1-\alpha}$ to denote 
the $100(1-\alpha)\%$ interior
points of the convex hull of 
$T^*_{(1)}, \dots, T^*_{(\infty)}$
based on Tukey's depth.
Consider \ref{DUPRO} with
\begin{align} \label{set:fP_hat B}
\mathcal{P}_{\hat{B}}
	(\alpha) =
	\left\{
	P \in \mathscr{P}(\mathbb{R}^I) \ |\ P(V \in \mathcal{V}
	) = 1, \ C \bE_P [V] =  \bar V_{I-1} -  S_{I-1}^{1/2} \tilde w / \sqrt{N},  \tilde w \in   \fW_{1-\alpha}
	\right\}.
\end{align}
The difference between 
$\mathcal{P}_{\hat{B}}(\alpha)$
and 
$\mathcal{P}_B(\alpha)$
in \eqref{set:fP_B} is that 
$\mathcal{P}_B(\alpha)$
is 
defined with 
$\fW^*_{1-\alpha}$ while 
$\mathcal{P}_{\hat{B}}(\alpha)$
is based on $\fW_{1-\alpha}$. 
The next lemma refers to Theorem 1 in \cite{yeh_balanced_1997} and the
comments following the theorem. 

\begin{lemma} \label{lem:bootstrap}
	If $P$, the true probability measure of $V$, is absolutely continuous in $\bR^I$, then there exists $\delta_{1-\alpha} > 0$ depending on $\alpha$ such that
	\begin{align}
 \label{eq:YeHS}
		\lim_{N \to \infty} \fW_{1-\alpha} \subseteq \fB(\delta_{1-\alpha}), \ \text{a.s.,}
	\end{align}
	where $\fB(\delta_{1-\alpha})$ is an $I$-dimensional
	ball 
	centered at $0$ with radius $\delta_{1-\alpha}$.
\end{lemma}

The lemma indicates that
$\fF(\mathcal{P}_{\hat{B}}(\alpha))$
is contained in an ellipsoid 
center at $\bar{V}_{I-1}$ when $N$ goes to infinity.
We will use this fact 
and Theorem \ref{thm:DUPRO-Ellip-optim}
to 
establish the convergence property of \ref{DUPRO} with $\mathcal{P}_{\hat{B}}(\alpha)$
in the 
next theorem.

	\begin{theorem}[\ref{DUPRO} with bootstrap ambiguity \eqref{set:fP_hat B}]
	\label{thm:DUPRO-Bootstrap-convg}
	Let
\begin{align}
\label{eq:minmax-BooTS}
	\vt_{{\rm bts}}:=	\max_{x \in \fX}
		\min_{	v\in \fF(
	\mathcal{P}_{\hat{B}}(\alpha)
		)}
		u(x; v)
	\end{align}
	and
	 $\tau \in (0,1)$,
  where 
	\begin{align}
	\label{eq:fF_PB_alpha-infinite-K}
		\fF(
		\mathcal{P}_{\hat{B}}(\alpha)
		)  = \left \{
		v \in \bR^I
		\left |
		\begin{array}{l}
			v \ge 0 \\
			e_I^Tv = 1 \\
			Cv=\bar V_{I-1}-S^{1/2}_{I-1}\tilde{w}/\sqrt{N},\; \tilde w \in \fW_{1-\alpha} 
		\end{array}
		\right.
		\right \}.
	\end{align}
    Suppose that $P$, the true probability 
	 measure 
	 induced by $V$, is absolutely continuous in $\bR^I$. 
	Then for any 
	small positive number $\epsilon$, there exist positive numbers $I_0$, $\hat{N}_0$  such that
		\bgeqn
		\label{eq:convg-opti-value-thm3}
		\prob\left(
		|\vt_{{\rm bts}}-\vt|\geq \epsilon
		\right)
		\leq \tau
		\edeqn
		for all $N\geq \hat{N}_0$, $I\geq I_0$.
	\end{theorem}

It is important to highlight that (\ref{eq:minmax-BooTS})
is a maximin robust optimization problem rather than a 
DRO
problem.
Moreover, 
by Lemma \ref{lem:bootstrap},
$\fW_{1-\alpha}$ is a bounded set a.s. as $N$ goes to infinity
and 
by Corollary 6 in \cite{shawe2003estimating},
$S^{1/2}_{I-1}$ converges to $\Sigma^{1/2}_{I-1}$ 
at an exponential rate
with the increase of $N$, then 
$S^{1/2}_{I-1}\tilde{w}/\sqrt{N}\to 0$ as $N\to \infty$ w.p.1. This means the ambiguity set
$\fF(
		\mathcal{P}_{\hat{B}}(\alpha)
		) $ 
  shrinks to a singleton w.p.1 as $N\to\infty$.

	\proof
The thrust of the proof is to
 show 
\bgeqn 
\label{eq:P_B--P_A}
\fF(
\mathcal{P}
_{\hat B}(\alpha)) \subseteq \fF(
\mathcal{P}
_A (\gamma^{(2)}_N ))
\edeqn 
for some $\gamma^{(2)}_N$
when $N$ is sufficiently large 
and then the conclusion follows from a similar proof to that of Theorem \ref{thm:DUPRO-Ellip-optim}.
To this end, we 
choose $\eta \in (0, \tau)$ and let $\beta = 1 - (\tau - \eta)$.
It follows by Lemma \ref{lem:bootstrap} that there exist $\delta_{1-\alpha} >0 $ and $N_1 > 0$ such that
\begin{align}
	\prob (\fW_{1-\alpha} \subseteq \fB(\delta_{1-\alpha})) \ge \beta,
\end{align}
for all $N > N_1$. Let 
 $\gamma^{(2)}_N :=\delta_{1-\alpha}^2/N$. 
Note that
\begin{align*}
	\fF(
\mathcal{P}_{\hat{B}}
	(\alpha)) = \{v \in \mathcal{V}
	\;|\; S^{-1/2}_{I-1} (Cv - \bar V_{I-1}) =  
	-\tilde{w}/\sqrt{N}, \tilde w \in \fW_{1-\alpha}\}
\end{align*}
and
\begin{align*}
	\fF(
\mathcal{P}_A
	(\gamma^{(2)}_N ))) = \{v \in \mathcal{V}
	\;|\; S^{-1/2}_{I-1} (Cv - \bar V_{I-1})= 
	-\tilde{w}/\sqrt{N}, \tilde w \in \fB(\delta_{1-\alpha}) \}.
\end{align*}
Let G denote the event that
(\ref{eq:P_B--P_A}) holds.
We then obtain
\bgeqn 
\label{eq:P_B--P_A-prob}
\prob\left(\fF(
\mathcal{P}
_{\hat B}(\alpha)) \subseteq \fF(
\mathcal{P}
_A (\gamma^{(2)}_N ))\right)\geq \beta
\edeqn
for $N > N_1$.
Consequently, for $N > N_1$, we have
	\begin{align}
		\prob\left(
		|	\vt_{\text{bts}}
		-\vt|\geq \epsilon
		\right) 
	\leq \ &
	\prob\left(
	\sup_{x\in 
		\fX, 
		v \in  \fF(
	\mathcal{P}
		_{\hat B} (\alpha))}
	|u(x;v)  -\bE [U(x)]|\geq \epsilon\right) \nonumber\\
	\leq \ & 	\prob\left(
	\sup_{x\in 
		\fX, 
		v \in  \fF(
	\mathcal{P}
		_{\hat B} (\alpha))}
	|u(x;v)  -\bE [U(x)]|\geq \epsilon \;\Bigg |\; G \right) \prob ( G )
 \nonumber\\
	& \qquad  + 	\prob\left(
	\sup_{x\in 
		\fX, 
		v \in  \fF(
		\mathcal{P}
		_{\hat B} (\alpha))}
	|u(x;v)  -\bE [U(x)]|\geq \epsilon \;\Bigg |\; \overline G \right) \prob ( \overline G )
 \nonumber\\
	\leq \ & 
	\prob\left(
	\sup_{x\in 
		\fX, v \in 	 \fF(
		\mathcal{P}
		_A(\gamma^{(2)}_N))}
	|u(x;v)  -\bbe[U(x)]|\geq \epsilon\right) + \prob ( \overline G ).
 \label{eq:boots-exponential-convg-a}
	\end{align}
It follows from a similar analysis to the proof of Theorem \ref{thm:DUPRO-Ellip-optim} that,
there exist
$N_2$ and $I_0$ 
such that
\begin{align}
\label{eq:boots-exponential-convg}
	\prob\left(
	\sup_{x\in 
		\fX, v \in 	 \fF(
	\mathcal{P}
		_A(\gamma^{(2)}_N))}
	|u(x;v)  -\bbe[U(x)]|\geq \epsilon\right) \le \eta,
\end{align}
for $N > N_2$ and $I > I_0$. 
Let $\hat{N}_0=\max\{N_1,N_2\}$.  Then for $N\geq \hat{N}_0$, we have
from (\ref{eq:P_B--P_A-prob})-(\ref{eq:boots-exponential-convg})
that
	\begin{align*}
	\prob\left(
	|
		\vt_{\text{bts}}-\vt|\geq \epsilon
	\right) 
	\leq  
 \eta + (1 - \beta) = \tau.
	\end{align*}
The proof is complete. \hfill $\Box$	

\begin{figure}[!htbp]
\vspace{-0.2cm}
\minipage{0.33\textwidth}
 \centering
 \includegraphics[width=4.5cm]{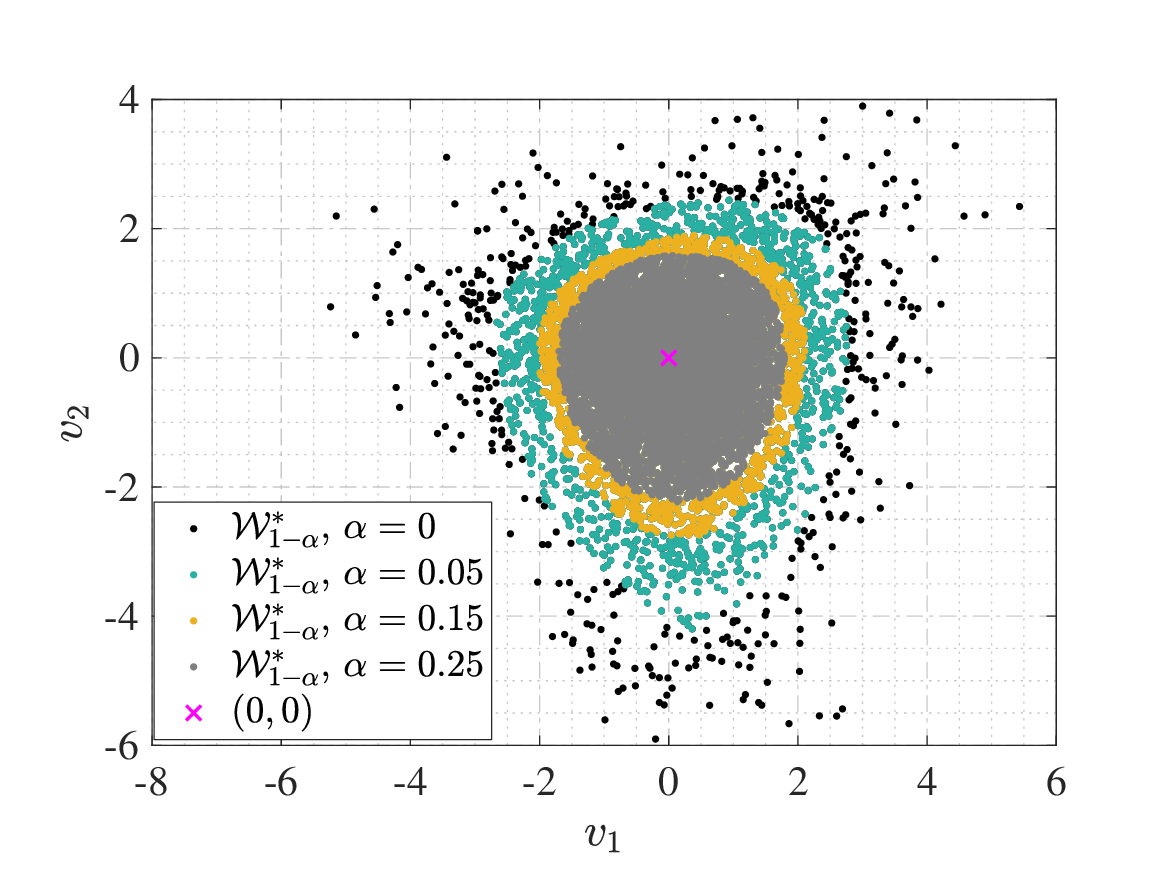}
  \textbf{\tiny (a) $N=20$}
\endminipage\hfill
\minipage{0.33\textwidth}
  \centering
 \includegraphics[width=4.5cm]{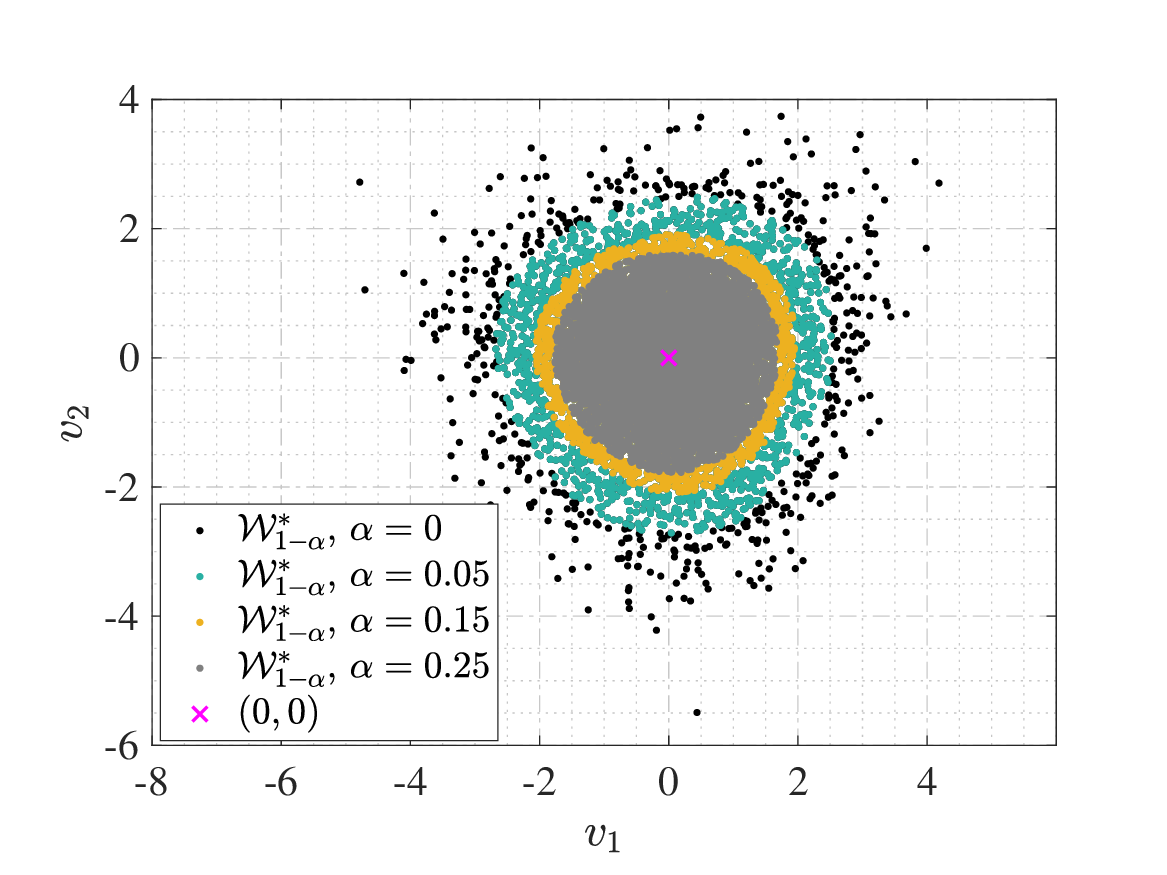}
  \textbf{\tiny (b) $N=50$}
\endminipage
\hfill
\minipage{0.33\textwidth}
  \centering
 \includegraphics[width=4.5cm]{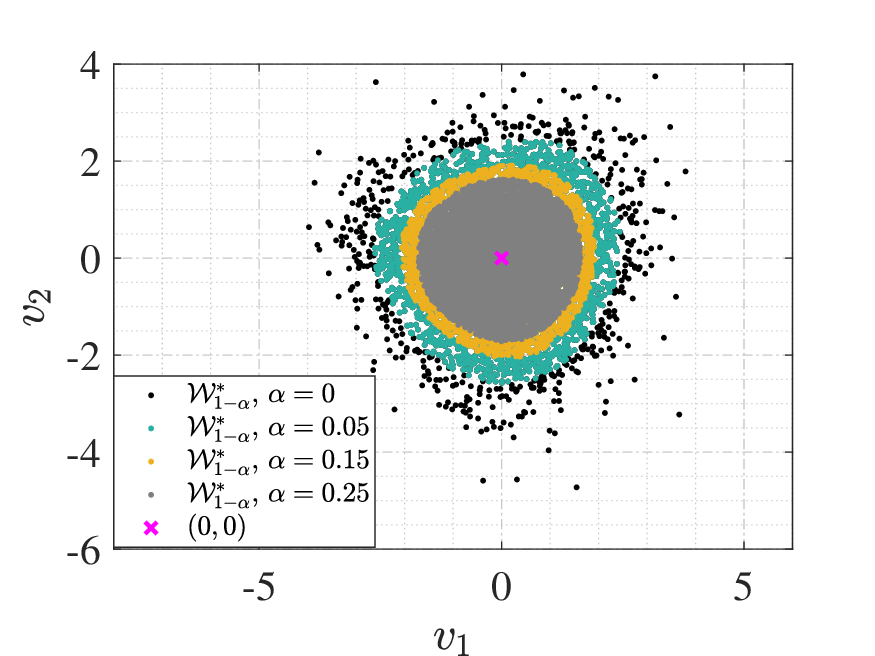}
  \textbf{\tiny (c) $N=100$}
\endminipage
\caption{\small $N=20$, $N=50$, $N=100$,  $K=10,000$.
The sets of coloured points represent
the set 
$\fW^*_{1-\alpha}$ 
with different values of $\alpha$.}
\label{Cstar_vs_ellipsoid_2}
\vspace{-0.5cm}
\end{figure}

\begin{figure}[!htbp]
\minipage{0.33\textwidth}
 \centering
 \includegraphics[width=4.5cm]{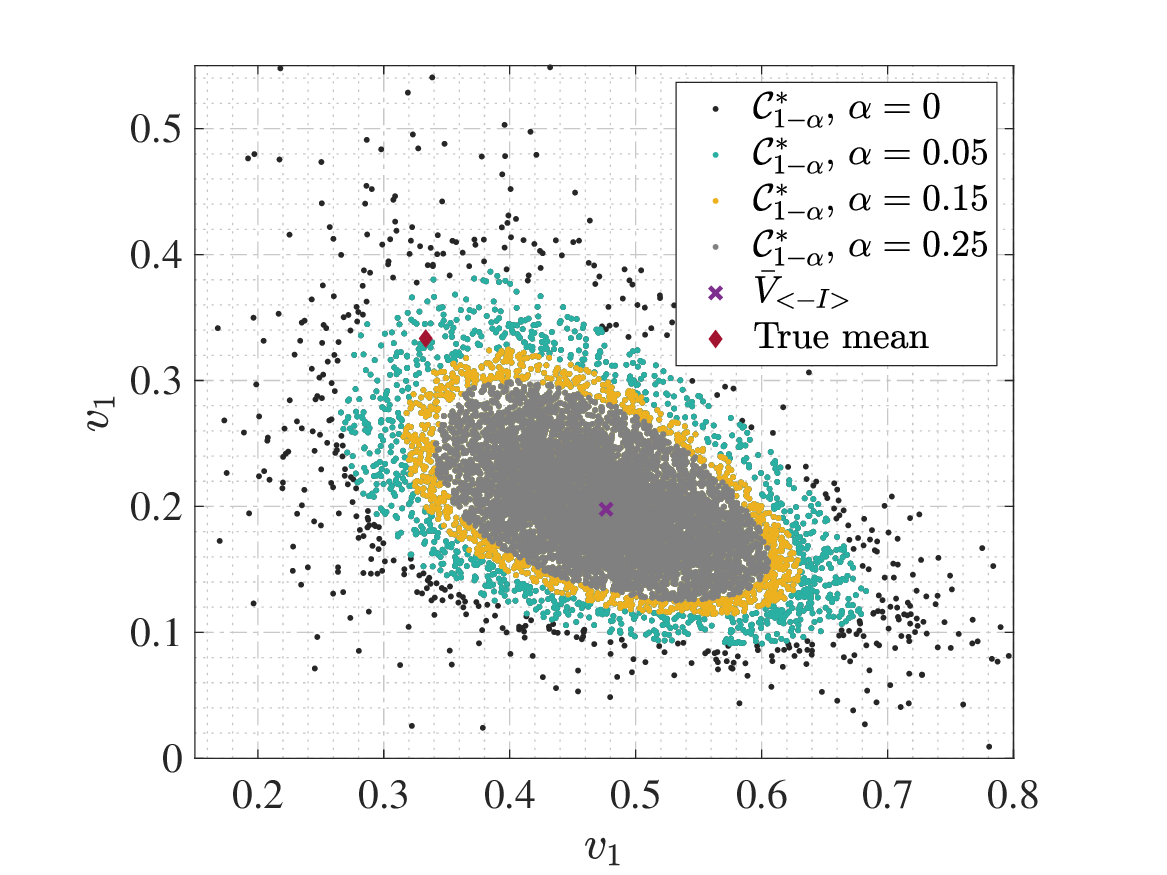}
  \textbf{\tiny (a) $N=20$}
\endminipage\hfill
\minipage{0.33\textwidth}
  \centering
 \includegraphics[width=4.5cm]{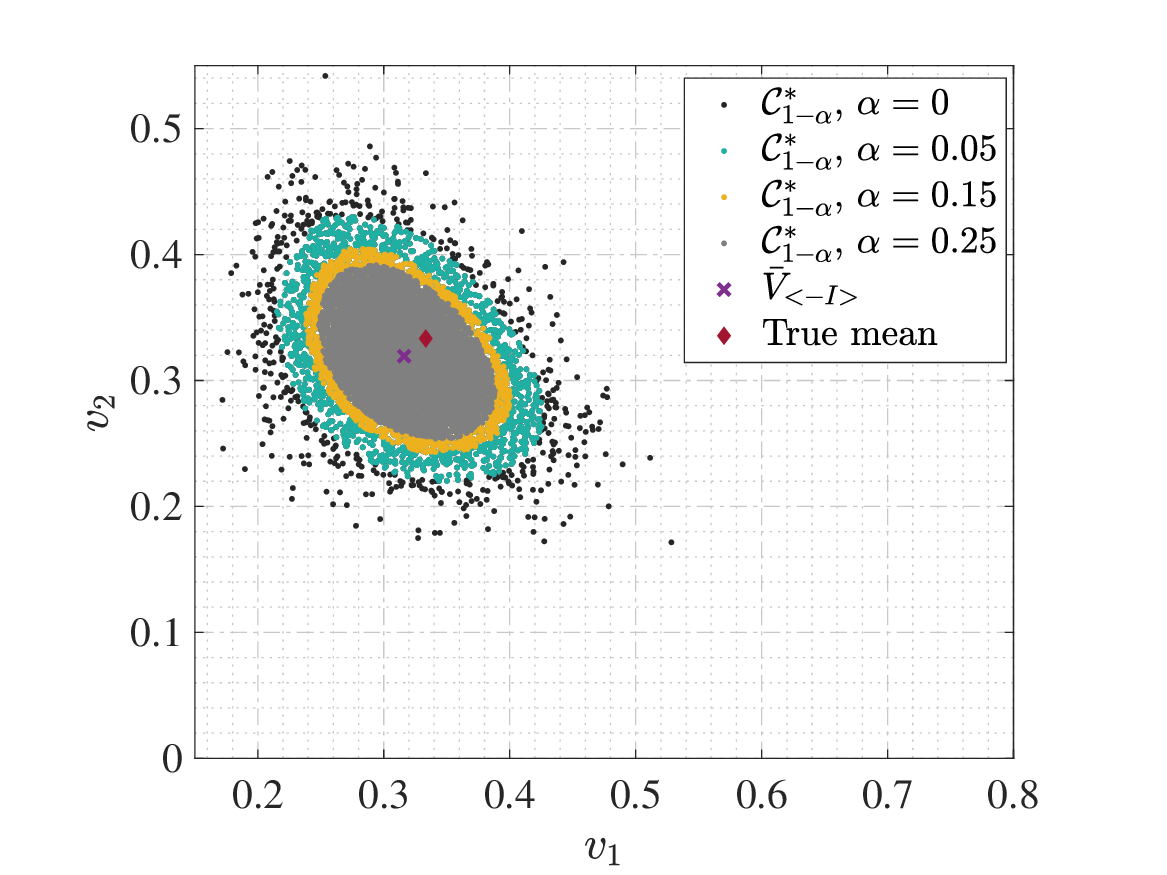}
  \textbf{\tiny (b) $N=50$}
\endminipage
\hfill
\minipage{0.33\textwidth}
  \centering
 \includegraphics[width=4.5cm]{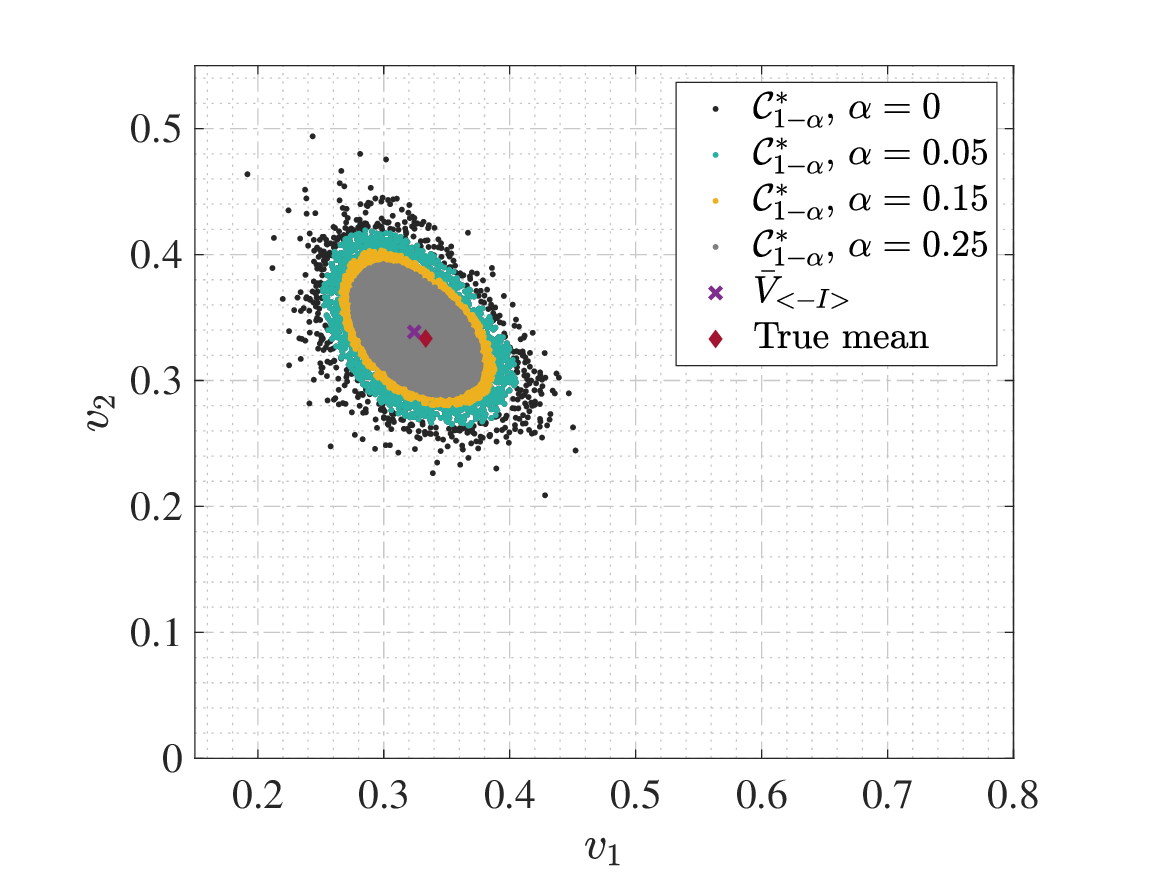}
  \textbf{\tiny (c) $N=100$}
\endminipage
\caption{\small $N=20$, $N=50$, $N=100$,  $K=10,000$.
The sets of coloured points represent
the set 
$\fC^*_{1-\alpha}$
with different values of $\alpha$.}
\label{Cstar_vs_ellipsoid_3}
\vspace{-0.5cm}
\end{figure}

The theorem guarantees that the optimal value obtained from solving 
problem (\ref{eq:minmax-BooTS}) converges to the true 
one
with a specified confidence as
the sample size $N$ of the original sample and the number $K$ of bootstrap resampling go to infinity.
The 
constant $\hat{N}_0$ 
depends heavily on $N_1$
to ensure (\ref{eq:P_B--P_A-prob}), and this is essentially down to 
the rate of convergence in (\ref{eq:YeHS}). In the literature of statistics, 
it is shown that the empirical cumulative distribution function 
of  statistical estimator based on bootstrap resamples 
converges to the one based on the original sample at a rate of 
$o(N^{-1/2})$
in single variate case, see Theorem 3.10 and follow-up discussions in \cite{dikta2021bootstrap}.
The situation here is much more complex since here $v_{I-1}$ is multivariate and 
the ambiguity set is constructed via Tukey's depth which trims down some of the resamples. It remains an open challenging question as to whether the classical result may be generalized to this case, we leave this for future research.
However, we have done some numerical studies 
on the convergence of $\fW^*_{1-\alpha}$ as $N$ increases.
We further set $K=10,000$ and such a large $K$ guarantees that $\fW^*_{1-\alpha}$ can closely approximates $\fW_{1-\alpha}$. On this basis, we examine the convergence 
of $\fW^*_{1-\alpha}$ 
as $N$ increases 
from 20 to 50 and 100 for different $\alpha$ values. Figure \ref{Cstar_vs_ellipsoid_2}
displays the tendency of the convergence.
Figure \ref{Cstar_vs_ellipsoid_3} depicts the respective ambiguity sets $\fC^*_{1-\alpha}$ in \eqref{eq:set:fC}.
We can see that when original sample size $N$ reaches $50$, the true mean lies within the ambiguity set $\fC^*_{1-\alpha}$ with $\alpha=0.25$.

  	\section{
  	Numerical tests}
	\label{sec:Case Studies}
    We 
have carried out some numerical 
studies of 
the proposed \ref{DUPRO} model.
In this section, we report the test results.
We begin 
by examining the performance of \ref{DUPRO} 
with
bootstrap ambiguity set 
 where we look into the  effect of the sample size of the random parameter $V$, the number of the bootstrap resamples, and the critical value determining the size of $\mathcal{P}_B$ (see (\ref{set:fP_B}))
and then move on to comparative analysis 
of \ref{DUPRO} 
relative to the sample average approximation model.
 

	\subsection{
 Performance studies of the bootstrap-based \ref{DUPRO} model}
 	\label{case_study1}

The first set of tests are aimed to analyse 
the performances of \ref{DUPRO} with the ambiguity set 
$\mathcal{P}_B(\alpha)$ defined as in \eqref{set:fP_B}:
\begin{enumerate}[label=(\alph*)]
    \item \label{test:a} Report the optimal values,
    worst-case expected utility functions, and computational time of
\ref{DUPRO}
based on 
different settings of the critical value
$\alpha$.

\item \label{test:b} Investigate the impact of increase 
of
the original sample size $N$ in relation to 
the convergence result in Theorem~\ref{thm:DUPRO-Bootstrap-convg}.

\item \label{test:c} Validate the reliability and performance of \ref{DUPRO} using out-of-sample tests. 
\end{enumerate}
The test problem has three attributes with a feasible set 
$\fX := \{x \in \mathbb{R}^3_+ \;|\; e^T_3 x = 1 \}$. 
The set of breakpoints of the piecewise utility functions can be set arbitrarily and we fix them here,
see  
Table~\ref{tab:breakpoints_3attributes}.
   \begin{table}[ht!]
\vspace{-0.3cm}
 	\centering
 	\begin{tiny}
 		\begin{tabular}{|p{1.2in}<{\centering} |c|p{0.4in}<{\centering}|c|p{0.4in}<{\centering}|c|p{0.4in}<{\centering}|}
 			\hline
 			\rule{0pt}{8pt}
 			Attributes
 			&  {$t_{m,1}$} 
 			&  {$t_{m,2}$ } 
 			& {$t_{m,3}$} 
 			&  {$t_{m,4}$}  
 			&  {$t_{m,5}$} 
 			&  {$t_{m,6}$} 
 			\\
 			\hline
 			\hline
	\tabincell{c}{Attribute 1}
 			& 0.0667 
 			&  0.4
 			& 0.6667 
 			& 1 
 			&  -
 			&  -
 			\\
 			\hline
            \tabincell{c}{Attribute 2}
 			& 0.05 
 			& 0.3
 			& 0.4 
 			& 0.5 
 			& 0.75 
 			&  1
 			\\
 			\hline
            \tabincell{c}{Attribute 3}
            & 0.0667 
 			& 0.2
 			& 0.5333 
 			& 0.6667 
 			&  1
 			& -
 			\\
 			\hline
 			
 		\end{tabular}
 	\end{tiny}
  	\vspace{-0.1cm}
 	\caption{The breakpoints of the tree attributes.}
 	\label{tab:breakpoints_3attributes}
 	\vspace{-0.5cm}
 \end{table}
 Thus $I_1=4$, $I_2=6$ and $I_3=5$ and $I=15$.
 We assume that the true probability distribution of vector
 $V$ (of dimension $I=15$)
 follows Dirichlet distribution with
 order $15$ with parameters $(0.5,0.5,\cdots,0.5)\in\mathbb{R}^{15}$. Thus, $\bbe[V]=(1/15,\cdots,1/15)\in \R^{15}$.
In this subsection, we report our findings.
 All of the tests are carried out by Matlab 2022a installed on a PC (16GB RAM, CPU 2.3 GHz) with Intel Core i7 processor. 
 We 
 use Gurobi solver to 
 solve the mixed integer problem
 (\ref{pp:PRO_B_equ}).

\subsubsection{Test case \ref{test:a}: Optimal value,
worst-case expected utility function and computational time}
We generate a random sample, $V^1, \cdots,V^{50}$, with Dirichlet distribution 
and then generate
  $K=10,000$
  bootstrap resamples $V^{k,1},\cdots,V^{k,50}$, $k=1,\cdots, 10,000.$
Table~\ref{tab:sample_mean_3attributes} displays the sample mean of $V$.

 \begin{table}[ht!]
 \vspace{-0.3cm}
 	\centering
 	\begin{tiny}
 		\begin{tabular}{|p{1.2in}<{\centering} |c|p{0.4in}<{\centering}|c|p{0.4in}<{\centering}|c|p{0.4in}<{\centering}|}
 			\hline
 			\rule{0pt}{8pt}
 			Attributes
 			&  {$\bar{V}_{m,1}$} 
 			&  {$\bar{V}_{m,2}$ } 
 			& {$\bar{V}_{m,3}$} 
 			&  {$\bar{V}_{m,4}$}  
 			&  {$\bar{V}_{m,5}$} 
 			&  {$\bar{V}_{m,6}$} 
 			\\
 			\hline
 			\hline
 			\tabincell{c}{Attribute 1}
 			& 0.0591
 			& 0.0537
 			& 0.0642
 			& 0.0719 
 			&  -
 			&  -
 			\\
 			\hline
            \tabincell{c}{Attribute 2}
 			& 0.0514
 			& 0.0592
 			& 0.0756 
 			& 0.0551 
 			& 0.1052 
 			& 0.0686
 			\\
 			\hline
            \tabincell{c}{Attribute 3}
 			& 0.0706
 			& 0.0675 
 			& 0.0594 
 			& 0.0830
 			& 0.0555
 			&  -
 			\\
 			\hline
 		\end{tabular}
 	\end{tiny}
  	\vspace{-0.1cm}
 	\caption{The sample mean $\bar{V}_{m,k}$ at the breakpoints.}
 	\label{tab:sample_mean_3attributes}
 \vspace{-0.5cm}
 \end{table}

\begin{figure}[!htbp]
\minipage{0.5\textwidth}
 \centering
  \includegraphics[width=5.0cm]{{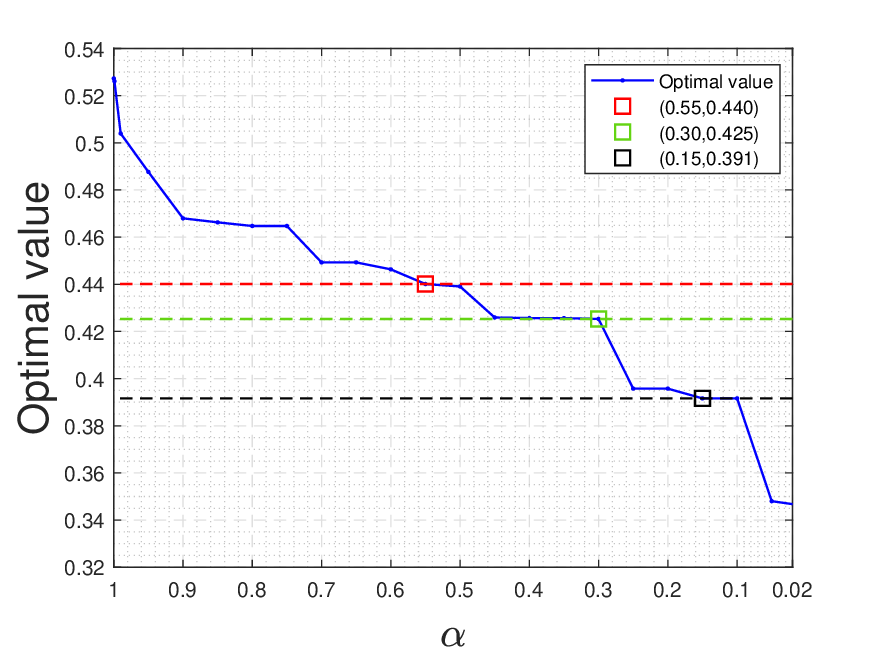}}
  \text{\tiny (a) optimal values}
\endminipage
\hfill
\minipage{0.5\textwidth}
  \centering
  \includegraphics[width=5.0cm]{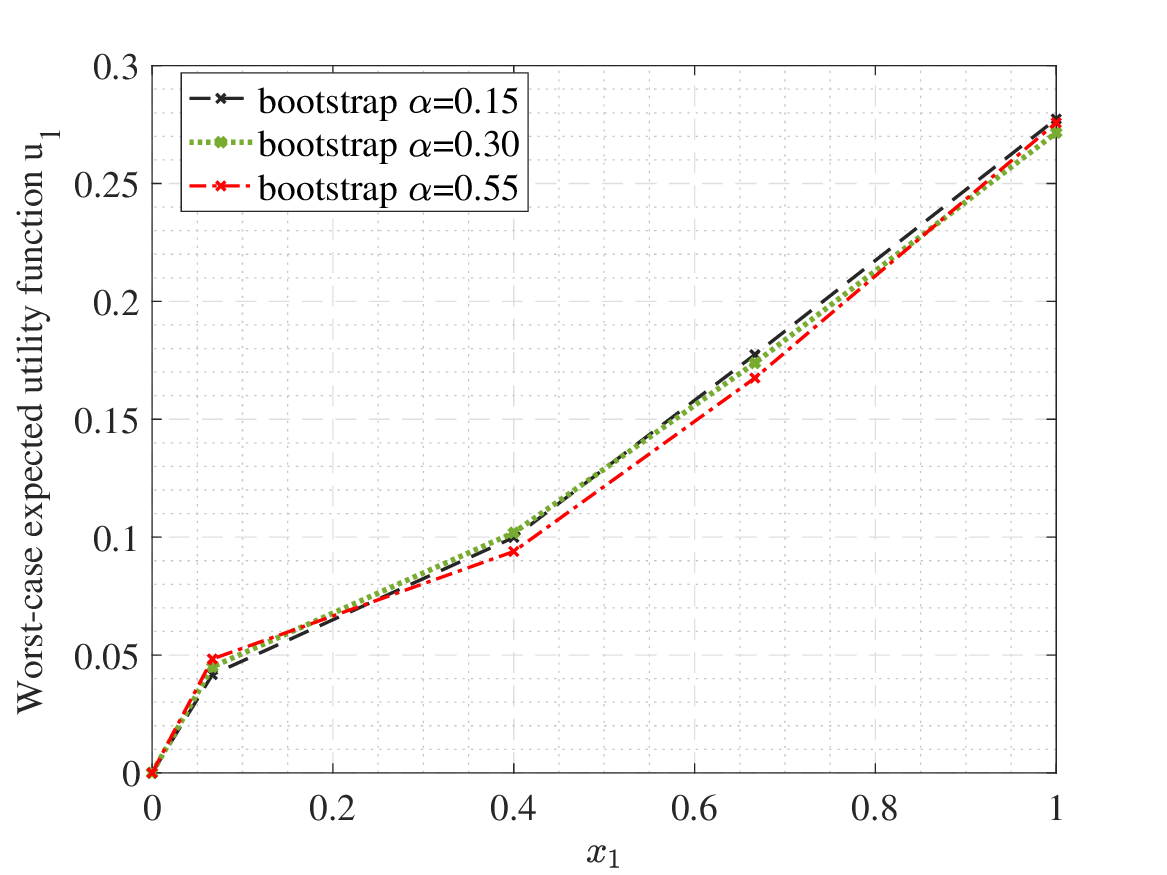}
  \text{\tiny (b) utility function of attribute $1$.}
\endminipage
\hfill
\minipage{0.5\textwidth}
  \centering
  \includegraphics[width=5.0cm]{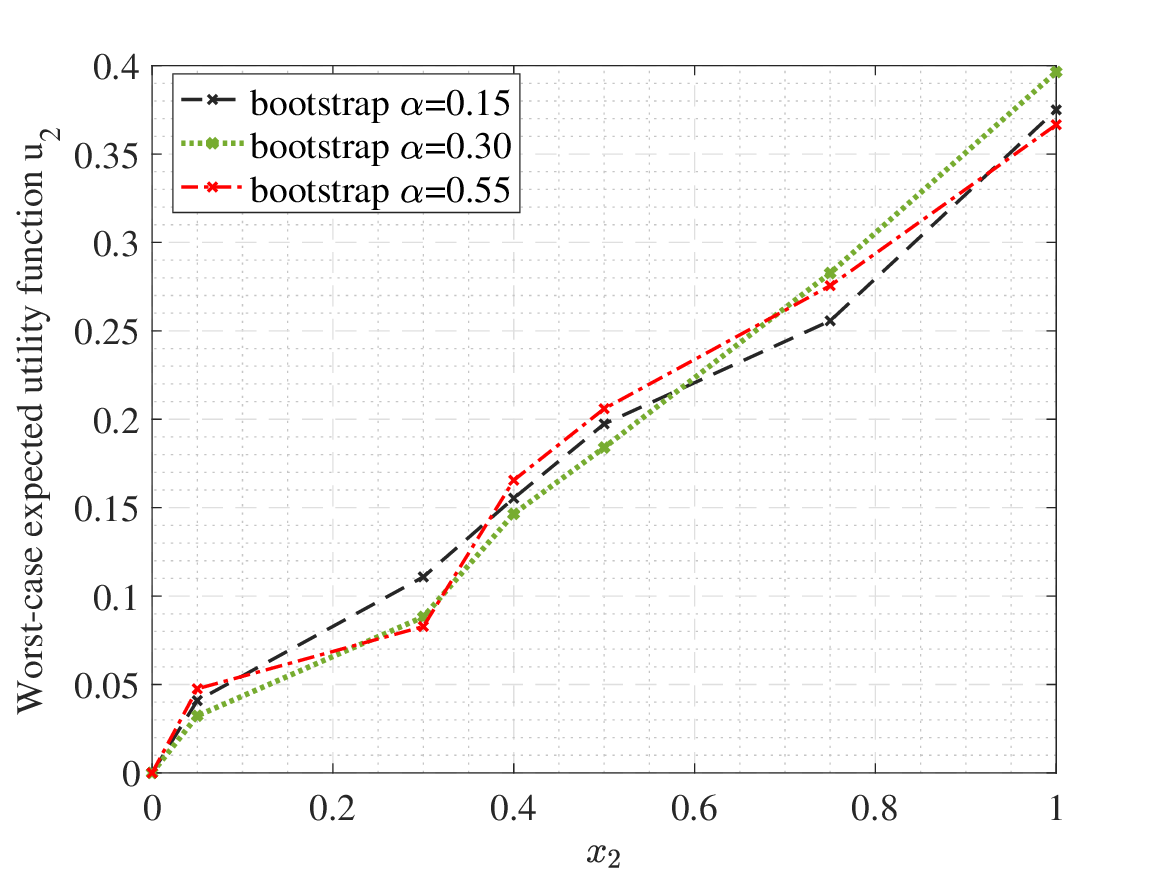}
  \text{\tiny (c) utility function of attribute $2$.}
\endminipage
\hfill
\minipage{0.5\textwidth}
  \centering
  \includegraphics[width=5.0cm]{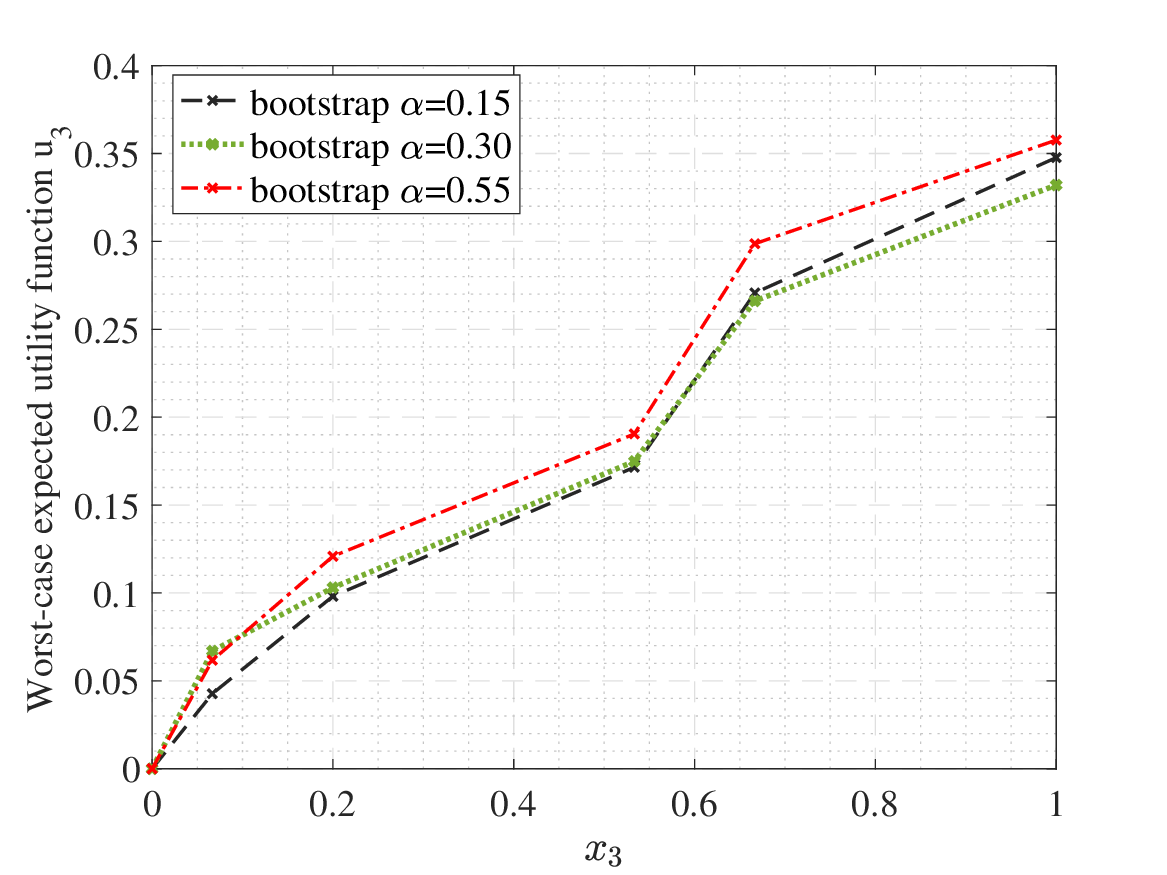}
  \text{\tiny (d) utility function of attribute $3$.}
\endminipage
\caption{\small (a) 
The optimal value of DUPRO with ambiguity set by bootstrap;
(b)-(d) The worst expected utility functions when 
$\alpha = 0.15$, $\alpha=0.30$ and $\alpha=0.55$.}
\label{utility_worst_1}
\vspace{-0.5cm}
\end{figure}

 With the sample data, 
 we are able to calculate Tukey's depth $d_{(1)},\cdots,d_{(K)}$ and then construct a $100(1-\alpha) \%$ bootstrap confidence region
 $\fC_{1-\alpha}^{*}$,
 based on which,
 we can obtain the ambiguity set $\mathcal{P}_B(\alpha)$ in \eqref{set:fP_B} (or
$\fF(
\mathcal{P}_B(\alpha)
)$
 in (\ref{eq:fF_PB_alpha})). 
 Next, we solve problem  (\ref{eq:DUPRO_FB_bootstrap}), equivalent to \ref{DUPRO},
 via solving program (\ref{pp:PRO_B_equ}) 
 and obtain an optimal solution $x^*$.
 By substituting $x^*$ into the inner problem $\min_{v\in \fF(
 \mathcal{P}_B(\alpha)
 )} v^Tf(x^*)$,
 we obtain the worst $v^*$ and subsequently 
 the worst expected utility function.
Figure~\ref{utility_worst_1}~(a)
depicts
 convergence of the  optimal value of
 problem (\ref{pp:PRO_B_equ}) 
 as $\alpha$ reduces
 from 
 $99.99\%$ to $2\%$.
 For example, we present   the worst-case expected utility 
 functions $u_1$, $u_2$ and $u_3$ when $\alpha$ = 0.15, 0.30, and 0.55, respectively, in Figure~\ref{utility_worst_1} (b)-(d).

\begin{figure}[!htbp]
 \centering
 \minipage{0.33\textwidth}
  \centering
  \includegraphics[width=4.0cm]{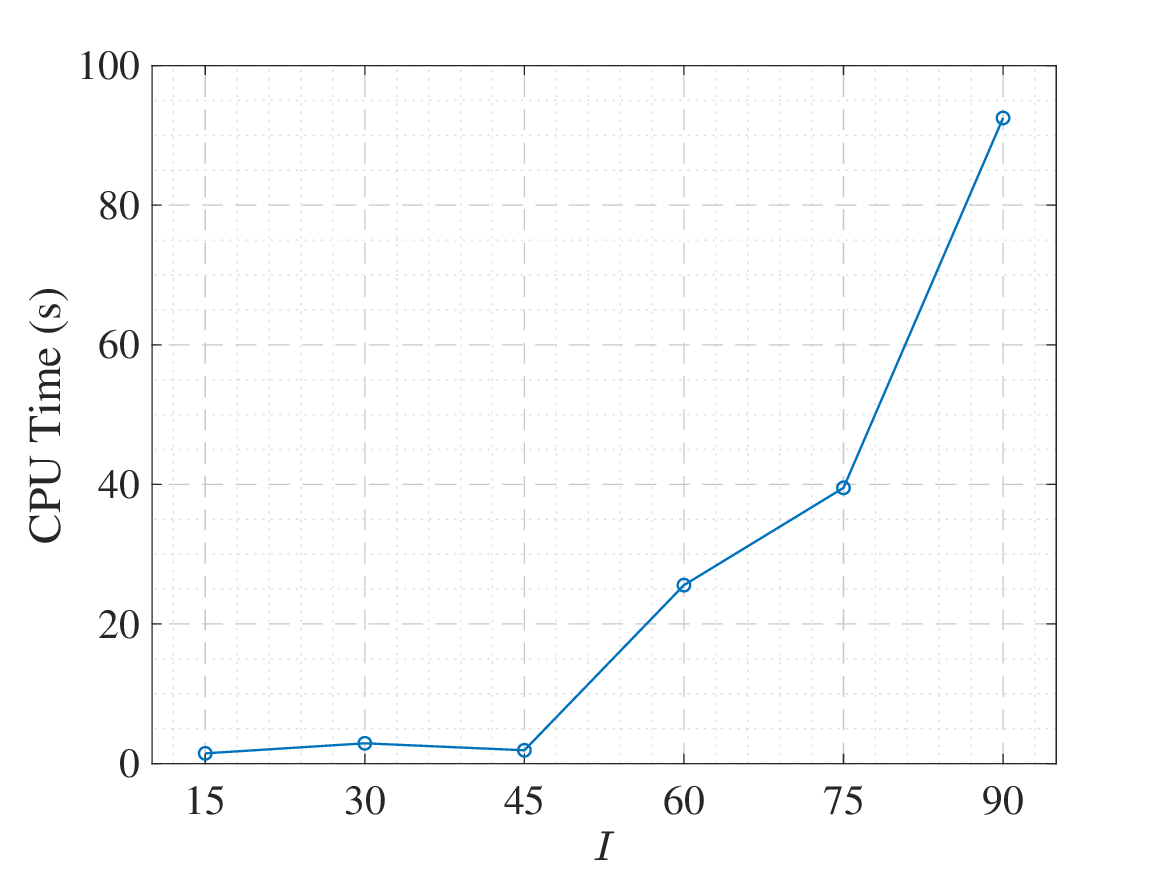}
  \text{\tiny (a) $\alpha = 0.15$}
\endminipage
 \minipage{0.33\textwidth}
  \centering
  \includegraphics[width=4.0cm]{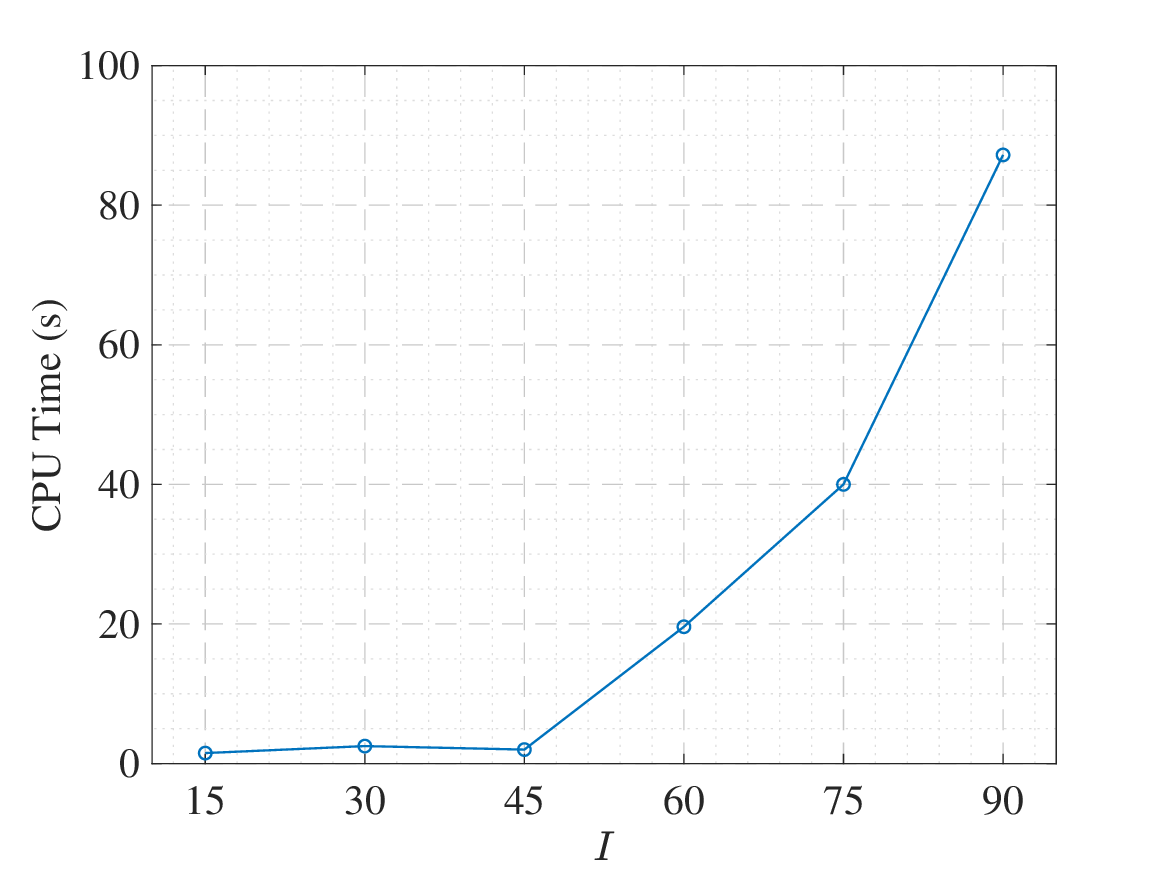}
   \text{\tiny (b) $\alpha = 0.30$}
\endminipage
 \minipage{0.33\textwidth}
  \centering
  \includegraphics[width=4.0cm]{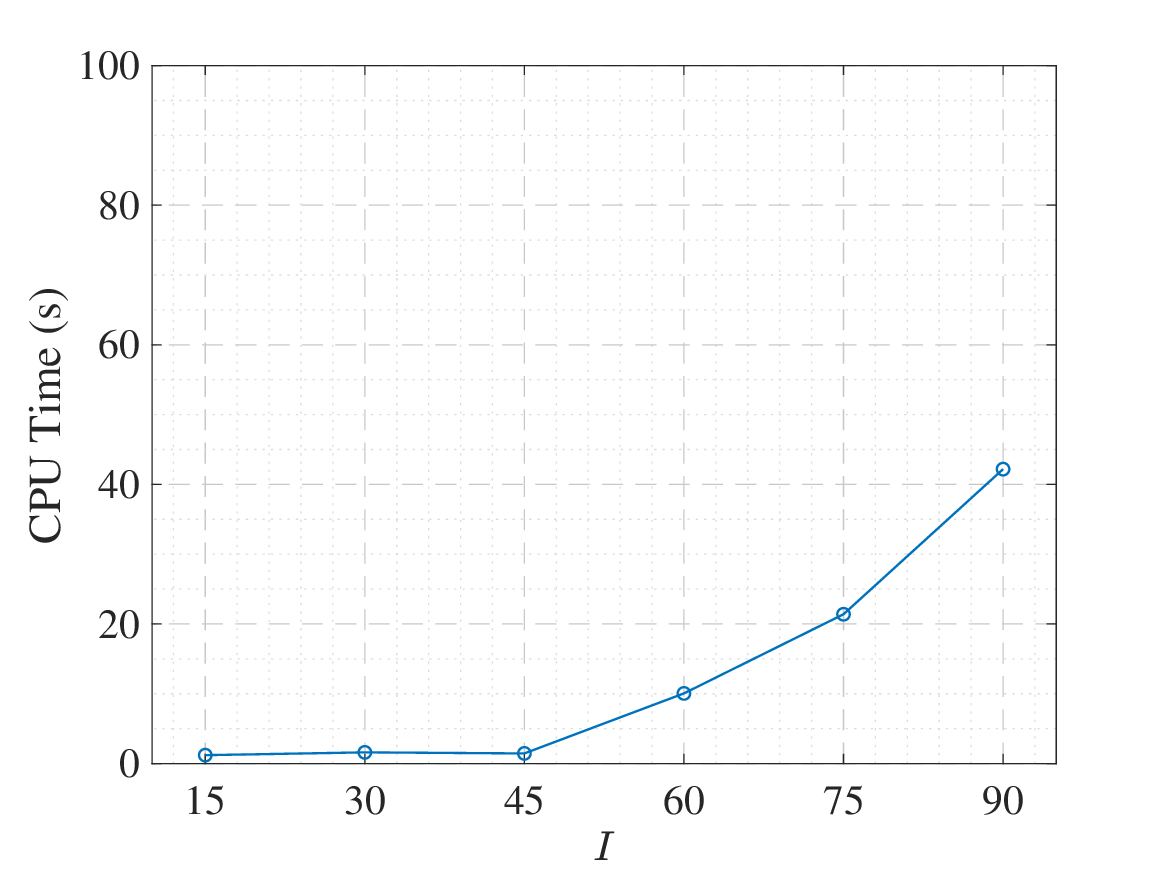}
   \text{\tiny (c) $\alpha = 0.55$}
\endminipage
\caption{Running time of \ref{DUPRO} (problem \eqref{pp:PRO_B_equ}).}
\label{fig:CPU}
\end{figure}

Next, we test the scalability of \ref{DUPRO} when $I$,
the total  number 
of the linear pieces of the utility functions, increases. Concomitantly, the complexity of problem \eqref{pp:PRO_B_equ} 
increases with more binary variables. We equally partition each of the intervals 
between the breakpoints 
given in Table~\ref{tab:breakpoints_3attributes} into $\tau$ subintervals 
with $I = 15 \tau$. 
In the study, we vary $\tau$ from 1 to 6.
Figure~\ref{fig:CPU}
depicts 
the running CPU time 
for solving problem \eqref{pp:PRO_B_equ} with $N = 50$, $K = 10,000$, and $\alpha$ varying from 
$0.15$ to $0.55$.
It shows that 
the 
increase in $I$ 
yields an exponential growth 
increase of the running time. 

\subsubsection{Test case \ref{test:b}: Convergence analysis}
\label{subsubsec:Convergence analysis}
We investigate the convergence result in Theorem~\ref{thm:DUPRO-Bootstrap-convg}. 
In the test problem, the true random utility function is step-like which means
\begin{align*} 
    \bbe[u(x;V)] = \bbe[V]^T f(x) = \frac{e_{15}^T f(x)}{15} ,
\end{align*}
where $f(x)$ is defined as in (\ref{eq:u-vec-rep}) and all of its breakpoints are given in Table \ref{tab:breakpoints_3attributes}. Let
$\vt = \max_{x \in \fX} e_I^T f(x) / 15$.
We use the bootstrap approach to construct 
the ambiguity set $\mathcal{P}_{B} (\alpha)$ in \eqref{set:fP_B}. Let $\tilde \vt_{{\rm bts}}$ be the optimal value of \ref{DUPRO} (problem (\ref{pp:PRO_B_equ})). The set 
${\cal P}_{B} (\alpha)$
is an approximation to 
${\cal P}_{\hat{B}} (\alpha)$
defined as in (\ref{eq:fF_PB_alpha-infinite-K}), and thus $\tilde \vt_{{\rm bts}}$ approximates $\vt_{{\rm bts}}$ defined as in 
 (\ref{eq:minmax-BooTS}).
We choose $K = 10,000$ replications for bootstrapping in this test. Such a large $K$ guarantees the quality of the approximation. Hence, we can obtain  
\bgeqn \label{eq:sim_prob}
\prob\left(|\vt_{{\rm bts}}-\vt|\geq \epsilon \right)
\simeq
\prob\left(|\tilde \vt_{\rm  bts}-\vt |\geq \epsilon \right),
\edeqn 
 where the right-hand side is the complementary cumulative distribution function (ccdf) of $| \tilde \vt_{\rm  bts} -\vt|$ with respect to $\epsilon$.

\begin{figure}[ht!]
\minipage{\textwidth}
 \centering
  \includegraphics[width=8.5cm]{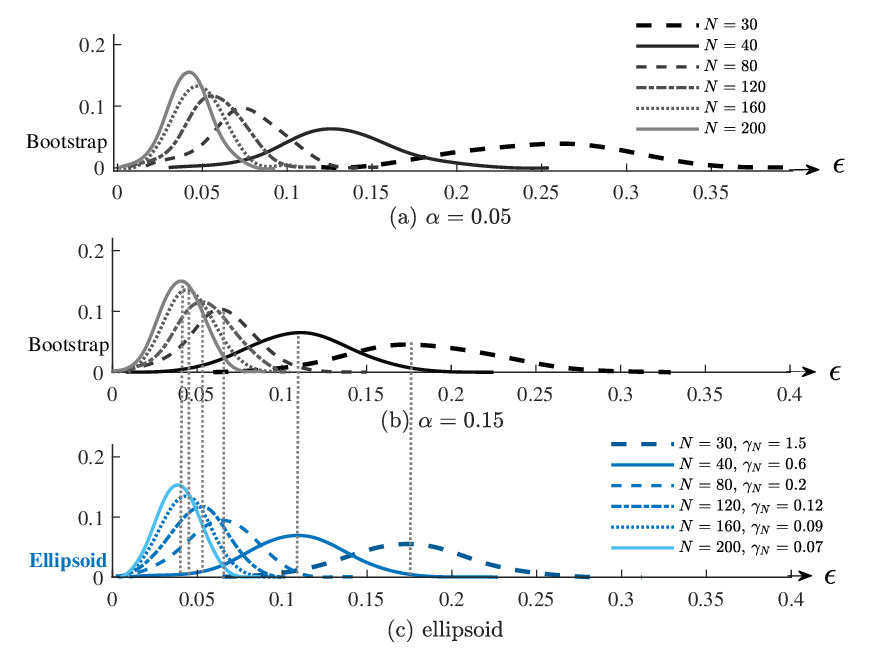}
  \label{fig:comparison-a}
\endminipage
\caption{\small  p.d.f. of 
$|\vt_{\rm bts}-\vt_I|$
and  p.d.f. of $|\vt_{\rm elp}-\vt_I|$}
\label{fig:comparison}
\vspace{-0.5cm}
\end{figure}
 
We run a Monte Carlo simulation to generate 1000 observations of $|\tilde \vt_{\rm  bts} - \vt|$. On this basis, we 
depict the empirical distribution to approximate the probability distribution function (pdf) of $| \tilde \vt_{\rm  bts} -\vt|$. Note that the ccdf is the area under the pdf on the right of 
a given critical value $\epsilon$. 
To ease the exposition, 
we omit ``empirical'' in the statement below. Figure \ref{fig:comparison} 
depicts the pdfs when $\alpha$ is 0.05 or 0.15, and $N$ (the size of the original random sample) varies from 30 to 200. In comparison, we also 
plot the pdf's
of $|\vt_{{\rm elp}} -\vt|$ in the case of using ellipsoidal ambiguity set in \eqref{eq:convg-opti-value-thm2}. For a fixed $\alpha$, an increase 
of $N$ results in a narrower and more concentrated pdf curve, which shifts to the left and its right tail gradually diminishing or becoming lighter. This shows the convergence of $\tilde \vt_{\rm bts}$ to $\vt$ in probability as Theorem~\ref{thm:DUPRO-Bootstrap-convg} claims. That is, the right-hand-side probability in \eqref{eq:sim_prob} dwindles for any given $\epsilon$. For example, when $\alpha=0.15$ and $\epsilon=0.05$, the probability is 
$0.999$, $0.994$, $0.882$, $0.758$, $0.558$
and $0.390$,
respectively, as $N$ increases from 30 to 200. Moreover, we observe that the pdf differs 
with respect to varying $\alpha$ value
when $N$ is small ($N$ = 30, 40 in this case), 
becomes indifferent 
$N$ is large. Theorem~\ref{thm:DUPRO-Ellip-optim} indicates that the probability of convergence of $\vt_{{\rm elp}}$ relies on a sequence of $\gamma_N$ in 
infinite descent. We choose $\gamma_N$ regarding $N$ to obtain the results similar to the case with $\alpha=0.15$. 
It is worth noting that Proposition \ref{Prop:bootstrap-confidence} recommends too large $\gamma_N$, 
under which we would have $\vt_{{\rm elp}} = 0$ in our tests.

\subsubsection{Test case \ref{test:c}: Out-of-Sample Test}
\begin{figure}[!htbp]
\minipage{0.33\textwidth}
 \centering
 \includegraphics[width=4.1cm]{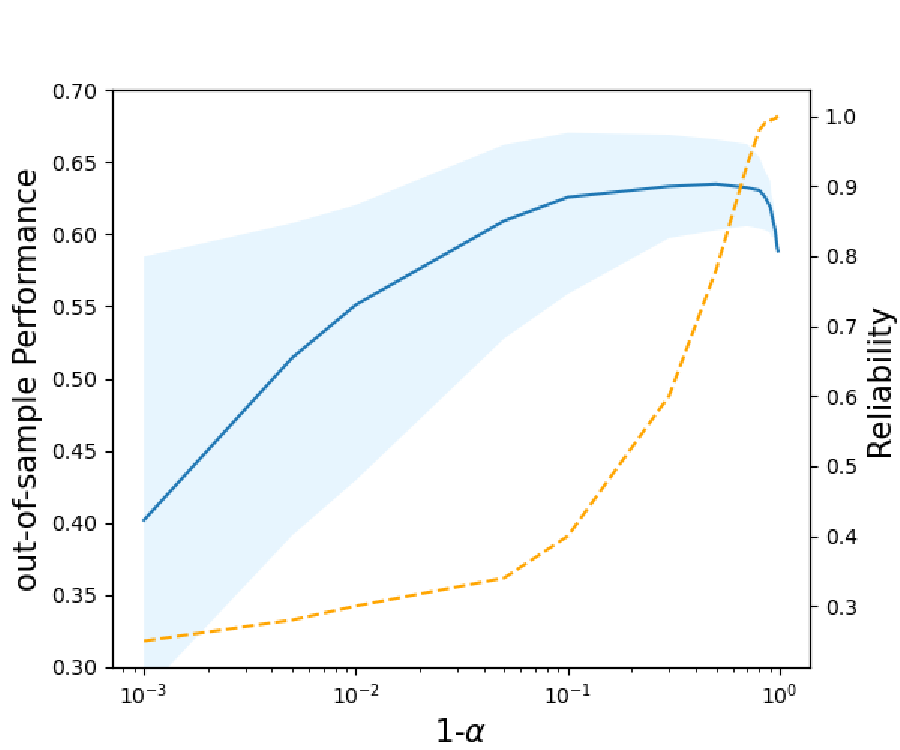}
  \textrm{\tiny (a) $N=20$ \\ (SAA: $0.3994$,  $[0.2766, 0.5803]$).}
\endminipage\hfill
\minipage{0.33\textwidth}
  \centering
 \includegraphics[width=4.1cm]{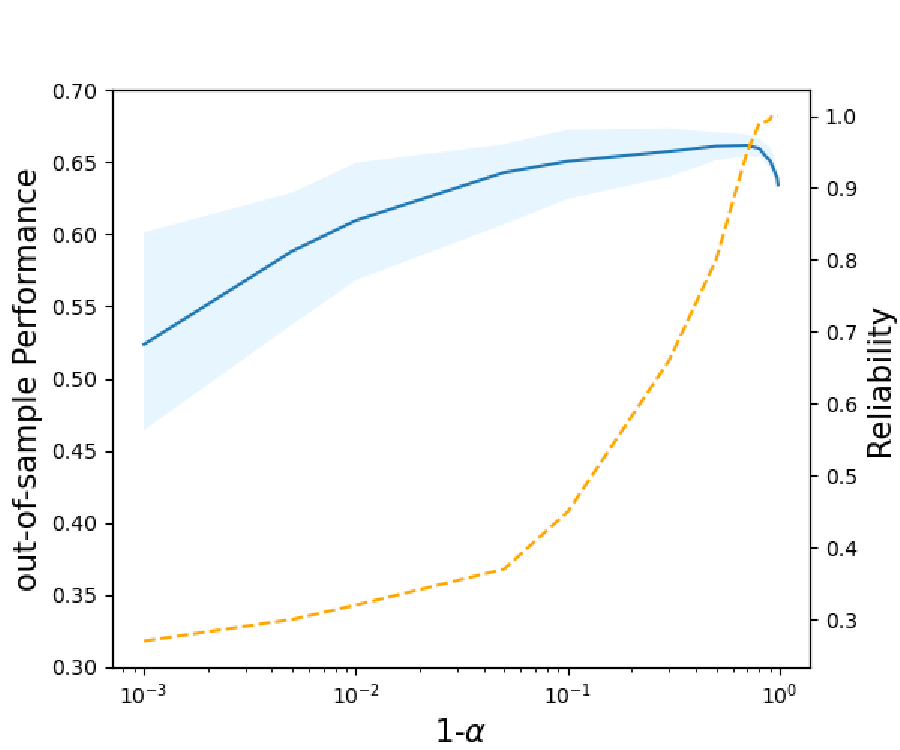}
  \textbf{\tiny (b) $N=100$ \\ (SAA: $0.5192$,  $[0.4539, 0.6002]$)}
\endminipage
\hfill
\minipage{0.33\textwidth}
  \centering
 \includegraphics[width=4.1cm]{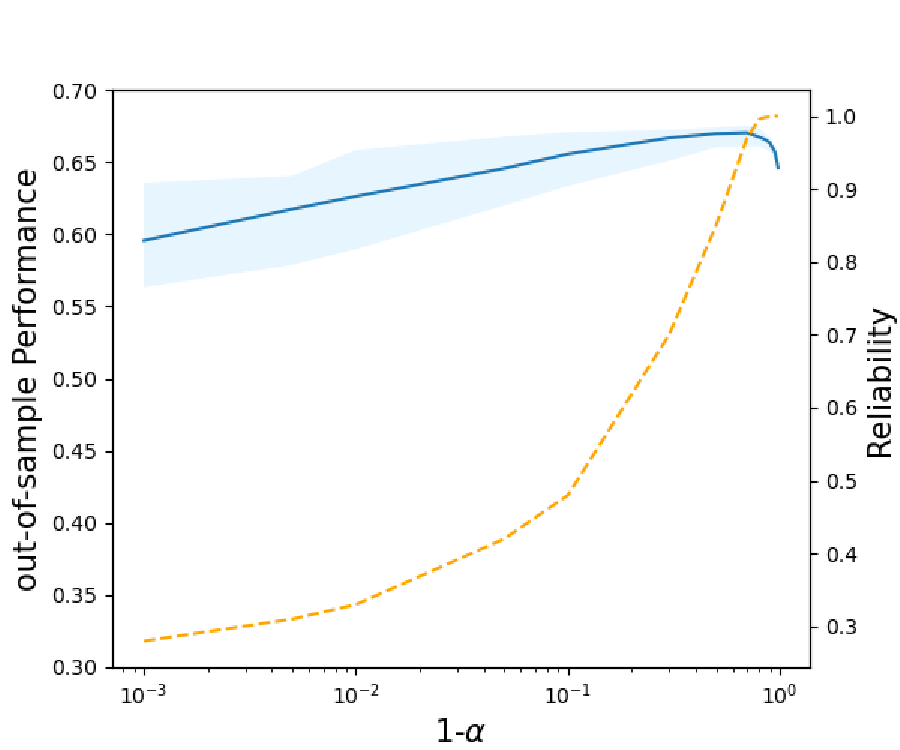}
  \textbf{\tiny (c) $N=200$ \\ (SAA: $0.5841$,  $[0.5615, 0.6341]$)}
\endminipage
\caption{\small The out-of-sample performance  of \ref{DUPRO}
((i) the solid line for the sample mean of $J_N(\alpha)$,  the shaded area for its 20-80\% sample quantiles, the dotted line for its reliability; (ii) the out-of-sample performance of the SAA solution given in the caption of each subfigure in the form of (sample mean, [20\% sample percentile, 80\% sample percentile]).)}  
\label{fig:oos-compar}
\vspace{-0.5cm}
\end{figure}

Analogous to \cite{Esfahani2018}, we analyze the out-of-sample performance of the bootstrap-based \ref{DUPRO}. In the out-of-sample performance tests, 
we change the true probability distribution of vector
 $V$ (of dimension $I=15$) to an asymmetric Dirichlet distribution with
 order $15$ with parameters $(0.25e^T_5, 0.75e^T_5, 0.50e^T_5)\in\mathbb{R}^{15}$.
For a given random sample, $V^1, \dots, V^N$, we solve
 \ref{DUPRO} (problem (\ref{pp:PRO_B_equ}))
with  a fixed $\alpha$ to obtain 
 an optimal solution $\hat x_N(\alpha)$ and 
 then examine its performance via 
the conditional expectation of the utility value
\begin{align}
 J_N(\alpha) =  \bE [u(\hat x_N(\alpha); V) \;|\; V^1, \dots, V^N].   
\end{align}
This study runs a Monte Carlo simulation
to generate $300$ observations of 
$J_N(\alpha)$. 
Figure \ref{fig:oos-compar}~(a)-(c), 
depicts the 20-80\% sample quantiles  of $J_N (\alpha)$ (shaded area)
and the sample means (solid curve in the area) when $1-\alpha$ varies from $10^{-3}$ to 1, for $N$ = 20, 100, and 200, respectively. Also, the dotted curves visualize the reliability that represents the empirical probability of the event $J_N(\alpha) \ge \tilde \vt_{{\rm bts}}$ via the simulation. Recall that $\tilde \vt_{{\rm bts}}$ is the optimal value of \ref{DUPRO} (problem \eqref{pp:PRO_B_equ}).

The out-of-sample performance improves until  $1-\alpha$ reaches a large critical value (close to 0.8 in this test) and then deteriorates. Note that Figure \ref{fig:oos-compar} uses a log-scale axis to draw the effectiveness of \ref{DUPRO} with a small-sized ambiguity set. The concave curves of sample mean indicate that the out-of-sample performance has a higher rate of improvement when $1-\alpha$ is small, while \ref{DUPRO} becomes too conservative when $1-\alpha$ approaches $1$. The robustness of \ref{DUPRO} is also reflected by the range between the 20-80\% sample quantiles shrinking as $1-\alpha$ increases. We also investigate reliability of $J_N(\alpha)$. In Figure \ref{fig:oos-compar}, the reliability attains the maximum when $1-\alpha=0.8$.

On the other hand, we can see that the increase of original sample size $N$ 
improves
the out-of-sample performance entirely, 
not only 
lifting 
the curve of the sample mean but also 
narrowing down 
the range 
of the  20-80\% sample quantiles. 
The results are consistent 
with 
the convergence analysis in Section \ref{subsubsec:Convergence analysis}. 
Moreover, Figure \ref{fig:oos-compar} shows that \ref{DUPRO}, 
as a data-driven approach, 
is 
affected 
by the level of distributional uncertainty. 
As $N$ increases, 
we deduce from (\ref{eq:set:fC})-(\ref{set:fP_B}) that
the gap between the ambiguity sets $\mathcal{P}_{B} (0)$ and $\mathcal{P}_{B} (\alpha)$ diminishes 
for any $\alpha \in (0, 1]$,
thereby compromising the superior performance of \ref{DUPRO} in comparison to the sample average approximation (SAA) 
approach (a trivial \ref{DUPRO} case with the singleton ambiguity set $\fF(
		\cP_B (0))$ including the sample mean only).  
Given that $1-\alpha$ changes in $[10^{-3}, 1]$, the largest relative improvement of \ref{DUPRO} over the SAA approach is $(0.6322-0.3994)/0.6322 = 36.82\%$
when $N= 20$ in Figure \ref{fig:oos-compar}(a) and is (0.6656-0.5841)/0.6656 = 7.23\% when $N= 200$ in Figure \ref{fig:oos-compar}(c).

\subsection{Project Investment}
\label{subsec:Project Investment}
	
We now illustrate how to
implement \ref{DUPRO} into a project investment problem. In this problem, an automotive manufacturer needs to learn 
consumer preference towards passenger vehicles, 
and on this basis,  to pursue the best portfolio investment 
with a fixed budget among the 10 projects 
listed in Table~\ref{tab:project-list}, see Appendix~\ref{App:table-project-case-2}. 
These candidate projects include
\begin{itemize}
\item {\em safety promotion:} 
design a new structure to alleviate the impacts of the collision;
\item {\em new car model development:} create a concept car with the aim of a more fashion style and higher market acceptance before actually producing it;
\item {\em engine upgrade:} upgrade the current engine and its vibration sensor system to achieve more horsepower and less engine noise;
\item {\em e-platform development:} improve the human-vehicle interaction experience
and visualize the performances of the car;
\item {\em computational fluid dynamics (CFD) testing system development:} implement
related fluid mechanics and numerical analysis into a testing platform to analyze the performance of concept cars;
\item {\em common modular platform (CMP) development:} develop a platform used for subcompact and compact car models with internal combustion engine and battery-electric cars;
\item {\em checking fixture promotion:} enhance the performance of checking fixture
to control the dimensions of auto parts (such as trim edge, surface profile, flatness, etc.) in a more convenient way for mass production of parts detection;
\item {\em noise, vibration, harshness (NVH) digitalization:}  incorporate  the characteristics of  the noise, vibration and harshness of vehicles to a digital platform 
with the target of more efficient computation in evaluating the driver satisfactions;
\item {\em driving assistance system development:} incorporate the latest interface standards and running multiple vision-based algorithms to support real-time multimedia, vision co-processing, and sensor fusion subsystems;
\item {\em digitalization of marketing network:}  build a digital platform used for all the dealers across different regions to share the customer resources, the service standards and the marketing strategies. 
\end{itemize}

These projects may lead to 
enhancement of
eight attributes 
related to vehicle performance, economy, after-sales service, etc, among which the attributes including {\em sale price, fuel consumption, depreciation rate (yearly), and estimated maintain and repair (M$\&$R) fee (in five years)} are related to consumers' major economic considerations when purchasing new cars, {\em wheel base, acceleration, comfort rating} are related 
to consumers' major consideration 
of vehicle performance, 
 and {\em dealership}, represented by the number of dealers able to deal with 
 certain services, directly influences consumers' satisfactions on after-sale services. 
 Table~\ref{tab:project-list} displays the estimated consequences and costs of the projects. The column of ``Base model'' gives the baseline attribute values before the investment, of which the vector is denoted by $x^0$, while the other columns show the vectors of attribute increment values, denoted by $y^d$ for $d = 1, \dots, 10$, after the investment on the projects.

\textbf{Application of \ref{DUPRO}:}
Denote by $z \in \{0, 1\}^{10}$ the decision vector of which a component $z_d = 1$ means the project $d$ is selected in the investment and otherwise $z_d = 0$. Assume that the improvement on the attributes yielded by each project is independent of the one by any other project. As a result, the attribute values achieved in the investment are obtained as
	
	\begin{align} \label{eq:cs_1}
		\hat x(z) = x^0 +  \sum_{d=1}^{10} z_d y^d.
	\end{align}
	Let $\Phi$ be the budget limit of the investment and $b_d$ the cost of project $d$ given in Table~\ref{tab:project-list}. A budget constraint is thus described as
	\begin{align} \label{eq:cs_2}
		\sum_{d=1}^{10} b_d z_d \le \Phi.
	\end{align} 
	Combining \eqref{eq:cs_1} and \eqref{eq:cs_2}, we formulate the feasible region of the attribute value achieved by the investment as
	\begin{align}
		\fX = \{ \hat x(z)  \;|\; \text{$z \in \{0, 1\}^{10}$ satisfies \eqref{eq:cs_2}}  \}.
	\end{align}
	Consequently, \ref{DUPRO} in this case is specified as
		\begin{align}
			\label{prob:project-car-max}
			\max_{x \in \fX} \min_{P \in \fP_B(\alpha)}  \bE_P [u (x; V)],
		\end{align}
	where $u (\cdot; V)$ represents the uncertain consumer preference and $\fP_B(\alpha)$ is the ambiguity set of the distribution specified by 
	means 
	of bootstrap given in \eqref{set:fP_B}.

\begin{figure}[htbp]
\vspace{-0.2cm}
\minipage{0.5\textwidth}
  \centering
  \includegraphics[width=6.0cm]{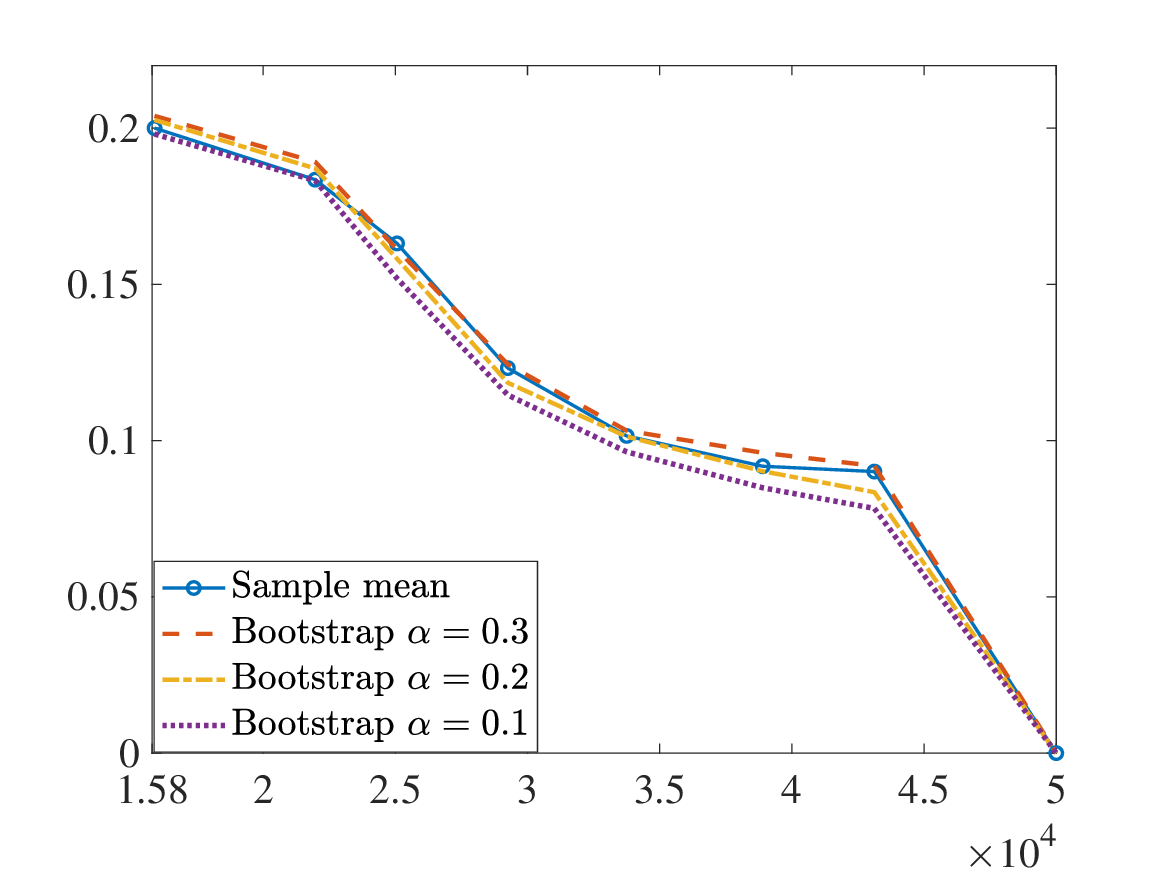}
  \textbf{\tiny (a) Price}
\endminipage
\hfill
\minipage{0.5\textwidth}
  \centering
  \includegraphics[width=6.0cm]{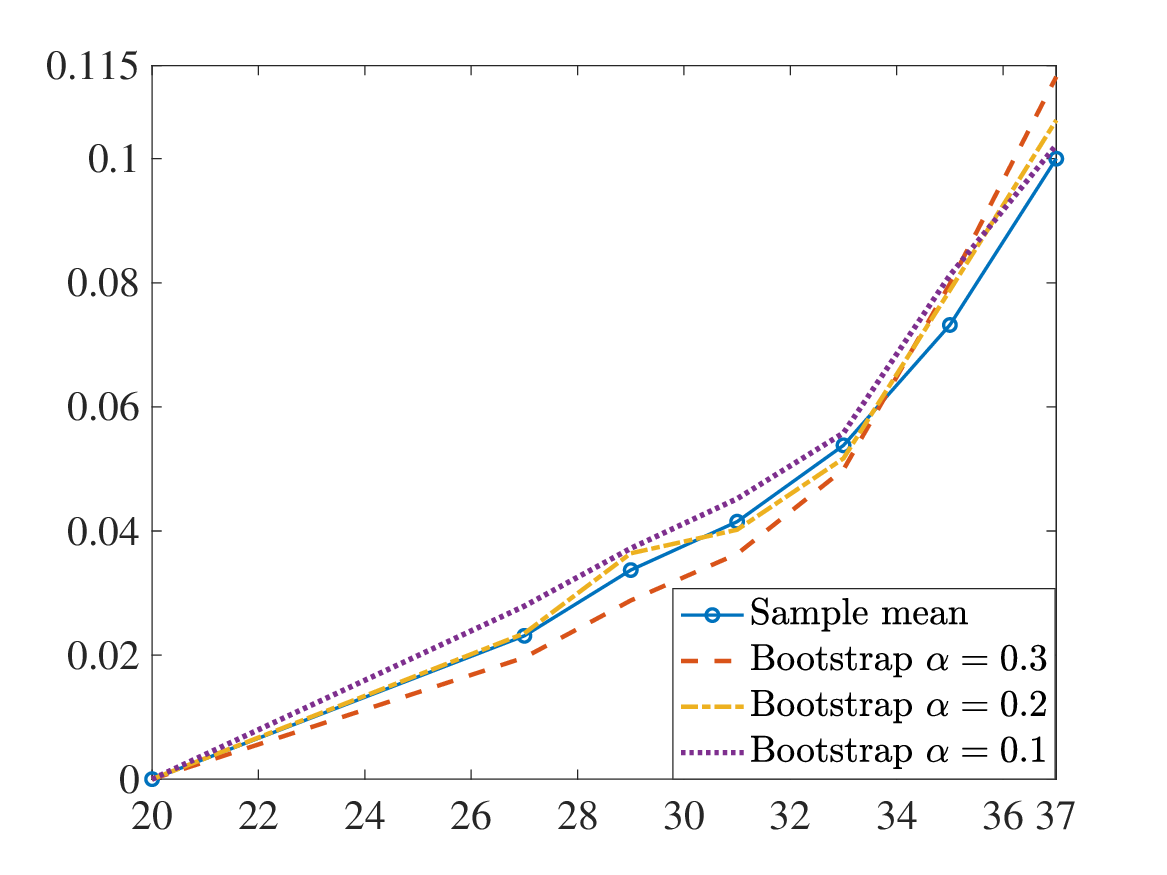}
  \textbf{\tiny (b) Fuel Consumption}
\endminipage
\hfill
\minipage{0.5\textwidth}
  \centering
  \includegraphics[width=6.0cm]{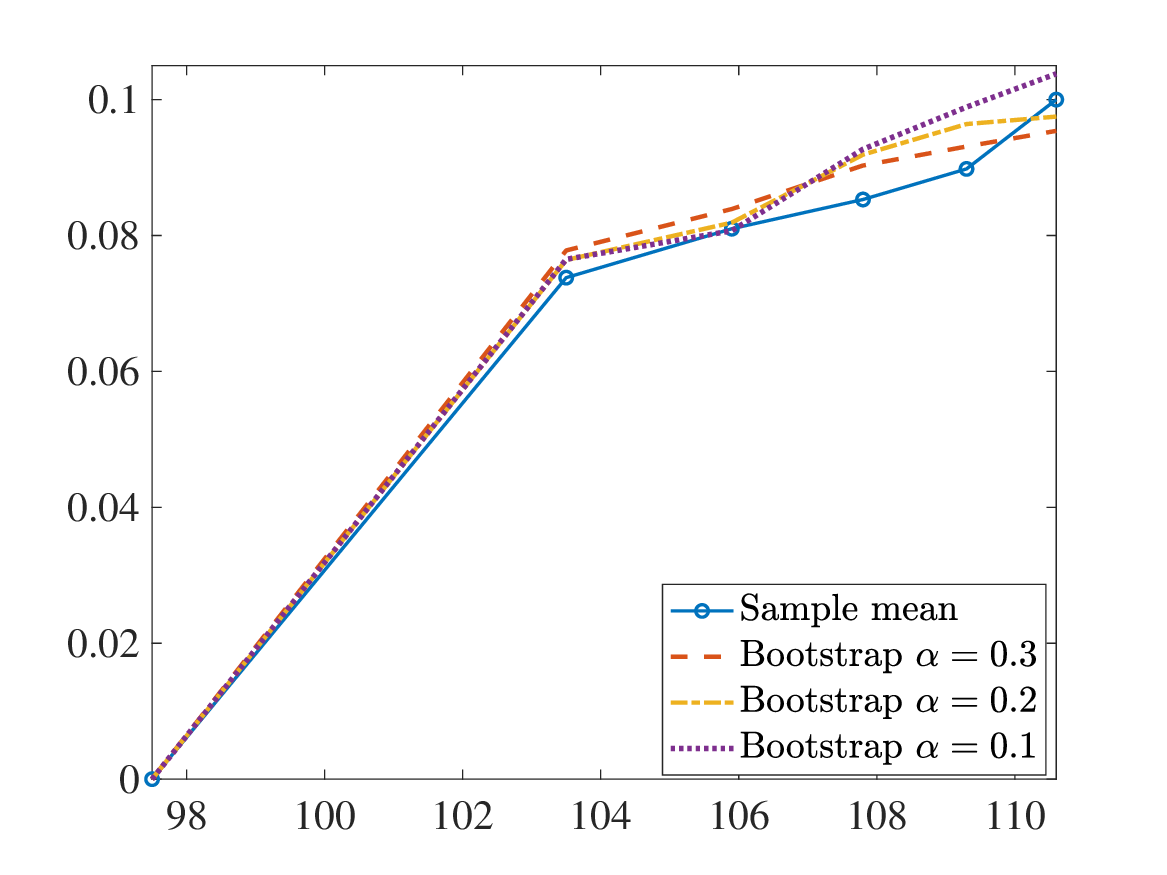}
  \textbf{\tiny (c) Wheelbase}
\endminipage
\hfill
\minipage{0.5\textwidth}
  \centering
  \includegraphics[width=6.0cm]{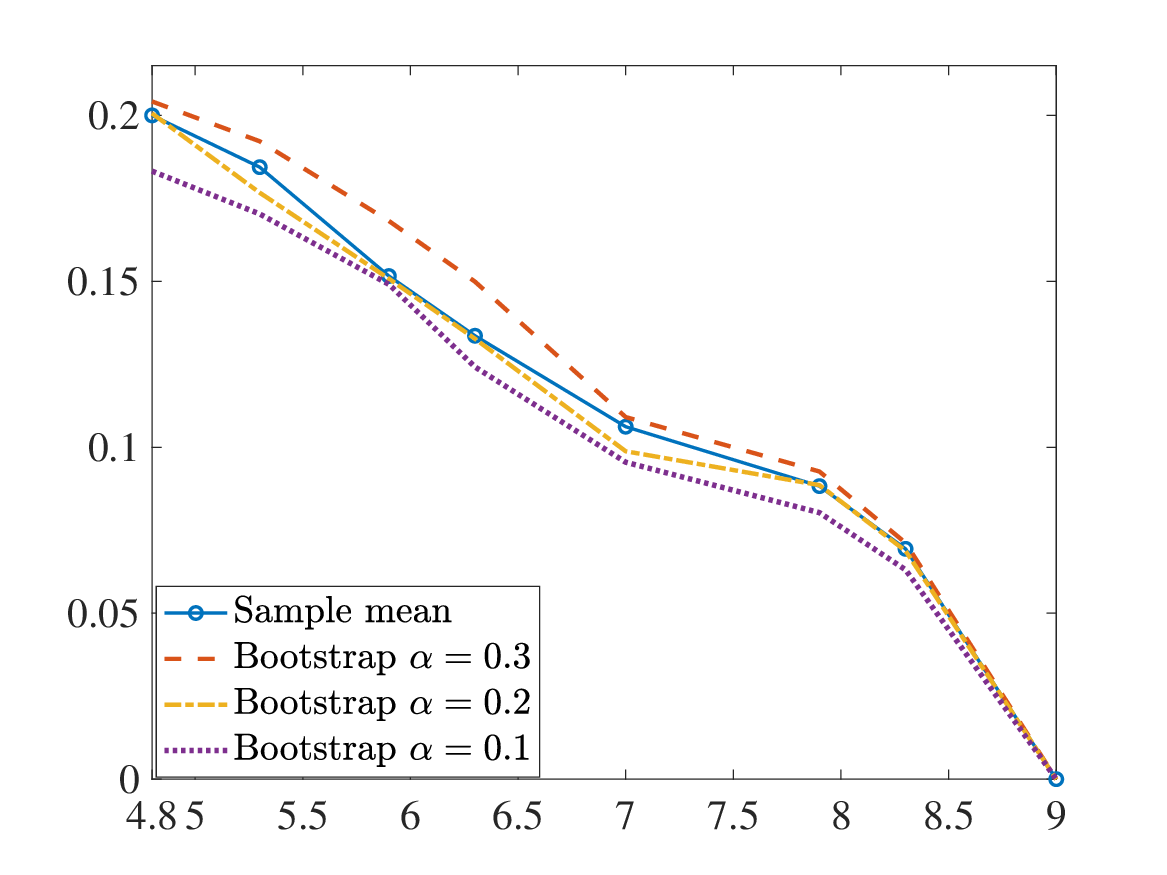}
  \textbf{\tiny (d) Acceleration}
\endminipage
\hfill
\minipage{0.5\textwidth}
  \centering
  \includegraphics[width=6.0cm]{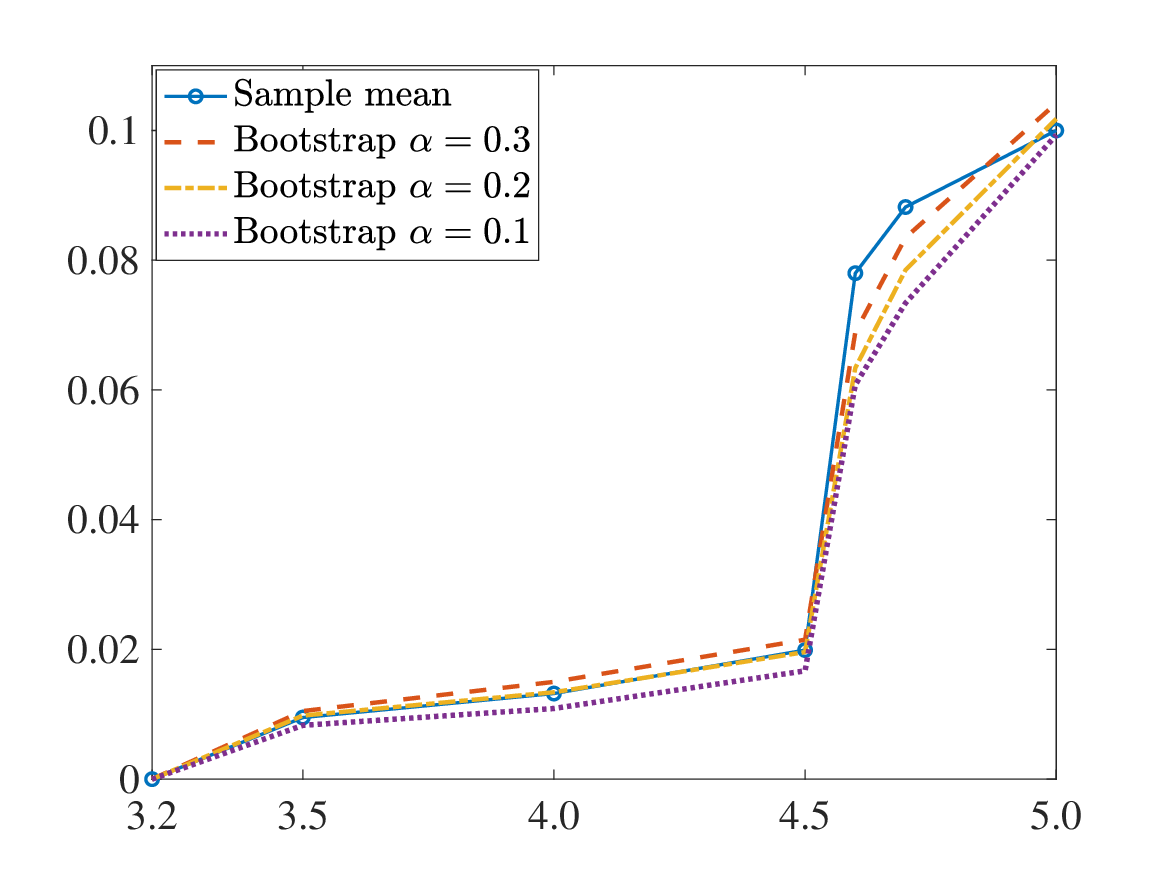}
  \textbf{\tiny (e) Comfort}
\endminipage
\hfill
\minipage{0.5\textwidth}
  \centering
  \includegraphics[width=6.0cm]{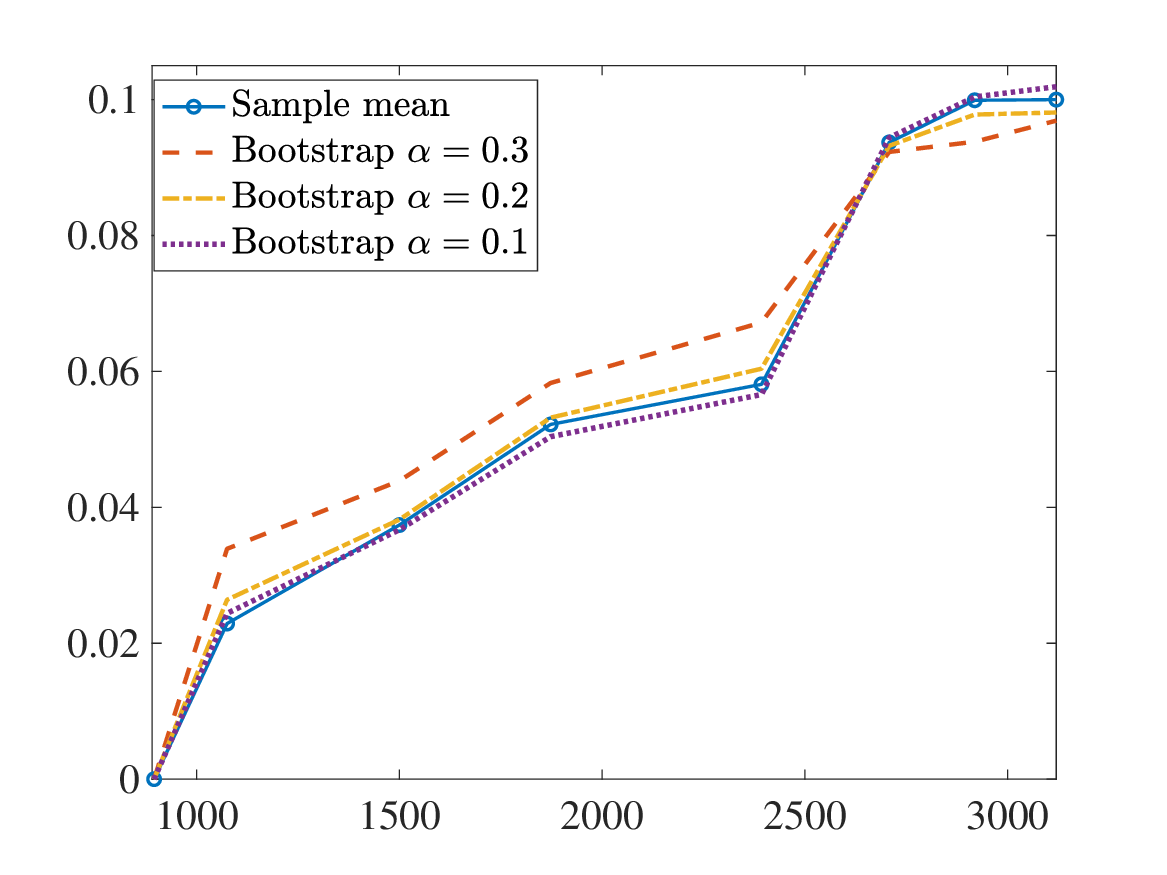}
  \textbf{\tiny (f) Dealership}
\endminipage
\hfill
\minipage{0.5\textwidth}
  \centering
  \includegraphics[width=6.0cm]{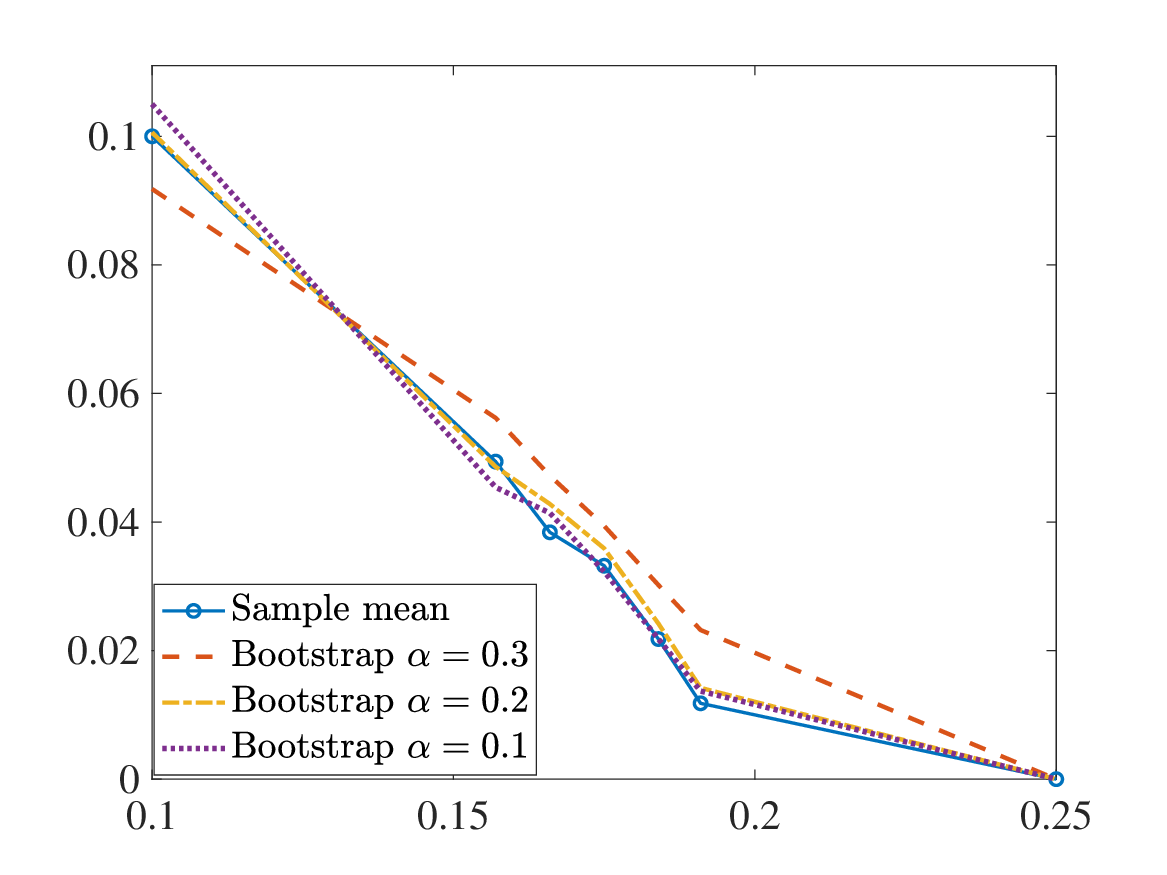}
  \textbf{\tiny (g) Depreciation Rate}
\endminipage
\hfill
\minipage{0.5\textwidth}
  \centering
  \includegraphics[width=6.0cm]{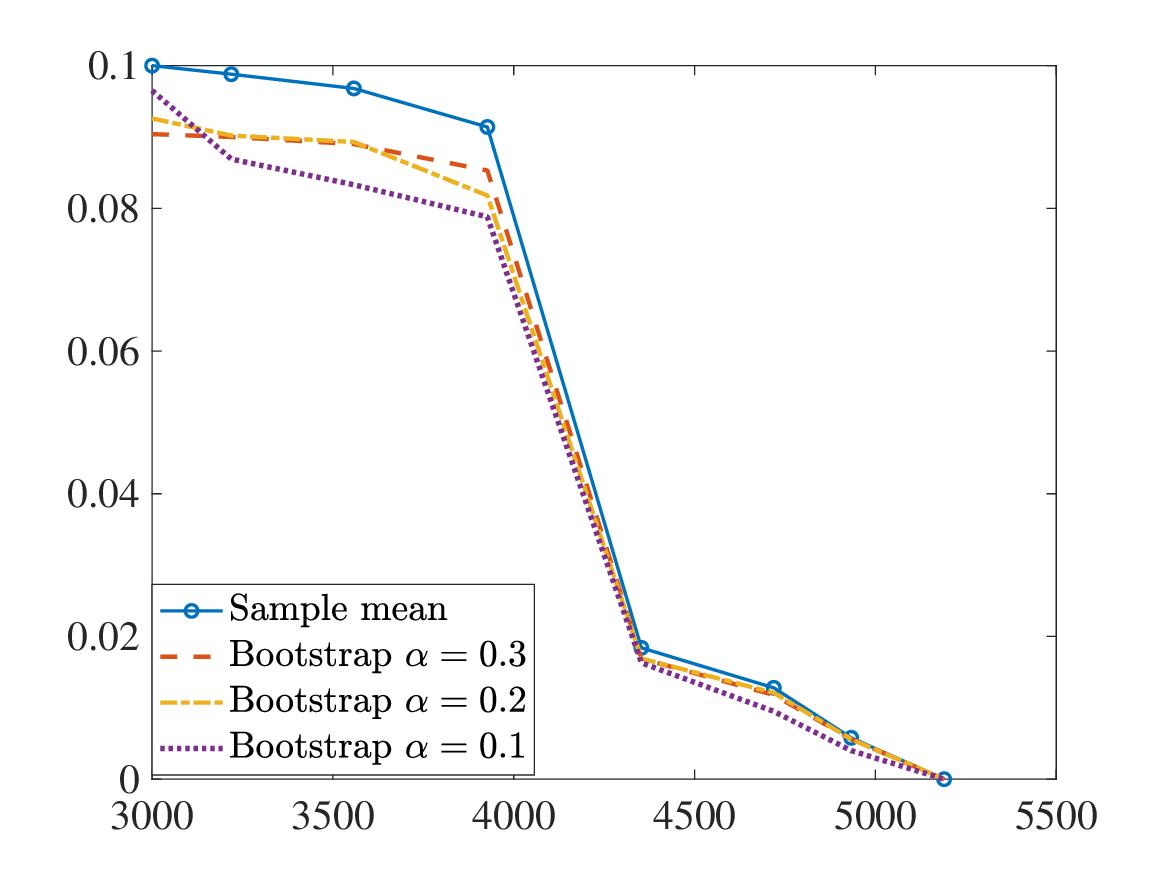}
  \textbf{\tiny (h) M$\&$R fee}
\endminipage
\caption{\small The single-attribute components of average consumer preference. These components include the disutility functions for {\em Price (a),  Acceleration (d), Depreciation Rate (g), and M$\&$R fee (h)}, and the utility functions for {\em Fuel Consumption (b), Wheelbase (c), Comfort (e), and Dealership (f).}}
\label{fig:oos-compar-2}
\vspace{-0.5cm}
\end{figure}

\textbf{Sample Generation:}
We next generate a random sample of $V$ using a logistic regression method and construct $\fP_B(\alpha)$. The logistic regression method is developed in \cite{feng_assessment_2016} to assess the consumer preference toward passenger vehicles in the US market. The dependent variable of the logistic regression model is the market share of the car industry in years 2013 and 2014, which is regarded as the probability distribution of consumers' choices among major vehicle brands. The independent variables are the performance of those brands regarding the attributes listed in Table~\ref{tab:project-list}. Following the approach in \cite{feng_assessment_2016}, we select a set of representative car brands and survey the monthly market shares and attribute values for each brand during two years. Let $S^\ell_j$ be the market share of brand $j$ at the $\ell$-th month, $x^j$ be the vector of the attribute values of brand $j$, and $u(\cdot; v^\ell)$, the utility function with the parameter vector $v^\ell$, represent the consumer preference at the $\ell$-th month. The probability that consumers purchase brand $j$ in month $\ell$ is 
	\begin{align} \label{mod:logistic}
		S^\ell_j = \frac{\exp  ( \eta u(x^j; v^\ell))}{\sum_{i=1}^J \exp( \eta u(x^i; v^\ell)) },
	\end{align}
	where $\eta$ is the parameter ensuring the normalization of $u$. 
 Estimating the parameters $\eta$ and $v^\ell$ in the logistic model \eqref{mod:logistic} using the survey data, we generate a random sample of $V$ including 24 month-wise observations, i.e., the estimate of $v^\ell$ for $\ell = 1, \dots, 24$. The case of ``Sample mean'' in Figure \ref{fig:oos-compar-2} draws each single-attribute component of $u(\cdot; \bar V)$ parameterized by the sample average $\bar V$.
Note that, with attributes such as {\em Price, Fuel Consumption, Acceleration, Depreciation Rate, and M$\&$R fee}, a lower utility value corresponds to a higher attribute value. To enhance readability, we present their disutility functions in Figure \ref{fig:oos-compar-2}. The utility functions, employed in \ref{DUPRO}, are the reflections of these disutility functions over the y-axis.
In addition, we interpret a population of customer's 
 preferences 
as 
a single  idiosyncratic 
 customer's random preference
 at different ``states'' whose average utility captures the one of the population.
 This kind of approach is used by
 \cite{liu2022stackelberg}
 in their recent work on Stackelberg risk preference design.

	      \begin{table}[h]
	    \centering
	    \begin{tiny}
	    \begin{tabular}{|c|c|c|c|c|c|c|c|c|c|c|}
	       \hline
	       \tabincell{c}{$ $} 
	       \tabincell{c}{} $\alpha$
	       &  \tabincell{c}{Safety \\promotion} 
	       &  \tabincell{c}{New car \\model \\develop-\\ment}   
	       &  \tabincell{c}{Engine \\upgrade} 
	       &  \tabincell{c}{E-platform \\develop-\\ment}
	       &  \tabincell{c}{CFD \\testing\\ system\\ development} 
	       &  \tabincell{c}{CMP \\devel-\\opment}   
	       &  \tabincell{c}{Checking\\ fixture \\ promotion} 
	       &  \tabincell{c}{NVH \\digitali-\\zation}
	       &  \tabincell{c}{Driving \\assistance\\  system\\ development} 
	       &  \tabincell{c}{Digitali-\\zation of\\ marketing \\network }
	       \\
	       \hline
	     \hline
	\tabincell{c}{$0.10$ \\ (Optimal solution $z$)}  &   $z_1$ = 0  &  $z_2$ = 0 &   $z_3$ = 0   & $z_4$ = 1
	                           &   $z_5$ = 1  &  $z_6$ = 1  &   $z_7$ = 1   & $z_8$ = 1
	                           &   $z_9$ = 0  &  $z_{10}$ = 0  
	   	\\
	     \hline
  \rule{0pt}{7pt}
$0.10$ 
	    &   
	    \multicolumn{10}{c|}{ Price~\$ $37.7$K, Fuel Consumption~$30$MPG, Wheelbase $107$in, 
	    Acceleration $6.5$ (0-60miles, sec)} \\
  \rule{0pt}{7pt}
(Attribute values $\hat{x}(z)$)	    & \multicolumn{10}{c|}{ Comfort rating $4.07$, Dealership $1250$, Depreciation rate $0.19$, M\&R fee \$ $4.9$K }  
	     \\
	     \hline
	  \tabincell{c}{$0.10$\\ (Optimal value)}
     & \multicolumn{10}{c|}{$0.5162$ }
	       \\
	     \hline
 
    \tabincell{c}{$0.20$\\ (Optimal solution $z$)} &   $z_1$ = 0  &  $z_2$ = 0 &   $z_3$ = 1   & $z_4$ = 0
	                           &   $z_5$ = 1  &  $z_6$ = 1 &   $z_7$ = 0   & $z_8$ = 1
	                           &   $z_9$ = 0  &  $z_{10}$ = 0 
	     \\
	     \hline
	\rule{0pt}{7pt}
$0.20$ 
	    &   
	    \multicolumn{10}{c|}{ Price~\$ $37.7$K, Fuel Consumption~$32$MPG, Wheelbase $107$in, 
	    Acceleration $5.5$ (0-60miles, sec)} \\
  \rule{0pt}{7pt}
(Attribute values $\hat{x}(z)$)	    & \multicolumn{10}{c|}{ Comfort rating $3.92$, Dealership $1200$, Depreciation rate $0.14$, M\&R fee \$ $4.8$K }

	     \\
	  \hline 
    \tabincell{c}{$0.20$\\ (Optimal value)}  & \multicolumn{10}{c|}{ $0.5381$ }
	       \\
	     \hline
	\tabincell{c}{$0.30$\\ (Optimal solution $z$)}  &   $z_1$ = 0  &  $z_2$ = 1&   $z_3$ = 0   & $z_4$ = 1
	                           &   $z_5$ = 1  &  $z_6$ = 1  &   $z_7$ = 0   & $z_8$ = 0
	                           &   $z_9$ = 0  &  $z_{10}$ = 0  
	     \\
	     \hline
  \rule{0pt}{7pt}
$0.30$ 
	    &   
	    \multicolumn{10}{c|}{ Price~\$ $35.2$K, Fuel Consumption~$36$MPG, Wheelbase $105$in, 
	    Acceleration $5.5$ (0-60miles, sec)} \\
  \rule{0pt}{7pt}
(Attribute values $\hat{x}(z)$)	    & \multicolumn{10}{c|}{ Comfort rating $3.98$, Dealership $1400$, Depreciation rate $0.19$, M\&R fee \$ $5.5$K }
	     \\
	 \hline
   \tabincell{c}{$0.30$\\ (Optimal value)}      & \multicolumn{10}{c|}{ $0.5507$ }
	       \\	       
	       \hline
	\tabincell{c}{SAA \\ (Optimal solution $z$)}  & $z_1$ = 0  &  $z_2$ = 0 &   $z_3$ = 0   & $z_4$ = 0
	                           &   $z_5$ = 1  &  $z_6$ = 1  &   $z_7$ = 0   & $z_8$ = 1
	                           &   $z_9$ = 0  &  $z_{10}$ = 1  
	     \\
	     \hline
  \rule{0pt}{7pt}
SAA  
	    &   
	    \multicolumn{10}{c|}{  Price~\$ $37.7$K, Fuel Consumption~$30$MPG, Wheelbase $107$in, 
	    Acceleration $6.5$ (0-60miles, sec)} \\
  \rule{0pt}{7pt}
(Attribute values $\hat{x}(z)$)	    & \multicolumn{10}{c|}{ Comfort rating $3.95$, Dealership $1200$, Depreciation rate $0.14$, M\&R fee \$ $4.5$K }
	     \\
	 \hline
   \tabincell{c}{SAA \\ (Optimal value)}      & \multicolumn{10}{c|}{ $0.5669$ }
	       \\	       
	       \hline
	     \rule{0pt}{7pt}
	    &   
	    \multicolumn{10}{c|}{ Price~\$ $38.0$K, Fuel Consumption~$30$MPG, Wheelbase $110$in, 
	    Acceleration $8$ (0-60miles, sec)} \\
  \rule{0pt}{7pt}
$x^0$	    & \multicolumn{10}{c|}{ Comfort rating $3.80$, Dealership $1050$, Depreciation rate $0.25$, M\&R fee \$ $5.5$K }
	     \\
	     \hline
	    \end{tabular}
	    \end{tiny}
        \vspace{0.2cm}
	    \caption{The optimal solutions and values of model \eqref{prob:project-car-max}  with different $\alpha$.}
	    \label{tab:opt-alpha-limit}
	\end{table}

\textbf{Computational Results:}
In this study,
we choose the size of bootstrap samples 
$K = 1000$ and budget limit $\Phi = \$200$ million. The value of $\alpha$ is varied from $0.10$ to $0.30$ in the setting of $\fP_B(\alpha)$. Note that the size of $\fP_B(\alpha)$ shrinks with the increase of $\alpha$ and 
subsequently the level of robustness dwindles. 
Table~\ref{tab:opt-alpha-limit} presents the results of model \eqref{prob:project-car-max} showing the best investment portfolio and the optimal value of consumer preference, and Figure \ref{fig:oos-compar-2} displays the worst-case utility functions. For comparative purposes, the last row in Table~\ref{tab:opt-alpha-limit} gives the baseline attribute values $x^0$ before the investment. 
When $\alpha$ is small, the ambiguity set is large. 
The robust approach with a high level of confidence
requires the DM to take a more conservative action
than the bottom line of the majority of customers
would prefer. Note that the SAA case represents the maximization of the observed sample mean of the random utility function, a method that does not incorporate considerations for distributional ambiguity. In the SAA case, the proposed strategic recommendations include the implementation of {\em ``CFD testing system development''}, {\em ``CMP development''}, {\em ``NVH digitalization''}, and {\em ``Digitalization of marketing network''}. These initiatives are aimed at enhancing attributes such as {\em ``Acceleration''}, {\em ``Comfort''}, {\em ``Dealership''}, {\em `` Depreciation rate''},  and {\em ``M$\&$R fee''}. It should be acknowledged, however, that these recommendations are made with the observation that there could be a potential trade-off, notably a sacrifice in the {\em ``Wheelbase''} attribute. In practical terms, the {\em ``Wheelbase''} attribute often receives 
fewer attentions during the vehicle purchasing process by customers. When the value of $\alpha$ is set to 0.30, there is a noticeable shift towards investments that aim to improve ``{\em Price}'' competitiveness and ``{\em Fuel consumption}'' through the ``New car model development'' project. Historical data shows that retail gasoline prices peaked in both current and constant dollars during the years 2013 and 2014, making fuel efficiency a priority at that time \cite{DOE_Fact915_2016}. This trend also supports the observed investments in enhancing acceleration and dealership services, which contribute to the overall optimal value of the expected car utility. As $\alpha$ is further decreased to 0.20, a more conservative strategy emerges, suggesting the replacement of ``{\em New car model development}'' with ``{\em Engine update}''. Although this approach yields a smaller improvement in ``{\em Fuel consumption}'', it simultaneously reduces the risks and costs associated with the development of a new car model. An even more conservative stance is taken when $\alpha = 0.10$, 
where the ``{\em Comfort rating}'' attribute gains prominence. This attribute enjoys stable demand from customers and is less susceptible to random economic fluctuations.


\section{Conclusions}
\label{sec:Conclusions}
The past decade has witnessed a surge of research interest in PRO which handles 
decision-making problems under ambiguity 
of the DM's utility preferences. The traditional utility theory, which existing PRO models rely on, 
assumes preference representations to be deterministic and consistent. 
Here we propose a novel PRO model which, combining the stochastic utility theory and distributionally robust optimization techniques, 
is capable of dealing with decision-making problems with inconsistent and mutable utility preferences.

We concentrate on the case that the samples of the 
random utility function are difficult to obtain and
use the well-known bootstrap method
for constructing the ambiguity with relatively small sample size. 
We reformulate the resulting DUPRO model 
as an
MILP  when the random utility function is piecewise linear 
and represented 
in the specific increment-based form, 
and LP when it is concave. 
We conduct some numerical tests to analyze the effects of the crucial parameters in the bootstrap method such as the size of original 
sample data, the number of the bootstrap resamples, the critical value of the confidence region, 
and the variance of the underlying true distribution. 
The preliminary results show that the proposed 
model and computational scheme work very well.
Of course, there is a gap between bootstrap confidence region (based on resamples) 
and the confidence region based on the original 
samples. In order for the two confidence regions to be close, the size of the original samples should be sufficiently large, see \cite[Theorem 1]{yeh_balanced_1997}.

To extend the scope of applicability of the proposed model and computational schemes,
we consider the PRO model with general random utility functions and  
discuss approximation of general random utility functions by piecewise linear random utility functions. 
Specifically, we quantify propagation of the error of the approximation to the optimal value and optimal solutions.
In the case when the data of the piecewise linear random utility functions are obtained from samples,
we demonstrate 
convergence of the optimal value obtained from solving the
sample-based piecewise linear approximated PRO model to its true counterpart as the sample size increases.

While this work effectively 
extends the existing research of PRO
from deterministic utility to random utility, it
also advances the literature of random utility theory
by
providing new modelling paradigms, computational schemes and underlying theory 
if we interpret a single DM's random preferences as 
a population of customer preferences. 
We leave it for future research as to how we may 
use them to study discrete-choice demand models.

	\section{Declaration of competing interest}
The authors declare that they have no known competing financial interests or personal relationships that could have appeared to influence the work reported in this paper.	\\
	

\noindent \textbf{Acknowledgments:} The authors thank two anonymous referees, whose comments have resulted in a significant improvement to the original draft.

	


	
	

	\vspace{1mm}








\begin{appendices}
\newpage 
\section{Data of Project Investment Problem -- Section \ref{subsec:Project Investment}.}
\label{App:table-project-case-2}

\begin{wraptable}{H}{0.5\textwidth}
\vspace{-8mm}
	\begin{sideways}
		\begin{tiny}
			\begin{tabular}{|l|c|c|c|c|c|c|c|c|c|c|c|}
				\hline
				\tabincell{c}{} 
				&  \tabincell{c}{Base \\model} 
				&  \tabincell{c}{Safety \\promotion} 
				&  \tabincell{c}{New car \\model \\develop-\\ment}   
				&  \tabincell{c}{Engine \\upgrade} 
				&  \tabincell{c}{E-platform \\develop-\\ment}
				&  \tabincell{c}{CFD testing\\ system\\ development} 
				&  \tabincell{c}{CMP \\ development}   
				&  \tabincell{c}{Checking\\ fixture \\ promotion} 
				&  \tabincell{c}{Noise, \\vibration, \\  harshness\\ (NVH) \\digitali-\\zation}
				&  \tabincell{c}{Driving \\assistance\\  system\\ development} 
				&  \tabincell{c}{Digitali-\\zation \\of marketing \\network }
				\\
				& ($x^0$) &  ($y^1$)  &  ($y^2$)  &   ($y^3$)   & ($y^4$)
				&   ($y^5$)  &  ($y^6$)  &   ($y^7$)  & ($y^8$)
				&   ($y^9$)  &  ($y^{10}$) 
				\\
				\hline
				\hline
				\tabincell{c}{Price\\ (\$k)} & 38 &   +7  &  -1  &   0   & 0
				&   0  &  -1.8  &   0   & +1.5
				&   +1.5  &  0 
				\\
				\hline
				\tabincell{c}{Fuel con-\\sumption\\ (MPG)} & 30 &   0     &  +6     &   +2      & 0
				&   0     &  0     &   0      & 0
				&   +2     &  0     
				\\
				\hline
				\tabincell{c}{Wheel-\\base (in)} & 110 &   0     &  -2     &   0      & 0
				&   0     &  -3     &   0      & 0
				&   0     &  0     
				\\
				\hline
				\tabincell{c}{Accelera-\\tion (0-60 \\miles, sec)} & 8 &  0 & -1 & -1 & 0	
				& -1.5 & 0 & 0 & 0 & 0 & 0 
				\\
				\hline
				\tabincell{c}{Comfort rating \\(0-5)} & 3.8 &   +0.02     &  +0.04     &   0      &  +0.1
				&  +0.04   & 0 & +0.05 & +0.08 & +0.02 & +0.03   
				\\
				\hline
				\tabincell{c}{Dealership \\(\# of dealers )} & 1050 &   0 &  +150 & +150 &  +200
				& 0 & 0 & 0 & 0 & 0 & +150   
				\\
				\hline
				\tabincell{c}{Depreciation \\ rate (0-1)} & 0.25 & -0.08 & -0.03 & -0.05 & 0 & 0
				& -0.03	& 0	& -0.03	& -0.02	& -0.05	
				\\
				\hline
				\tabincell{c}{M\&R \\ fee (\$K)} & 5.5 & -0.5	& -0.5	& +0.3 & +0.5	
				& -0.5	& 0	& -0.1	& -0.5	& -0.2	& 0	
				\\
				\hline
				\hline
				\tabincell{c}{Project \\ costs\\ (\$ million)} & --- & 50	& 100 & 70	& 50 & 20	
				& 30 & 20 & 80 & 50 &	70\\
				\hline
			\end{tabular}
		\end{tiny}
   \end{sideways}
    \caption{The list of candidate projects.}
	\label{tab:project-list}
	\vspace{-20mm}
  \end{wraptable}

\end{appendices}

\end{document}